\documentstyle[
11pt,amsfonts]{article}

\newcommand {\rel} {{\mathbb R}}

\newcommand {\nat} {{\mathbb N}}

\newtheorem{proposition}{Proposition}[section]
\newtheorem{theorem}{Theorem}[section]
\newtheorem{corollaryth}[theorem]{Corollary}
\newtheorem{lemma}[proposition]{Lemma}

\newtheorem{definition}[proposition]{Definition}

\renewcommand{\theequation}{\mbox{\arabic{section}.\arabic{equation}}}

\newcommand{\nc}{\newcommand}

\newcommand{\bcite}[1] {\cite{#1}}

\newcommand {\pr} {\bf}

\newcommand{\proof} {   \begin{flushright}
                        ///
                        \end{flushright}
                }

\newcommand{\defin} { \hspace*{\fill} $\Box$ }

\newcommand{\mint}[1]{\mbox{$\displaystyle{
-\hspace{-1.05em}\int_{#1}} $}}

\def    \mean   {{ \bf \vec H }}

\def    \Ho	{{ {\cal H}^1 }}

\def    \Lt	{{ {\cal L}^2 }}

\def	\d	{{ \ {\rm d} }}

\def	\pro	{{ \pi }}
\def	\funda	{{ \nu }}

\def    \W      {{ {\mathcal W} }}
\def	\mini	{{ {\mathcal M}_{p,n} }}

\def	\sur	{{ \Sigma }}
\def	\fu	{{ \sigma }}
\def	\poin	{{ poin }}

\def	\geu	{{ g_{euc} }}
\def	\gp	{{ g_{\poin} }}
\def	\tgp	{{ \tilde g_{\poin} }}
\nc{\gpo}[1]	{{ g_{\poin}^{#1} }}
\nc{\gpu}[1]	{{ g_{\poin,{#1}} }}
\def	\gpm	{{ g_{\poin,m} }}
\nc{\gpmo}[1]	{{ g_{\poin,m}^{#1} }}
\nc{\tgpmo}[1]	{{ \tilde g_{\poin,m}^{#1} }}

\nc{\gplo}[1]	{{ g_{\poin,\lambda}^{#1} }}

\def	\tgpm	{{ \tilde g_{\poin,m} }}

\def	\tgpml 	{{ \tilde g_{\poin,m,\lambda} }}
\def	\tgpmlm	{{ \tilde g_{\poin,m,\lambda,\mu} }}

\def	\vf	{{ \mathcal{X} }}
\def	\moeb	{{ \Phi }}
\def	\diffeo	{{ \phi }}

\def	\vari	{{ V }}
\def	\tp	{{ \pi }}
\def	\tpc	{{ \hat \tp }}

\def	\teich	{{ \cal T }}
\def	\M	{{ M }}
\def	\pos	{{ \gamma }}

\def	\nor	{{ \cal N }}

\def	\mark	{{ * }}

\def    \rechtsarrow {{ \searrow }}
\def    \linksarrow {{ \nearrow }}

\nc{\doub}[1]{{ \ddot{#1} }}
\nc{\dd}{ \begin{displaymath} }
\nc{\df}{ \end{displaymath} }
\nc{\dcd}{ \begin{displaymath} \begin{array}{c}}
\nc{\dcf}{ \end{array} \end{displaymath} }
\nc{\ee}{ \begin{equation} }
\nc{\ef}{ \end{equation} }
\nc{\ad}{ \begin{array}{c} }
\nc{\af}{ \end{array} }

\begin{document}

\begin{center}
{\huge \bf Minimizers of the Willmore functional \\
\vspace{.15cm}
under fixed conformal class}
\\ \ \\
Ernst Kuwert
\footnote{E.\,Kuwert and R.\,Sch\"atzle
were supported by the DFG Forschergruppe 469.
Main parts of this article were discussed
during a stay of both authors
at the Centro Ennio De~Giorgi in Pisa.} \\
Mathematisches Institut der
Albert-Ludwigs-Universit\"at Freiburg,
Eckerstra\ss e 1, \\
D-79104 Freiburg,
Germany,
email: kuwert@mathematik.uni-freiburg.de \\
\ \\
Reiner Sch\"atzle$^{\ 1}$ \\
Mathematisches Institut der
Eberhard-Karls-Universit\"at T\"ubingen, \\
Auf der Morgenstelle 10,
D-72076 T\"ubingen, Germany, \\
email: schaetz@everest.mathematik.uni-tuebingen.de \\
 
\end{center}
\vspace{1cm}

\begin{quote}

{\bf Abstract:}
We prove the existence of a smooth minimizer
of the Willmore energy in the class of
conformal immersions of a given closed Riemann surface
into $\ \rel^n, n = 3,4,\ $
if there is one conformal immersion with Willmore energy
smaller than a certain bound $\ \W_{n,p}\ $
depending on codimension and genus $\ p\ $ of the Riemann surface.
For tori in codimension 1, we know $\ \W_{3,1} = 8 \pi\ $.
\end{quote}

\vspace{1cm}



\setcounter{equation}{0}

\section{Introduction} \label{intro}

For an immersed closed surface $\ f : \sur
\rightarrow \rel^n\ $ the Willmore functional
is defined by
\begin{displaymath}
	\W(f) = \frac{1}{4} \int \limits_\sur |\mean|^2 \d \mu_g,
\end{displaymath}
where $\ \mean\ $ denotes the mean curvature vector of
$\ f, g = f^* \geu\ $ the pull-back metric
and $\ \mu_g\ $ the induced
area measure on $\ \sur\ $.
The Gau\ss\ equations and the Gau\ss-Bonnet Theorem
give rise to equivalent expressions
\begin{equation} \label{intro.gauss}
	\W(f) = \frac{1}{4} \int \limits_\sur
	|A|^2 \d \mu_g + 2 \pi (1 - p(\sur))
	= \frac{1}{2} \int \limits_\sur
	|A^\circ|^2 \d \mu_g + 4 \pi (1 - p(\sur))
\end{equation}
where $\ A\ $ denotes the second fundamental form,
$\ A^\circ = A - \frac{1}{2} g \otimes H\ $ its
tracefree part and $\ p(\sur)\ $ is the genus of $\ \sur\ $.

We always have $\ \W(f) \geq 4 \pi\ $ with equality
only for round spheres, see \bcite{will}
in codimension one that is $\ n = 3\ $.
On the other hand, if $\ \W(f) < 8 \pi\ $
then $\ f\ $ is an embedding by an inequality of Li and Yau
in \bcite{li.yau}.

Critical points of $\ \W\ $ are called Willmore surfaces
or more precisely Willmore immersions.
They satisfy the Euler-Lagrange equation
which is the fourth order, quasilinear geometric equation
\begin{displaymath}
        \Delta_g \mean + Q(A^0) \mean = 0
\end{displaymath}
where the Laplacian of the normal bundle along $\ f\ $ is used
and $\ Q(A^0)\ $ acts linearly on normal vectors
along $\ f\ $ by
\begin{displaymath}
        Q(A^0) \phi := g^{ik} g^{jl} A^0_{ij}
        \langle A^0_{kl} , \phi \rangle.
\end{displaymath}
The Willmore functional is scale invariant
and moreover invariant under the full M\"obius group of $\ \rel^n\ $.
As the M\"obius group is non-compact,
minimizers of the Willmore energy cannot
be found via the direct method.

In \bcite{kuw.schae.will6}, we investigated the relation
of the pull-back metric $\ g\ $ to constant curvature
metrics on $\ \sur\ $ after dividing out the M\"obius group.
More precisely, we proved that for immersions
$\ f: \sur \rightarrow \rel^n, n = 3,4,\ $
and genus $\ p = p(\sur) \geq 1\ $
satisfying $\ \W(f) \leq \W_{n,p} - \delta
\mbox{ for some } \delta > 0\ $, where
\begin{equation} \label{intro.moeb.energ}
\begin{array}{c}
	\W_{3,p} = \min(8 \pi ,
	4 \pi + \sum_k (\beta_{p_k}^3 - 4 \pi)\ |
	\ \sum_k p_k = p, 0 \leq p_k < p), \\
	\W_{4,p} = \min(8 \pi , \beta_p^4 + 8 \pi / 3 ,
	4 \pi + \sum_k (\beta_{p_k}^4 - 4 \pi)\ |
	\ \sum_k p_k = p, 0 \leq p_k < p),
\end{array}
\end{equation}
and
\begin{equation} \label{intro.mini}
	\beta_p^n := \inf \{ \W(f)\ |
	\ f: \sur \rightarrow \rel^n \mbox{ immersion},
	p(\sur) = p\ \},
\end{equation}
there exists a M\"obius transformation $\ \moeb\ $
of the ambient space $\ \rel^n\ $
such that the pull-back metric $\ \tilde g = (\moeb \circ f)^* \geu\ $
differs from a constant curvature metric $\ e^{-2u} \tilde g\ $
by a bounded conformal factor, more precisely
\begin{displaymath}
        \parallel u \parallel_{L^\infty(\sur)},
        \parallel \nabla u
	\parallel_{L^2(\sur,\tilde g)}
	\leq C(p,\delta).
\end{displaymath}
In this paper, we consider conformal immersions
$\ f: \sur \rightarrow \rel^n\ $
of a fixed closed Riemann surface $\ \sur\ $
and prove existence of minimizers in this
conformal class under the above energy assumptions.
\\ \ \\
{\bf Corollary \ref{euler.conf}}
{\it
Let $\ \sur\ $ be a closed Riemann surface
of genus $\ p \geq 1\ $
with
\begin{displaymath}
	\inf \{ \W(f)\ |
	\ f: \sur \rightarrow \rel^n
	\mbox{ conformal immersion}\ \}
	< \W_{n,p},
\end{displaymath}
where $\ \W_{n,p}\ $
is defined in (\ref{intro.moeb.energ}) above
and $\ n = 3,4\ $.

Then there exists a smooth conformal immersion
$\ f: \sur \rightarrow \rel^n\ $
which minimizes the Willmore energy
in the set of all conformal immersions.
Moreover $\ f\ $ satisfies
the Euler-Lagrange equation
\begin{displaymath}
	\Delta_g \mean + Q(A^0) \mean = g^{ik} g^{jl} A_{ij}^0 q_{kl}
	\quad \mbox{on } \sur,
\end{displaymath}
where $\ q\ $ is a smooth transverse traceless
symmetric 2-covariant tensor with respect
to $\ g = f^* \geu\ $.
}
\defin
\\ \ \\
{\large \bf Acknowledgement:} \\
The main ideas of this paper came out
during a stay of both authors
at the Centro Ennio De~Giorgi in Pisa.
Both authors thank very much
the Centro Ennio De~Giorgi in Pisa
for the hospitality and for providing
a fruitful scientific atmosphere.


\setcounter{equation}{0}

\section{Direct method} \label{direct}

Let $\ \sur\ $ be a closed orientable surface
of genus $\ p \geq 1\ $ with smooth metric $\ g\ $
satisfying
\begin{displaymath}
	\W(\sur,g,n) :=
	\inf \{ \W(f)\ |
	\ f: \sur \rightarrow \rel^n
	\mbox{ conformal immersion to } g\ \}
	< \W_{n,p},
\end{displaymath}
as in the situation of Corollary \ref{euler.conf}
where $\ \W_{n,p}\ $
is defined in (\ref{intro.moeb.energ}) above
and $\ n = 3,4\ $.
To get a minimizer, we consider a minimizing sequence
of immersions $\ f_m: \sur \rightarrow \rel^n\ $
conformal to $\ g\ $ in the sense
\begin{displaymath}
	\W(f_m) \rightarrow
	\inf \{ \W(f)\ |
	\ f: \sur \rightarrow \rel^n
	\mbox{ conformal immersion to } g\ \}.
\end{displaymath}
After applying suitable M\"obius transformations
according to \bcite{kuw.schae.will6} Theorem 4.1
we will be able to estimate $\ f_m \mbox{ in }
W^{2,2}(\sur) \mbox{ and } W^{1,\infty}(\sur)\ $,
see below, and after passing to a subsequence, we get a limit
$\ f \in W^{2,2}(\sur) \cap W^{1,\infty}(\sur)\ $
which is an immersion in a weak sense,
see (\ref{direct.prop.imm}) below.
To prove that $\ f\ $ is smooth,
which implies that it is a minimizer,
and that $\ f\ $ satisfies
the Euler-Lagrange equation in Corollary \ref{euler.conf}
we will consider variations,
say of the form $\ f + \vari\ $.
In general, these are not conformal to $\ g\ $ anymore,
and we want to correct it by $\ f + \vari + \lambda_r \vari_r\ $
for suitable selected variations $\ \vari_r\ $.
Now even these are not conformal to $\ g\ $
since the set of conformal metrics is quite small
in the set of all metrics.
To increase the set of admissible pull-back metrics,
we observe that it suffices for
$\ (f + \vari + \lambda_r \vari_r) \circ \diffeo\ $
being conformal to $\ g\ $ for some diffeomorphism
$\ \diffeo \mbox{ of } \sur\ $.
In other words, the pullback metric
$\ (f + \vari + \lambda_r \vari_r)^* \geu\ $
need not be conformal to $\ g\ $,
but has to coincide only in the modul space.
Actually, we will consider the Teichm\"uller space,
which is coarser than the modul space,
but is instead a smooth open manifold,
and the bundle projection $\ \tp: {\cal M} \rightarrow \teich\ $
of the sets of metrics $\ {\cal M}\ $
into the Teichm\"uller space $\ \teich\ $,
see \bcite{fisch.trom.conf}, \bcite{trom.teich}.
Clearly $\ \W(\sur,g,n)\ $ depends only
on the conformal structure defined by $\ g\ $,
in particular descends to Teichm\"uller space
and leads to the following definition.

\begin{definition} \label{direct.def}

We define $\ \mini: \teich \rightarrow [0,\infty]
\mbox{ for } p \geq 1, n \geq 3,\ $
by selecting a closed, orientable
surface $\ \sur \mbox{ of genus } p\ $ and
\begin{displaymath}
	\mini(\tau)
	:= \inf \{ \W(f)\ |
	\ f: \sur \rightarrow \rel^n
	\mbox{ smooth immersion},
	\tp(f^* \geu) = \tau\ \}.
\end{displaymath}
\defin
\end{definition}
We see
\begin{displaymath}
	\mini(\tp(g)) = \W(\sur,g,n).
\end{displaymath}
Next, $\ \inf_\tau \mini(\tau) = \beta_p^n\ $
for the infimum under fixed genus defined in (\ref{intro.mini}),
and, as the minimum is attained and $\ 4 \pi < \beta_p^n < 8 \pi\ $,
see \bcite{sim.will} and \bcite{bau.kuw},
\begin{displaymath}
	4 \pi < \min \limits_{\tau \in \teich} \mini(\tau)
	= \beta_p^n < 8 \pi.
\end{displaymath}
In the following proposition, we consider a slightly
more general situation than above.

\begin{proposition} \label{direct.prop}

Let $\ f_m: \sur \rightarrow \rel^n, n = 3,4,\ $ be smooth immersions
of a closed, orientable surface $\ \sur \mbox{ of genus } p \geq 1\ $
satisfying
\begin{equation} \label{direct.prop.ener}
	\W(f_m) \leq \W_{n,p} - \delta
\end{equation}
and
\begin{equation} \label{direct.prop.teich}
	\tp(f_m^* \geu) \rightarrow \tau_0 \mbox{ in } \teich.
\end{equation}
Then replacing $\ f_m \mbox{ by }
\moeb_m \circ f_m \circ \diffeo_m\ $
for suitable M\"obius transformations $\ \moeb_m\ $
and diffeomorphisms $\ \diffeo_m \mbox{ of } \sur\ $
homotopic to the identity, we get
\begin{equation} \label{direct.prop.w22}
	\limsup \limits_{m \rightarrow \infty}
	\parallel f_m \parallel_{W^{2,2}(\sur)}
	\leq C(p,\delta,\tau_0)
\end{equation}
and $\ f_m^* \geu = e^{2 u_m} \gpm\ $
for some unit volume constant curvature metrics $\ \gpm\ $ with
\begin{equation} \label{direct.prop.met}
\begin{array}{c}
	\parallel u_m \parallel_{L^\infty(\sur)},
	\parallel \nabla u_m \parallel_{L^2(\sur,\gpm)}
	\leq C(p,\delta), \\
	\gpm \rightarrow \gp
	\quad \mbox{smoothly}
\end{array}
\end{equation}
with $\ \tp(\gp) = \tau_0\ $.
After passing to a subsequence
\begin{equation} \label{direct.prop.conv}
\begin{array}{c}
	f_m \rightarrow f \mbox{ weakly in } W^{2,2}(\sur),
	\mbox{weakly}^* \mbox{ in } W^{1,\infty}(\sur), \\
	u_m \rightarrow u \mbox{ weakly in } W^{1,2}(\sur),
	\mbox{weakly}^* \mbox{ in } L^\infty(\sur),
\end{array}
\end{equation}
and
\begin{equation} \label{direct.prop.imm}
	f^* \geu = e^{2u} \gp
\end{equation}
where $\ (f^* \geu)(X,Y) := \langle \partial_X f , \partial_Y f \rangle
\mbox{ for } X,Y \in T \sur\ $.
\end{proposition}
{\pr Proof:} \\
Clearly, replacing $\ f_m \mbox{ by }
\moeb_m \circ f_m \circ \diffeo_m\ $ as above
does neither change the Willmore energy
nor the projection into the Teichm\"uller space.

By \bcite{kuw.schae.will6} Theorem 4.1
after applying suitable M\"obius transformations,
the pull-back metric $\ g_m := f_m^* \geu\ $
is conformal to a unique constant curvature metric
$\ e^{-2 u_m} g_m =: \gpm\ $ of unit volume
with
\begin{displaymath}
	osc_\sur u_m,
	\parallel \nabla u_m \parallel_{L^2(\sur,g_m)}
	\leq C(p,\delta).
\end{displaymath}
Combining the M\"obius transformations with suitable
homotheties, we may further assume
that $\ f_m\ $ has unit volume.
This yields
\begin{displaymath}
	1 = \int \limits_\sur \d \mu_{g_m}
	= \int \limits_\sur e^{2u_m} \d \mu_\gpm,
\end{displaymath}
and, as $\ \gpm\ $ has unit volume as well,
we conclude that $\ u_m\ $ has a zero on $\ \sur\ $,
hence
\begin{displaymath}
	\parallel u_m \parallel_{L^\infty(\sur)},
	\parallel \nabla u_m \parallel_{L^2(\sur,\gpm)}
	\leq C(p,\delta)
\end{displaymath}
and, as $\ f_m^* \geu = g_m = e^{2 u_m} \gpm\ $,
\begin{displaymath}
	\parallel \nabla f_m \parallel_{L^\infty(\sur,\gpm)}
	\leq C(p,\delta).
\end{displaymath}
Next
\begin{displaymath}
	\Delta_\gpm f_m
	= e^{2 u_m} \Delta_{g_m} f_m
	= e^{2 u_m} \mean_{f_m}
\end{displaymath}
and
\begin{displaymath}
	\parallel \Delta_\gpm f_m \parallel_{L^2(\sur,\gpm)}^2
	\leq e^{2 \max u_m} \int \limits_\sur |\mean_{f_m}|^2
	e^{2 u_m} \d \mu_\gpm =
\end{displaymath}
\begin{equation} \label{direct.prop.lap}
	= 4 e^{2 \max u_m} \W(f_m) \leq C(p,\delta).
\end{equation}
To get further estimates, we employ the convergence
in Teichm\"uller space (\ref{direct.prop.teich}).
We consider a slice $\ {\cal S}\ $ of unit volume
constant curvature metrics for $\ \tau_0 \in \teich\ $,
see \bcite{fisch.trom.conf}, \bcite{trom.teich}.
There exist unique $\ \tgpm \in {\cal S}
\mbox{ with } \tp(\tgpm) = \tp(g_m) \rightarrow \tau_0\ $
for $\ m\ $ large enough,
hence
\begin{displaymath}
	\diffeo_m^* \tgpm = \gpm
\end{displaymath}
for suitable diffeomorphisms $\ \diffeo_m
\mbox{ of } \sur\ $ homotopic to the identity.
Replacing $\ f_m \mbox{ by } f_m \circ \diffeo_m\ $,
we get $\ \gpm = \tgpm \in {\cal S}\ $ and
\begin{displaymath}
	\gpm \rightarrow \gp
	\quad \mbox{smoothly}
\end{displaymath}
with $\ \gp \in {\cal S}, \tp(\gp) = \tau_0\ $.
Then translating $\ f_m\ $ suitably, we obtain
\begin{displaymath}
	\parallel f_m \parallel_{L^\infty(\sur)}
	\leq C(p,\delta).
\end{displaymath}
Moreover standard elliptic theory, see \bcite{gil.tru} Theorem 8.8,
implies by (\ref{direct.prop.lap}),
\begin{displaymath}
	\parallel f_m \parallel_{W^{2,2}(\sur)}
	\leq C(p,\delta,\gp)
\end{displaymath}
for $\ m\ $ large enough.

Selecting a subsequence, we get $\ f_m \rightarrow f
\mbox{ weakly in } W^{2,2}(\sur),
f \in W^{1,\infty}(\sur),
D f_m \rightarrow D f\ $ pointwise almost everywhere,
$\ u_m \rightarrow u \mbox{ weakly in } W^{1,2}(\sur),\ $
pointwise almost everywhere, and
$\ u \in L^\infty(\sur)\ $.
Putting $\ g := f^* \geu\ $,
that is $\ g(X,Y)
:= \langle \partial_X f , \partial_Y f \rangle\ $,
we see
\begin{displaymath}
	g(X,Y) \leftarrow \langle \partial_X f_m , \partial_Y f_m \rangle
	= g_m(X,Y) = e^{2u_m} \gpm(X,Y)
	\rightarrow e^{2u} \gp(X,Y)
\end{displaymath}
for $\ X,Y \in T \sur\ $
and pointwise almost everywhere on $\ \sur\ $,
hence
\begin{displaymath}
	f^* \geu = g = e^{2u} \gp.
\end{displaymath}
\proof
We call a mapping $\ f \in W^{1,\infty}(\sur,\rel^n)
\mbox{ with } g := f^* \geu = e^{2u} \gp, u \in L^\infty(\sur)\ $
as in (\ref{direct.prop.imm}) a weak local bilipschitz immersion.
If further $\ f \in W^{2,2}(\sur,\rel^n)\ $,
we see $\ g, u \in W^{1,2}(\sur)\ $,
hence we can define weak Christoffel symbols $\ \Gamma \in L^2\ $
in local charts,
a weak second fundamental form $\ A\ $
and a weak Riemann tensor $\ R\ $
via the equations of Weingarten and Gau\ss
\begin{displaymath}
\begin{array}{c}
	\partial_{ij} f = \Gamma_{ij}^k \partial_k f + A_{ij}, \\
	R_{ijkl} = \langle A_{ik} , A_{jl} \rangle
	- \langle A_{ij} , A_{kl} \rangle.
\end{array}
\end{displaymath}
Of course this defines $\ \mean, A^0, K\ $ for $\ f\ $ as well.
Moreover we define the tangential and normal projections
for $\ \vari \in W^{1,2}(\sur,\rel^n) \cap L^\infty(\sur,\rel^n)\ $
\begin{equation} \label{direct.pro}
\begin{array}{c}
	\pro_f . \vari := g^{ij} \langle \partial_i f ,
	\vari \rangle \partial_j f, \quad
	\pro_f^\perp . \vari := \vari - \pro_\sur . \vari
	\in W^{1,2}(\sur,\rel^n) \cap  L^\infty(\sur,\rel^n).
\end{array}
\end{equation}
Mollifing $\ f\ $ as in \bcite{schoe.uhl.harmb} \S 4 Proposition,
we get smooth $\ f_m: \sur \rightarrow \rel^n\ $ with
\begin{equation} \label{direct.imm.strong}
	f_m \rightarrow f
	\quad \mbox{strongly in } W^{2,2}(\sur),
	\mbox{weakly}^* \mbox{ in } W^{1,\infty}(\sur)
\end{equation}
and, as $\ D f \in W^{1,2}\ $,
that locally uniformly $\ \sup_{|x-y| \leq C/m}
d(D f_m(x),D f(y)) \rightarrow 0\ $.
This implies for that the pull-backs
are uniformly bounded from below and above
\begin{equation} \label{direct.imm.comp}
	c_0 g \leq f_m^* \geu \leq C g
\end{equation}
for some $\ 0 < c_0 \leq C < \infty
\mbox{ and } m\ $ large,
in particular $\ f_m\ $ are smooth immersions.


\setcounter{equation}{0}

\section{The full rank case} \label{full}

For a smooth immersion $\ f: \sur \rightarrow \rel^n\ $
of a closed, orientable surface $\ \sur \mbox{ of genus } p \geq 1\ $
and $\ \vari \in C^\infty(\sur,\rel^n)\ $,
we see that the variations $\ f + t \vari\ $
are still immersions for small $\ t\ $.
We put $\ g_t := (f + t \vari)^* \geu,
g = g_0 = f^* \geu\ $
and define
\begin{equation} \label{full.vari}
	\delta \tp_f . \vari := d \tp_g . {\partial_t g_t}_{|t = 0}.
\end{equation}
The elements of
\begin{displaymath}
	{\cal V}_f := \delta \tp_f . C^\infty(\sur,\rel^n)
	\subseteq T_{\tp g} \teich
\end{displaymath}
can be considered
as the infinitessimal variations of $\ g\ $
in Teichm\"uller space obtained by ambient variations of $\ f\ $.
We call $\ f\ $ of full rank in Teichm\"uller space,
if $\ \dim {\cal V}_f = \dim\ \teich\ $.
In this case, the necessary corrections in Teichm\"uller space
mentioned in \S \ref{direct} can easily be achieved
by the inverse function theorem,
as we will see in this section.

Writing $\ g = e^{2u} \gp\ $
for some unit volume constant curvature metric $\ \gp\ $
by Poincar\'e's Theorem,
see \bcite{fisch.trom.conf}, \bcite{trom.teich},
we see $\ \tp(g_t) = \tp(e^{-2u} g_t)\ $,
hence for an orthonormal basis
$\ q^r(\gp), r = 1, \ldots, \dim\ \teich,\ $
of transvere traceless tensors in
$\ S_2^{TT}(\gp)\ $ with respect to $\ \gp\ $
\begin{equation} \label{full.vari-poin}
	\delta \tp_f . \vari
	= d \tp_g . {\partial_t g_t}_{|t = 0}
	= d \tp_\gp . {e^{-2u} \partial_t g_t}_{|t = 0} =
\end{equation}
\begin{displaymath}
	= \sum \limits_{r=1}^{\dim\ \teich}
	\langle {e^{-2u} \partial_t g_t}_{|t = 0} . q^r(\gp) \rangle_\gp
	d \tp_\gp . q^r(\gp).
\end{displaymath}
Calculating in local charts
\begin{equation} \label{full.metric-vari}
	g_{t,ij} = g_{ij}
	+ t \langle \partial_i f , \partial_j \vari \rangle
	+ t \langle \partial_j f , \partial_i \vari \rangle
	+ t^2 \langle \partial_i \vari , \partial_j \vari \rangle,
\end{equation}
we get
\begin{displaymath}
	\langle {e^{-2u} \partial_t g_t}_{|t = 0} , q^r(\gp) \rangle_\gp 
	= \int \limits_\sur \gpo{ik} \gpo{jl}
	{e^{-2 u} \partial_t g_{t,ij}}_{|t = 0} q^r_{kl}(\gp) \d \mu_\gp =
\end{displaymath}
\begin{displaymath}
	= 2 \int \limits_\sur g^{ik} g^{jl}
	\langle \partial_i f , \partial_j \vari \rangle
	q^r_{kl}(\gp) \d \mu_g =
\end{displaymath}
\begin{displaymath}
	= - 2 \int \limits_\sur g^{ik} g^{jl}
	\langle \nabla^g_j \nabla^g_i f , \vari \rangle
	q^r_{kl}(\gp) \d \mu_g
	- 2 \int \limits_\sur g^{ik} g^{jl}
	\langle \partial_i f , \vari \rangle
	\nabla^g_j q^r_{kl}(\gp) \d \mu_g =
\end{displaymath}
\begin{displaymath}
	= - 2 \int \limits_\sur g^{ik} g^{jl}
	\langle A_{ij}^0 , \vari \rangle
	q^r_{kl}(\gp) \d \mu_g,
\end{displaymath}
as $\ q \in S_2^{TT}(\gp) = S_2^{TT}(g)\ $ is divergence-
and tracefree with respect to $\ g\ $.
Therefore
\begin{displaymath}
	\delta \tp_f . \vari
	= \sum \limits_{r=1}^{\dim\ \teich}
	2 \int \limits_\sur g^{ik} g^{jl}
	\langle \partial_i f , \partial_j \vari \rangle
	q^r_{kl}(\gp) \d \mu_g
	\ d \tp_\gp . q^r(\gp) =
\end{displaymath}
\begin{equation} \label{full.vari-form}
	= \sum \limits_{r=1}^{\dim\ \teich}
	- 2 \int \limits_\sur g^{ik} g^{jl}
	\langle A_{ij}^0 , \vari \rangle
	q^r_{kl}(\gp) \d \mu_g
	\ d \tp_\gp . q^r(\gp),
\end{equation}
and $\ \delta \tp_f \mbox{ and } {\cal V}_f\ $
are well defined for
weak local bilipschitz immersions $\ f
\mbox{ and } V \in W^{1,2}(\sur)\ $.

First, we select variations whoose image via $\ \delta \tp_f\ $
form a basis of $\ {\cal V}_f\ $.

\begin{proposition} \label{full.basis}

For a weak local bilipschitz immersion $\ f\ $
and finitely many points $\ x_0, \ldots, x_N \in \sur\ $,
there exist $\ \vari_1, \ldots, \vari_{\dim {\cal V}_f}
\in C^\infty_0(\sur - \{x_0, \ldots, x_N\},\rel^n)\ $ such that
\begin{equation} \label{full.basis.span}
	{\cal V}_f = span\ \{ \delta \tp_f . \vari_s\ |
	\ s = 1, \ldots, \dim\ {\cal V}_f\ \}.
\end{equation}
\end{proposition}
{\pr Proof:} \\
Clearly $\ \delta \tp_f: C^\infty(\sur,\rel^n)
\rightarrow T_{\tp \gp} \teich\ $ is linear.
For suitable $\ W \subseteq C^\infty(\sur,\rel^n)\ $,
we obtain a direct sum decomposition
\begin{displaymath}
	C^\infty(\sur,\rel^n) = ker\ \delta \tp_f \oplus W
\end{displaymath}
and see that $\ \delta \tp_f | W \rightarrow im\ \delta \tp_f\ $
is bijective, hence $\ \dim\ W = \dim\ {\cal V}_f\ $
and a basis $\ \vari_1, \ldots, \vari_{\dim {\cal V}_f} \mbox{ of } W\ $
satisfies (\ref{full.basis.span}).

Next we choose cutoff functions $\ \varphi_\varrho
\in C^\infty_0(\cup_{k=0}^N B_\varrho(x_k))
\mbox{ such that } 0 \leq \varphi_\varrho \leq 1,
\varphi_\varrho = 1 \mbox{ on } \cup_{k=0}^N B_{\varrho/2}(x_k)
\mbox{ and } |\nabla \varphi_\varrho|_g \leq C \varrho^{-1}\ $,
hence $\ 1 - \varphi_\varrho \in C^\infty_0(\sur - \{x_0, \ldots, x_N\}),
1 - \varphi_\varrho \rightarrow 1 \mbox{ on }
\sur - \{x_0, \ldots, x_N\},
\int_\sur |\nabla \varphi_\varrho|_g \d \mu_g \leq C \varrho \rightarrow 0
\mbox{ for } \varrho \rightarrow 0\ $.
Clearly $\ \varphi_m \vari_s \in
C^\infty_0(\sur - \{x_0, \ldots, x_N\},\rel^n)\ $
and by (\ref{full.vari-form})
\begin{displaymath}
	\delta \tp_f . (\varphi_m \vari_s)
	\rightarrow \delta \tp_f . \vari_s,
\end{displaymath}
hence $\ \varphi_m \vari_1, \ldots, \varphi_m \vari_{\dim {\cal V}_f}\ $
satisfies for large $\ m\ $
all conclusions of the proposition.
\proof
We continue with a convergence criterion
for the first variation.

\begin{proposition} \label{full.conv}

Let $\ f: \sur \rightarrow \rel^n\ $ be a weak local bilipschitz immersion
approximated by smooth immersions $\ f_m\ $
with pull-back metrics $\ g = f^* \geu = e^{2u} \gp,
g_m = f_m^* \geu = e^{2u_m} \gp_m\ $
for some smooth unit volume constant
curvature metrics $\ \gp, \gpm\ $
and satisfying
\begin{equation} \label{full.conv.ass}
\begin{array}{c}
	f_m \rightarrow f
	\quad \mbox{weakly in } W^{2,2}(\sur),
	\mbox{weakly}^* \mbox{ in } W^{1,\infty}(\sur), \\
	\Lambda^{-1} \gp \leq g_m \leq \Lambda \gp, \\
	\parallel u_m \parallel_{L^\infty(\sur)} \leq \Lambda
\end{array}
\end{equation}
for some $\ \Lambda < \infty\ $.
Then for any $\ W \in L^2(\sur,\rel^n)\ $
\begin{displaymath}
	\tp_{f_m} . W \rightarrow \tp_f . W
\end{displaymath}
as $\ m \rightarrow \infty\ $.
\end{proposition}
{\pr Proof:} \\
By (\ref{full.conv.ass})
\begin{equation} \label{full.conv.met-conv}
\begin{array}{c}
	g_m = f_m^* \geu \rightarrow f^* \geu = g
	\quad \mbox{weakly in } W^{1,2}(\sur),
	\mbox{weakly}^* \mbox{ in } L^\infty(\sur), \\
	\Gamma_{g_m,ij}^k \rightarrow \Gamma_{g,ij}^k
	\quad \mbox{weakly in } L^2, \\
	A_{f_m,ij}, A_{f_m,ij}^0
	\rightarrow A_{f,ij}, A_{f,ij}^0
	\quad \mbox{weakly in } L^2, \\
	\int \limits_\sur |K_{g_m}| \d \mu_{g_m}
	\leq \frac{1}{2} \int \limits_\sur
	|A_{f_m}|^2 \d \mu_{g_m} \leq C,
\end{array}
\end{equation}
hence by \bcite{fisch.trom.conf}, \bcite{trom.teich},
\begin{equation} \label{full.conv.conv}
\begin{array}{c}
	q^r(\gpm) \rightarrow q^r(\gp)
	\left\{
	\begin{array}{c}
		\mbox{weakly in } W^{1,2}(\sur), \\
		\mbox{weakly}^*
		\mbox{ in } L^\infty(\sur),
	\end{array}
	\right. \\
	d \tp_\gpm . q^r(\gpm)
	\rightarrow d \tp_\gp . q^r(\gp).
\end{array}
\end{equation}
Then by (\ref{full.vari-form})
\begin{displaymath}
	\delta \tp_{f_m} . W =
\end{displaymath}
\begin{displaymath}
	= \sum \limits_{r=1}^{\dim\ \teich}
	- 2 \int \limits_\sur g_m^{ik} g_m^{jl}
	\langle A_{f_m,ij}^0 , W \rangle
	q^r_{kl}(\gpm) \d \mu_{g_m}
	\ d \tp_\gpm . q^r(\gpm) \rightarrow
\end{displaymath}
\begin{displaymath}
	\rightarrow \delta \tp_f . W.
\end{displaymath}
\proof
In the full rank case, the necessary corrections in Teichm\"uller space
mentioned in \S \ref{direct} are achieved
in the following lemma by the inverse function theorem.

\begin{lemma} \label{full.corr}

Let $\ f: \sur \rightarrow \rel^n\ $ be a weak local bilipschitz immersion
approximated by smooth immersions $\ f_m\ $
satisfying (\ref{direct.prop.teich}) - (\ref{direct.prop.imm}).

If $\ f\ $ is of full rank in Teichm\"uller space,
then for arbitrary $\ x_0 \in \sur,
\mbox{ neighbourhood } U_\mark(x_0)
\subseteq \sur \mbox{ of } x_0
\mbox{ and } \Lambda < \infty\ $,
there exists a neighbourhood $\ U(x_0) \subseteq U_\mark(x_0)
\mbox{ of } x_0\ $,
variations $\ \vari_1, \ldots, \vari_{\dim \teich}
\in C^\infty_0(\sur - \overline{U(x_0)},\rel^n)\ $,
satisfying (\ref{full.basis.span}),
and $\ \delta > 0, C < \infty, m_0 \in \nat\ $
such that for any $\ \vari \in C_0^\infty(U(x_0),\rel^n)\ $
with $\ f_m + \vari\ $ a smooth immersion
for some $\ m \geq m_0\ $, and $\ \vari = 0\ $ or
\begin{equation} \label{full.corr.ass}
\begin{array}{c}
	\Lambda^{-1} \gp \leq (f_m + \vari)^* \geu \leq \Lambda \gp, \\
	\parallel \vari \parallel_{W^{2,2}(\sur)}
	\leq \Lambda, \\
	\int \limits_{U_\mark(x_0)} |A_{f_m + \vari}|^2 \d \mu_{f_m + \vari}
	\leq \varepsilon_0(n),
\end{array}
\end{equation}
where $\ \varepsilon_0(n)\ $ is as in Lemma \ref{conf.lemma},
and any $\ \tau \in \teich\ $ with
\begin{equation} \label{full.corr.teich-point}
	d_\teich(\tau,\tau_0) \leq \delta,
\end{equation}
there exists $\ \lambda \in \rel^{\dim \teich}\ $ satisfying
\begin{displaymath}
	\tp( (f_m + \vari
	+ 
	\lambda_r \vari_r)^* \geu)
	= \tau
\end{displaymath}
and
\begin{displaymath}
	|\lambda| \leq C d_\teich \Big( \tp( (f_m + \vari)^* \geu),
	\tau \Big).
\end{displaymath}
\end{lemma}
{\pr Proof:} \\
By (\ref{direct.prop.w22}), (\ref{direct.prop.met}) and $\ \Lambda\ $
large enough, we may assume
\begin{equation} \label{full.corr.d2}
\begin{array}{c}
	\parallel u_m, D f_m \parallel_{W^{1,2}(\sur)
	\cap L^\infty(\sur)} \leq \Lambda, \\
	\Lambda^{-1} \gp \leq f_m^* \geu \leq \Lambda \gp,
\end{array}
\end{equation}
in particular
\begin{equation} \label{full.corr.fund}
	\int \limits_\sur |A_{f_m}|^2 \d \mu_{f_m}
	\leq C(\sur,\gp,\Lambda).
\end{equation}
Putting $\ \funda_m := |\nabla^2_\gp f_m|_\gp^2 \mu_{\gp}\ $,
we see $\ \funda_m(\sur) \leq C(\Lambda,\gp) $
and for a subsequence $\ \funda_m \rightarrow \funda\ $ 
weakly$^* \mbox{ in } C_0^0(\sur)^*\ $.
Clearly $\ \funda(\sur) < \infty\ $,
and there are at most finitely many $\ y_1, \ldots, y_N \in \sur
\mbox{ with } \funda(\{y_i\}) \geq \varepsilon_1\ $,
where we choose $\ \varepsilon_1 > 0\ $ below.

As $\ f\ $ is of full rank,
we can select $\ \vari_1, \ldots, \vari_{\dim \teich}
\in C^\infty_0(\sur - \{x_0, y_1, \ldots, y_N\},\rel^n)
\mbox{ with } span \{ \delta \tp_f . \vari_r \} = T_\gp \teich\ $
by Proposition \ref{full.basis}.
We choose a neighbourhood $\ U_0(x_0)
\subseteq U_\mark(x_0) \mbox{ of } x_0\ $
with a chart $\ \varphi_0: U_0(x_0)
\stackrel{\approx}{\longrightarrow} B_2(0),
\varphi_0(x_0) = 0\ $,
\begin{displaymath}
	supp\ \vari_r \cap \overline{U_0(x_0)} = \emptyset
	\quad \mbox{for } r = 1, \ldots, \dim \teich,
\end{displaymath}
put $\ x_0 \in U_\varrho(x_0) = \varphi_0^{-1}(B_\varrho(0))
\mbox{ for } 0 < \varrho \leq 2\ $
and choose $\ x_0 \in U(x_0) \subseteq U_{1/2}(x_0)\ $
small enough, as we will see below.

Next for any $\ x \in \cup_{r=1}^{\dim \teich} supp\ \vari_r\ $,
there exists a neighbourhood $\ U_0(x) \mbox{ of } x\ $
with a chart $\ \varphi_x: U_0(x)
\stackrel{\approx}{\longrightarrow} B_2(0),
\varphi_x(x) = 0,
\overline{U_0(x)} \cap \overline{U_0(x_0)} = \emptyset,
\funda(U_0(x)) < \varepsilon_1\ $ and
in the coordinates of the chart $\ \varphi_x\ $
\begin{equation} \label{full.corr.christ}
	\int \limits_{U_0(x)} \gpo{ik} \gpo{jl} \gpu{rs}
	\Gamma_{\gp,ij}^r \Gamma_{\gp,kl}^s \d \mu_\gp
	\leq \varepsilon_1.
\end{equation}
Putting $\ x \in U_\varrho(x) = \varphi_x^{-1}(B_\varrho(0))
\subset \subset U_0(x) \mbox{ for } 0 < \varrho \leq 2\ $,
we see that there are finitely many $\ x_1, \ldots, x_\M
\in \cup_{r=1}^{\dim \teich} supp\ \vari_r\ $
such that
\begin{displaymath}
	\cup_{r=1}^{\dim \teich} supp\ \vari_r
	\subseteq \cup_{k=1}^\M U_{1/2}(x_k).
\end{displaymath}
Then there exists $\ m_0 \in \nat\ $
such that for $\ m \geq m_0\ $
\begin{displaymath}
	\int \limits_{U_1(x_k)} |\nabla^2_\gp f_m|^2_\gp
	\d \mu_{\gp} < \varepsilon_1
	\quad \mbox{for } k = 1, \ldots, \M.	
\end{displaymath}
For $\ \vari \mbox{ and } m \geq m_0\ $ as above,
we put $\ \tilde f_{m,\lambda} := f_m + \vari + \lambda_r \vari_r\ $.
Clearly
\begin{displaymath}
	supp(f_m - \tilde f_{m,\lambda})
	\subseteq \cup_{k=0}^\M U_{1/2}(x_k).
\end{displaymath}
By (\ref{full.corr.ass}), (\ref{full.corr.d2}), (\ref{full.corr.christ}),
and $\ |\lambda| < \lambda_0 \leq 1/4\ $ small enough
independent of $\ m \mbox{ and } \vari\ $,
$\ \tilde f_{m,\lambda}\ $ is a smooth immersion with
\begin{equation} \label{full.corr.met}
	(2 \Lambda)^{-1} \gp \leq \tilde g_{m,\lambda}
	:= \tilde f_{m,\lambda}^* \geu
	\leq 2 \Lambda \gp,
\end{equation}
and if $\ \vari \neq 0\ $ by (\ref{full.corr.ass})
and the choice of $\ U_0(x)\ $ that
\begin{displaymath}
	\int \limits_{U_1(x_k)} |A_{\tilde f_{m,\lambda}}|^2
	\d \mu_{\tilde g_{m,\lambda}} \leq \varepsilon_0(n)
	\quad \mbox{for } k = 0, \ldots, \M
\end{displaymath}
for $\ C(\Lambda,\gp) (\varepsilon_1 + \lambda_0) \leq \varepsilon_0(n)\ $.
If $\ \vari = 0\ $, we see
$\ supp(f_m - \tilde f_{m,\lambda})
\subseteq \cup_{k=1}^\M U_{1/2}(x_k)\ $.
Further by (\ref{full.corr.fund})
\begin{displaymath}
	\int \limits_\sur |K_{\tilde g_{m,\lambda}}| \d \mu_{\tilde g_{m,\lambda}}
	\leq \frac{1}{2} \int \limits_\sur |A_{\tilde f_{m,\lambda}}|^2
	\d \mu_{\tilde g_{m,\lambda}} \leq
\end{displaymath}
\begin{displaymath}
	\leq \frac{1}{2} \int \limits_\sur |A_{f_m}|^2 \d \mu_{f_m}
	+  \frac{1}{2} \sum \limits_{k=0}^\M \int \limits_{U_{1/2}(x_k)}
	|A_{\tilde f_{m,\lambda}}|^2 \d \mu_{\tilde g_{m,\lambda}}
	\leq C(\sur,\gp,\Lambda) + (\M + 1) \varepsilon_0(n).
\end{displaymath}
This verifies (\ref{conf.lemma.metric2})
and (\ref{conf.lemma.fund})
for $\ f = f_m, \tilde f = \tilde f_{m,\lambda}, g_0 = \gp\ $
and different, but appropriate $\ \Lambda\ $.
(\ref{conf.lemma.ub}) follows from (\ref{direct.prop.met})
and (\ref{full.corr.fund}).
Then for the unit volume constant curvature metric
$\ \tgpml = e^{-2 \tilde u_{m,\lambda}} \tilde g_{m,\lambda}\ $
conformal to $\ \tilde g_{m,\lambda}\ $
by Poincar\'e's Theorem,
see \bcite{fisch.trom.conf}, \bcite{trom.teich},
we get from Lemma \ref{conf.lemma} that
\begin{equation} \label{full.corr.conf}
	\parallel \tilde u_{m,\lambda} \parallel_{L^\infty(\sur)},
	\parallel \nabla \tilde u_{m,\lambda}
	\parallel_{L^2(\sur,\gp)} \leq C
\end{equation}
with $\ C < \infty\ $ independent of $\ m \mbox{ and } \vari\ $.

From (\ref{full.corr.ass}),
we have a $\ W^{2,2} \cap W^{1,\infty}-$bound on $\ \tilde f_0\ $,
hence for $\ \tilde f_{m,\lambda}\ $.
On $\ \sur - U(x_0)\ $, we get
$\ \tilde f_{m,\lambda} = f_m + \lambda_r \vari_r
\rightarrow f
\mbox{ weakly in } W^{2,2}(\sur - U(x_0))
\mbox{ and weakly}^* \mbox{ in } W^{1,\infty}(\sur - U(x_0))
\mbox{ for } m_0 \rightarrow \infty, \lambda_0 \rightarrow 0\ $
by (\ref{direct.prop.conv}).
If $\ \vari = 0\ $,
then $\ \tilde f_{m,\lambda} = f_m + \lambda_r \vari_r \rightarrow f
\mbox{ weakly in } W^{2,2}(\sur)
\mbox{ and weakly}^* \mbox{ in } W^{1,\infty}(\sur)
\mbox{ for } m_0 \rightarrow \infty, \lambda_0 \rightarrow 0\ $
by (\ref{direct.prop.conv}).
Hence letting $\ m_0 \rightarrow \infty,
\lambda_0 \rightarrow 0, U(x_0) \rightarrow \{x_0\}\ $,
we conclude
\begin{equation} \label{full.corr.f-conv}
	\tilde f_{m,\lambda} \rightarrow f
	\quad \mbox{weakly in } W^{2,2}(\sur),
	\mbox{weakly}^* \mbox{ in } W^{1,\infty}(\sur),
\end{equation}
in particular
\begin{displaymath}
	\tilde g_{m,\lambda} \rightarrow f^* \geu = e^{2u} \gp =: g
	\left\{
	\begin{array}{c}
		\mbox{weakly in } W^{1,2}(\sur), \\
		\mbox{weakly}^*
		\mbox{ in } L^\infty(\sur).
	\end{array}
	\right.
\end{displaymath}
Together with (\ref{full.corr.met}),
this implies by \bcite{fisch.trom.conf}, \bcite{trom.teich},
\begin{equation} \label{full.corr.teich-conv}
	\tp(\tilde g_{m,\lambda}) \rightarrow \tau_0.
\end{equation}
We select a chart $\ \psi: U(\tp(\gp)) \subseteq \teich
\rightarrow \rel^{\dim \teich}\ $
and put $\ \tpc := \psi \circ \tp,
\delta \tpc = d \psi \circ d \tp,
\hat {\cal V}_f := d \psi_{\tp \gp} . {\cal V}_f\ $.
By (\ref{full.corr.teich-conv}) for $\ m_0\ $ large enough,
$\ \lambda_0 \mbox{ and } U(x_0)\ $ small enough independent of $\ \vari\ $,
we get $\ \tp(\tilde g_{m,\lambda}) \in U(\tau_0)\ $
and define
\begin{displaymath}
        \Phi_m(\lambda) := \tpc(\tilde g_{m,\lambda}).
\end{displaymath}
This yields by (\ref{full.vari-form}),
(\ref{full.corr.met}), (\ref{full.corr.conf}),
(\ref{full.corr.f-conv})
and Proposition \ref{full.conv}
\begin{displaymath}
	D \Phi_m(\lambda)
	= (\delta \tpc_{\tilde f_{m,\lambda}} .
	\vari_r)_{r = 1, \ldots, \dim \teich}
	\rightarrow (\delta \tpc_f . \vari_r)_{r = 1, \ldots, \dim \teich}
	=: A \in \rel^{\dim \teich \times \dim \teich}.
\end{displaymath}
As $\ span \{ \delta \tp_f . \vari_r \} = T_\gp \teich
\cong \rel^{\dim \teich}\ $,
the matrix $\ A\ $ is invertible.
Choosing $\ m_0\ $ large enough,
$\ \lambda_0 \mbox{ and } U(x_0)\ $ small enough,
we obtain
\begin{displaymath}
	\parallel D \Phi_m(\lambda) - A \parallel
	\leq 1/(2 \parallel A^{-1} \parallel),
\end{displaymath}
hence by standard inverse function theorem
for any $\ \xi \in \rel^{\dim \teich}
\mbox{ with } |\xi - \Phi_m(0)| < \lambda_0 / (2 \parallel A^{-1} \parallel)\ $
there exists $\ \lambda \in B_{\lambda_0}(0)\ $ with
\begin{displaymath}
\begin{array}{c}
	\Phi_m(\lambda) = \xi, \\
	|\lambda| \leq 2 \parallel A^{-1} \parallel\ |\xi - \Phi_m(0)|.
\end{array}
\end{displaymath}
As
\begin{displaymath}
	d_\teich \Big( \tp( (f_m + \vari)^* \geu), \tau \Big)
	\leq d_\teich(\tp(\tilde g_0),\tau_0)
	+ d_\teich ( \tau_0,\tau )
	< d_\teich(\tp(\tilde g_0),\tau_0) + \delta
\end{displaymath}
by (\ref{full.corr.teich-point}),
we see for $\ \delta\ $ small enough,
$\ m_0\ $ large enough
and $\ U(x_0)\ $ small enough
by (\ref{full.corr.teich-conv})
that there exists $\ \lambda \in \rel^{\dim \teich}\ $ satisfying
\begin{displaymath}
\begin{array}{c}
	\tp( (f_m + \vari + \lambda_r \vari_r)^* \geu)
	= \tp(\tilde g_{m,\lambda})
	= \tau, \\
	|\lambda| \leq C d_\teich \Big( \tp( (f_m + \vari)^* \geu),
	\tau \Big),
\end{array}
\end{displaymath}
and the lemma is proved.
\proof


\setcounter{equation}{0}

\section{The degenerate case} \label{deg}

In this section, we consider the degenerate case
when the immersion is not of full rank
in Teichm\"uller space.
First we see that the image in Teichm\"uller space
looses at most one dimension.

\begin{proposition} \label{deg.image}

For a weak local bilipschitz immersion $\ f \in W^{2,2}(\sur,\rel^n)\ $
we always have
\begin{equation} \label{deg.image.dim}
	\dim\ {\cal V}_f \geq \dim\ \teich - 1.
\end{equation}
If $\ f\ $ is not of full rank in Teichm\"uller space
then $\ f\ $ is isothermic locally
around all but finitely many points of $\ \sur\ $,
that is around all but finitely many points of $\ \sur\ $
there exist local conformal principal curvature coordinates.
\end{proposition}
{\pr Proof:} \\
For $\ q \in S_2^{TT}(\gp)\ $,
we put $\ \Lambda_q: C^\infty(\sur,\rel^n) \rightarrow \rel\ $
\begin{displaymath}
	\Lambda_q . \vari := - 2 \int \limits_\sur g^{ik} g^{jl}
	\langle A_{ij}^0, \vari \rangle
	q_{kl} \d \mu_g
\end{displaymath}
and define the annihilator
\begin{displaymath}
	Ann := \{ q \in S_2^{TT}(\gp)\ |\ \Lambda_q = 0\ \}.
\end{displaymath}
As $\ d \pi_\gp | S_2^{TT}(\gp) \rightarrow T_\gp \teich\ $
is bijective, we see by (\ref{full.vari-form})
and elementary linear algebra
\begin{equation} \label{deg.image.rank}
	\dim\ \teich = \dim\ {\cal V}_f + \dim\ Ann.
\end{equation}
Clearly,
\begin{equation} \label{deg.image.ann}
	q \in Ann \Longleftrightarrow
	g^{ik} g^{jl} A_{ij}^0 q_{kl} = 0.
\end{equation}
Choosing a conformal chart with respect to the smooth metric $\ \gp\ $,
we see $\ g_{ij} = e^{2v} \delta_{ij}
\mbox{ for some } v \in W^{1,2} \cap L^\infty\ $
and $\ A^0_{11} = - A^0_{22}, A^0_{12} = A^0_{21},
q_{11} = - q_{22}, q_{12} = q_{21}\ $,
as both $\ A^0 \mbox{ and } q \in S_2^{TT}(\gp)\ $
are symmetric and tracefree with respect to $\ g = e^{2u} \gp\ $.
This rewirtes (\ref{deg.image.ann}) into
\begin{equation} \label{deg.image.ann2}
	q \in Ann \Longleftrightarrow
	A_{11}^0 q_{11} + A_{12}^0 q_{12} = 0.
\end{equation}
The correspondance between $\ S_2^{TT}(\gp)\ $
and the holomorphic quadratic differentials
is exactly that in conformal coordinates
\begin{equation} \label{deg.image.holo}
	h := q_{11} - i q_{12} \mbox{ is holomorphic}.
\end{equation}
Now if (\ref{deg.image.dim}) were not true,
there would be two linearly independent $\ q^1, q^2 \in Ann\ $
by (\ref{deg.image.rank}).
Likewise the holomorphic functions $\ h^k := q_{11}^k - i q_{12}^k\ $
are linearly independant over $\ \rel\ $,
in particular neither of them vanishes identically,
hence these vanish at most at finitely many points,
as $\ \sur\ $ is closed.
Then $\ h^1/h^2\ $ is meromorphic
and moreover not a real constant.
This implies that $\ Im(h^1/h^2)\ $ does not vanish identically,
hence vanishes at most at finitely many points.
Outside these finitely many points, we calculate
\begin{displaymath}
	Im(h^1 / h^2) = |h^2|^{-2} Im(h^1 \overline{h^2})
	= |h^2|^{-2} \det \left(
	\begin{array}{cc}
		q_{11}^1 & q_{12}^1 \\
		q_{11}^2 & q_{12}^2 \\
	\end{array}
	\right),
\end{displaymath}
hence
\begin{displaymath}
        \det \left(
        \begin{array}{cc}
                q_{11}^1 & q_{12}^1 \\
                q_{11}^2 & q_{12}^2 \\
        \end{array}
        \right)
        \mbox{ vanishes at most at finitely many points}
\end{displaymath}
and by (\ref{deg.image.ann2})
\begin{equation} \label{deg.image.nab}
	A^0 = 0.
\end{equation}
Approximating $\ f\ $ by smmoth immersions
as in (\ref{direct.imm.comp}), we get
\begin{displaymath}
	\nabla_k \mean = 2 g^{ij} \nabla_i A^0_{jk} = 0
	\quad \mbox{ weakly},
\end{displaymath}
where $\ \nabla = D^\perp\ $ denotes
the normal connection in the normal bundle along $\ f\ $. 
Therefore $\ |\mean|\ $ is constant and
\begin{displaymath}
	\partial_k \mean
	= g^{ij} \langle \partial_k \mean ,
	\partial_j f \rangle \partial_i f
	= - g^{ij} \langle \mean , A_{jk} \rangle \partial_i f
	\quad \mbox{ weakly}.
\end{displaymath}
Using $\ \Delta_\gp f = e^{2u} \Delta_g f = e^{2u} \mean
\mbox{ and } u \in W^{1,2} \cap L^\infty\ $,
we conclude successively that $\ f \in C^\infty\ $.
Then (\ref{deg.image.nab}) implies
that $\ f\ $ parametrizes a round sphere or a plane,
contradicting $\ p \geq 1\ $,
and (\ref{deg.image.dim}) is proved.
%

Next if $\ f\ $ is not of full rank in Teichm\"uller space,
there exits $\ q \in Ann - \{0\} \neq \emptyset\ $
by (\ref{deg.image.rank}),
and the holomorphic function $\ h\ $ in (\ref{deg.image.holo})
vanishes at most at finitely many points.
In a neighbourhood of a point where $\ h\ $ does not vanish,
there is a holomorphic function $\ w \mbox{ with } (w')^2 = i h\ $.
Then $\ w\ $ has a local inverse $\ z\ $
and using $\ w\ $ as new local conformal coordinates,
$\ h\ $ transforms into $\ (h \circ z) (z')^2 = -i\ $,
hence $\ q_{11} = 0, q_{12} = 1\ $
in $\ w-$coordinates.
By (\ref{deg.image.ann2})
\begin{displaymath}
	A_{12} = 0,
\end{displaymath}
and $\ w\ $ are local conformal principal curvature coordinates.
\proof
Since we loose at most one dimension in the degenerate case,
we will do the necessary corrections in Teichm\"uller space
mentioned in \S \ref{direct} by investigating
the second variation in Teichm\"uller space.
To be more precise, for a smooth immersion
$\ f: \sur \rightarrow \rel^n\ $
with pull-back metric $\ f^* \geu\ $
conformal to a unit volume constant curvature metric
$\ \gp = e^{-2u} f^* \geu\ $
by Poincar\'e's Theorem,
see \bcite{fisch.trom.conf}, \bcite{trom.teich},
we select a chart $\ \psi: U(\tp(\gp)) \subseteq \teich
\rightarrow \rel^{\dim \teich}\ $,
put $\ \tpc := \psi \circ \tp,
\hat {\cal V}_f := d \psi_{\tp \gp} . {\cal V}_f\ $,
and define the second variation in Teichm\"uller space
of $\ f\ $ with respect to the chart $\ \psi\ $ by
\begin{equation} \label{deg.sec-def}
	\delta^2 \tpc_f(\vari)
	:= \Big( \frac{d}{dt} \Big)^2
	\tpc( (f + t \vari)^* \geu)_{|t = 0}.
\end{equation}
Putting $\ f_t := f + t \vari, g_t := f_t^* \geu, g = g_0\ $,
we see $\ \tpc(g_t) = \tpc(e^{-2u} g_t)\ $ and calculate
\begin{displaymath}
	\delta^2 \tpc_f(\vari)
	= \Big( \frac{d}{dt} \Big)^2
	\tpc( e^{-2u} g_t )_{|t = 0} =
\end{displaymath}
\begin{displaymath}
	= d \tpc_\gp . (e^{-2u} (\partial_{tt} g_t)_{|t = 0})
	+ d^2 \tpc_\gp(e^{-2u} (\partial_t g_t)_{|t = 0},
	e^{-2u} (\partial_t g_t)_{|t = 0}).
\end{displaymath}
Decomposing
\begin{equation} \label{deg.decomp}
	e^{-2u} (\partial_t g_t)_{|t = 0}
	= \fu \gp + {\cal L}_X \gp + q,
\end{equation}
with $\ \fu \in C^\infty(\sur),
X \in \vf(\sur), q \in S_2^{TT}(\gp)\ $,
we continue recalling $\ \{ \fu \gp + {\cal L}_X \gp \}
= ker\ d \tpc_\gp\ $
by \bcite{fisch.trom.conf}, \bcite{trom.teich},
and
\begin{displaymath}
	\delta^2 \tpc_f(\vari)
	= d \tpc_\gp . (e^{-2u} (\partial_{tt} g_t)_{|t = 0})
	+ d^2 \tpc_\gp(\fu \gp + {\cal L}_X \gp + q,
	\fu \gp + {\cal L}_X \gp + q) =
\end{displaymath}
\begin{displaymath}
	= d \tpc_\gp . (e^{-2u} (\partial_{tt} g_t)_{|t = 0})
	+ d^2 \tpc_\gp \hspace{-.1cm} (q,q)
	+ d^2 \tpc_\gp \hspace{-.1cm} (\fu \gp + {\cal L}_X \gp,
	\fu \gp + {\cal L}_X \gp + 2 q) =
\end{displaymath}
\begin{displaymath}
	= d \tpc_\gp . \Big( e^{-2u} (\partial_{tt} g_t)_{|t = 0}
	- {\cal L}_X {\cal L}_X \gp
	- 2 \fu {\cal L}_X \gp
	- 2 \fu q - 2 {\cal L}_X q \Big)
	+ d^2 \tpc_\gp(q,q).
\end{displaymath}
For an orthonormal basis
$\ q^r(\gp), r = 1, \ldots, \dim\ \teich,\ $
of transverse traceless tensors in
$\ S_2^{TT}(\gp)\ $ with respect to $\ \gp\ $,
we obtain
\begin{equation} \label{deg.sec}
	\delta^2 \tpc_f(\vari)
	= \sum \limits_{r=1}^{\dim\ \teich}
	\alpha_r\ d \tpc_\gp . q^r(\gp)
	+ \sum \limits_{r,s=1}^{\dim\ \teich}
	\beta_r \beta_s\ d^2 \tpc_\gp(q^r(\gp),q^s(\gp)),
\end{equation}
where
\begin{displaymath}
\begin{array}{c}
	\alpha_r := \\
	= \int \limits_\sur \gpo{ik} \gpo{jl}
	\Big( e^{-2u} (\partial_{tt} g_t)_{|t=0}
	\hspace{-.1cm} - \hspace{-.1cm} {\cal L}_X {\cal L}_X \gp
	\hspace{-.1cm} - \hspace{-.1cm} 2 \fu {\cal L}_X \gp
	\hspace{-.1cm} - \hspace{-.1cm} 2 \fu q
	\hspace{-.1cm} - \hspace{-.1cm} 2 {\cal L}_X q \Big)_{ij}
	q^r_{kl}(\gp) \d \mu_\gp,
\end{array}
\end{displaymath}
\begin{equation} \label{deg.sec.coeff}
	\beta_r
	:= \int \limits_\sur \gpo{ik} \gpo{jl}
	q_{ij} q^r_{kl}(\gp) \d \mu_\gp.
\end{equation}
Since
\begin{displaymath}
	{\cal L}_X {\cal L}_X \gpu{ij}
	= \gpo{mo} \Big( X_m \nabla_o {\cal L}_X \gpu{ij}
	+ \nabla_i X_m {\cal L}_X \gpu{jo}
	+ \nabla_j X_m {\cal L}_X \gpu{oi} \Big),
\end{displaymath}
we get integrating by parts
\begin{displaymath}
	\alpha_r
	= \int \limits_\sur \gpo{ik} \gpo{jl}
	\Big( e^{-2u} (\partial_{tt} g_t)_{|t=0}
	- 2 \fu {\cal L}_X \gp
	- 2 \fu q - 2 {\cal L}_X q \Big)_{ij}
	q^r_{kl}(\gp) \d \mu_\gp +
\end{displaymath}
\begin{displaymath}
	+ \hspace{-.1cm} \int \limits_\sur \hspace{-.05cm}
	\gpo{ik} \gpo{jl} \gpo{mo} \hspace{-.1cm}
	\Big( \hspace{-.1cm} \nabla_o X_m {\cal L}_X \gpu{ij}
	\hspace{-.1cm} - \hspace{-.1cm}
	\nabla_i X_m {\cal L}_X \gpu{jo}
	\hspace{-.1cm} - \hspace{-.1cm}
	\nabla_j X_m {\cal L}_X \gpu{oi} \hspace{-.1cm} \Big)
	\hspace{-.05cm} q^r_{kl}(\gp)
	\hspace{-.1cm} \d \mu_\gp
	\hspace{-.1cm} +
\end{displaymath}
\begin{equation} \label{deg.sec.coeff2}
	+ \int \limits_\sur \gpo{ik} \gpo{jl} \gpo{mo}
	X_m {\cal L}_X \gpu{ij}
	\nabla_o q^r_{kl}(\gp) \d \mu_\gp.
\end{equation}
For a weak local bilipschitz immersion $\ f \in W^{2,2}(\sur,\rel^n)
\mbox{ and } \vari \in W^{1,2}(\sur,\rel^n)\ $,
we see $\ g \in W^{1,2}(\sur), u \in W^{1,2}(\sur) \cap L^\infty(\sur),
(\partial_t g_t)_{|t=0} \in L^2(\sur),
\mbox{ and } (\partial_{tt} g_t)_{|t=0}
\in L^1(\sur)\ $ by (\ref{full.metric-vari}).
Then we get a decomposition as in (\ref{deg.decomp})
with $\ \fu \in L^2(\sur), X \in \vf^1(\sur),
q \in S_2^{TT}(\sur) \subseteq S_2(\sur)\ $.
Observing that $\ q^r(\gp) \in S_2^{TT}(\sur) \subseteq S_2(\sur)\ $,
we conclude that $\ \delta^2 \tpc_f\ $ is well defined for
weak local bilipschitz immersions $\ f \in W^{2,2}(\sur,\rel^n)
\mbox{ and } \vari \in W^{1,2}(\sur,\rel^n)\ $
by (\ref{deg.sec}), (\ref{deg.sec.coeff})
and (\ref{deg.sec.coeff2}).

\begin{proposition} \label{deg.basis}

For a weak local bilipschitz immersion $\ f \in W^{2,2}(\sur,\rel^n)\ $,
which is not of full rank in Teichm\"uller space,
and finitely many points $\ x_1, \ldots, x_N \in \sur\ $,
there exist $\ \vari_1, \ldots, \vari_{\dim \teich - 1}, \vari_\pm
\in C^\infty_0(\sur - \{x_1, \ldots, x_N\},\rel^n)\ $ such that
\begin{equation} \label{deg.basis.span}
	\hat {\cal V}_f = span\ \{ \delta \tpc_f . \vari_s\ |
	\ s = 1, \ldots, \dim\ \teich - 1\ \}
\end{equation}
and for some $\ e \perp \hat {\cal V}_f, |e| = 1,\ $
\begin{equation} \label{deg.basis.pm}
	\pm \langle \delta^2 \tpc_f(\vari_\pm) , e \rangle > 0,
	\quad \delta \tpc_f . \vari_\pm = 0.
\end{equation}
\end{proposition}
{\pr Proof:} \\
By Proposition \ref{deg.image},
there exists $\ x_0 \in \sur - \{x_1, \ldots, x_N\}\ $
such that $\ f\ $ is isothermic around $\ x_0\ $.
Moreover, since $\ \sur\ $ is not a sphere,
hence $\ A^0\ $ does not vanish almost everywhere
with respect to $\ \mu_g\ $,
as we have seen in the argument after (\ref{deg.image.nab})
in Proposition \ref{deg.image},
we may assume that
\begin{equation} \label{deg.basis.pm.nab}
	x_0 \in supp\ |A^0|^2 \mu_g.
\end{equation}
By Proposition \ref{full.basis},
there exist $\ \vari_1, \ldots, \vari_{\dim \teich - 1}
\in C^\infty_0(\sur - \{x_0, \ldots, x_N\},\rel^n)\ $
which satisfy (\ref{deg.basis.span}).
Next we select a chart $\ \varphi: U(x_0)
\stackrel{\approx}{\longrightarrow} B_1(0)\ $
of conformal principal curvature coordinates
that is
\begin{equation} \label{deg.basis.pm.princ}
	g = e^{2v} \geu, \quad A_{12} = 0,
\end{equation}
in local coordinates of $\ \varphi\ $
and where $\ v \in W^{1,2}(B_1(0)) \cap L^\infty(B_1(0))\ $.
Moreover, we choose $\ U(x_0)\ $ so small that
$\ \overline{U(x_0)} \cap supp\ \vari_s = \emptyset
\mbox{ for } s = 1, \ldots, \dim\ \teich - 1,
\mbox{ and } \overline{U(x_0)} \cap \{ x_1, \ldots, x_N \}
= \emptyset\ $.
For any $\ \vari_0 \in W^{1,2}(\sur,\rel^n)
\cap L^\infty(\sur)
\mbox{ with } supp\ \vari_0 \subseteq U(x_0)\ $,
there exists a unique $\ \gamma \in \rel^{\dim \teich - 1}\ $
such that for $\ \vari := \vari_0 - \gamma_s \vari_s
\in W^{1,2}(\sur,\rel^n),
supp\ \vari \subseteq \sur - \{x_1, \ldots, x_N\}\ $
\begin{equation} \label{deg.basis.pm.null}
	\delta \tpc_f . \vari = 0.
\end{equation}
By (\ref{full.vari-form})
\begin{equation} \label{deg.basis.pm.null-esti}
	|\gamma| \leq C |\delta \tpc_f . \vari_0|
	\leq C \parallel A^0 \parallel_{L^2(supp \vari_0,g)}
	\ \parallel \vari_0 \parallel_{L^\infty(U(x_0))},
\end{equation}
where $\ C\ $ does not depend on $\ \vari_0\ $.
By (\ref{full.vari-poin}),
we get for the decomposition in (\ref{deg.decomp})
that $\ q = 0\ $.
We select the orthonormal basis
$\ q^r(\gp), r = 1, \ldots, \dim\ \teich,
\mbox{ of } S_2^{TT}(\gp)\ $ with respect to $\ \gp\ $,
in such a way that $\ d \tpc_\gp . q^r(\gp)
\in \hat {\cal V}_f \mbox{ for } r = 2, \ldots, \dim\ \teich,
\mbox{ and } \langle d \tpc_\gp . q^1(\gp) , e \rangle > 0\ $.
Putting
\begin{displaymath}
	I(\vari_0) := \int \limits_\sur \gpo{ik} \gpo{jl}
	\Big( e^{-2u} (\partial_{tt} g_t)_{|t=0}
	- 2 \fu {\cal L}_X \gp \Big)_{ij}
	q^1_{kl}(\gp) \d \mu_\gp +
\end{displaymath}
\begin{displaymath}
	+ \hspace{-.1cm} \int \limits_\sur \hspace{-.05cm}
	\gpo{ik} \gpo{jl} \gpo{mo} \hspace{-.1cm}
	\Big( \hspace{-.1cm} \nabla_o X_m {\cal L}_X \gpu{ij}
	\hspace{-.1cm} - \hspace{-.1cm}
	\nabla_i X_m {\cal L}_X \gpu{jo}
	\hspace{-.1cm} - \hspace{-.1cm}
	\nabla_j X_m {\cal L}_X \gpu{oi} \hspace{-.1cm} \Big)
	\hspace{-.05cm} q^1_{kl}(\gp)
	\hspace{-.1cm} \d \mu_\gp
\end{displaymath}
\begin{equation} \label{deg.basis.pm.i}
	+ \int \limits_\sur \gpo{ik} \gpo{jl} \gpo{mo}
	X_m {\cal L}_X \gpu{ij}
	\nabla_o q^r_{kl}(\gp) \d \mu_\gp.
\end{equation}
we obtain from (\ref{deg.sec}), (\ref{deg.sec.coeff})
and (\ref{deg.sec.coeff2})
\begin{displaymath}
	\langle \delta^2 \tpc_f(\vari) , e \rangle
	=  \langle d \tpc_\gp . q^1(\gp) , e \rangle
	\ I(\vari_0).
\end{displaymath}
As $\ I(\vari_{0,m}) \rightarrow I(\vari_0)
\mbox{ for } V_{0,m} \rightarrow \vari_0
\mbox{ in } W^{1,2}(\sur,\rel^n)\ $,
it suffices to find
$\ \vari_0 \mbox{ respectively } \vari \in W^{1,2}(\sur)\ $
such that
\begin{equation} \label{deg.basis.pm.equ}
	\pm I(\vari_0) > 0.
\end{equation}
Recalling $\ \delta \tpc_f . W \in \hat {\cal V}_f \perp e
\mbox{ for any } W \in C^\infty(\sur,\rel^n)\ $,
we see in the same way by (\ref{full.vari-form}) that
\begin{displaymath}
	0 = \langle \delta \tpc_f . W , e \rangle
	= - 2 \int \limits_\sur g^{ik} g^{jl}
	\langle A_{ij}^0 , W \rangle
	q^1_{kl}(\gp) \d \mu_g
	\ \langle d \tp_\gp . q^1(\gp) , e \rangle,
\end{displaymath}
hence, as $\ \langle d \tpc_\gp . q^1(\gp) , e \rangle > 0\ $
and $\ W \in C^\infty(\sur,\rel^n)\ $ is arbitrary,
$\ g^{ik} g^{jl} A_{ij}^0 q^1_{kl}(\gp) = 0\ $.
We get by (\ref{deg.basis.pm.princ}) that
$\ A_{11}^0 q^1_{11}(\gp) = 0
\mbox{ in } B_1(0) \cong U(x_0)\ $
and by (\ref{deg.basis.pm.nab}),
as $\ q^1_{11}(\gp) - i q^1_{12}(\gp)\ $
is holomorphic,
\begin{equation} \label{deg.basis.pm.q}
\left.
\begin{array}{c}
	q^1_{11}(\gp) = - q^1_{22}(\gp) = 0, \\
	q^1_{12} \equiv: q^1 \in \rel - \{0\}
\end{array}
\right\}
	\quad \mbox{in } B_1(0) \cong U(x_0)
\end{equation}
in local coordinates of $\ \varphi\ $.

Next we consider $\ \vari_0\ $ to be normal at $\ f\ $
and calculate by (\ref{full.metric-vari}) that
\begin{displaymath}
	(\partial_t g_{t,ij})_{|t=0}
	= \langle \partial_i f , \partial_j \vari \rangle
	+ \langle \partial_j f , \partial_i \vari \rangle =
\end{displaymath}
\begin{displaymath}
	= - 2 \langle A_{ij} , \vari_0 \rangle
	- \gamma_s \langle \partial_i f , \partial_j \vari_s \rangle
	- \gamma_s \langle \partial_j f , \partial_i \vari_s \rangle,
\end{displaymath}
hence by (\ref{deg.basis.pm.q})
\begin{displaymath}
	\parallel (\partial_t g_t)_{|t=0} \parallel_{L^2(\sur)}
	\leq 2 \parallel A \parallel_{L^2(supp \vari_0,g)}
	\parallel \vari_0 \parallel_{L^\infty(\sur)}
	+ C |\gamma| \leq
\end{displaymath}
\begin{displaymath}
	\leq C \parallel A \parallel_{L^2(supp \vari_0,g)}
	\parallel \vari_0 \parallel_{L^\infty(\sur)}.
\end{displaymath}
As above
we can select $\ \fu \in L^2(\sur), X \in \vf^1(\sur)\ $
in (\ref{deg.decomp}) such that
\begin{displaymath}
	\parallel \fu \parallel_{L^2(\sur)},
	\parallel X \parallel_{W^{1,2}(\sur)}
	\leq C \parallel A \parallel_{L^2(supp \vari_0,g)}
	\parallel \vari_0 \parallel_{L^\infty(\sur)},
\end{displaymath}
where $\ C\ $ does not depend on $\ \vari_0\ $.
We get from (\ref{deg.basis.pm.i})
\begin{equation} \label{deg.basis.pm.aux}
	\Big| I(\vari_0) - \int \limits_\sur \gpo{ik} \gpo{jl}
	e^{-2u} (\partial_{tt} g_{t,ij})_{|t=0}
	\ q^1_{kl}(\gp) \d \mu_\gp \Big|
	\leq C \parallel A \parallel_{L^2(supp \vari_0,g)}^2
	\parallel \vari_0 \parallel_{L^\infty(\sur)}^2,
\end{equation}
where $\ C\ $ does not depend on $\ \vari_0\ $.

We continue with (\ref{full.metric-vari}) and get,
as $\ supp\ \vari_0 \cap supp\ \vari_s = \emptyset
\mbox{ for } s = 1, \ldots, \dim\ \teich-1\ $,
\begin{displaymath}
	(\partial_{tt} g_{t,ij})_{|t=0}
	= 2 \langle \partial_i \vari , \partial_j \vari \rangle
	= 2 \langle \partial_i \vari_0 , \partial_j \vari_0 \rangle
	+ \gamma_r \gamma_s \langle \partial_i \vari_r,
	\partial_j \vari_s \rangle.
\end{displaymath}
As
\begin{displaymath}
	\Big| \int \limits_\sur \gpo{ik} \gpo{jl}
	e^{-2u} 2 \gamma_r \gamma_s
	\langle \partial_i \vari_r,
	\partial_j \vari_s \rangle
	q^1_{kl}(\gp) \d \mu_\gp \Big| \leq
\end{displaymath}
\begin{displaymath}
	\leq C |\gamma|^2
	\leq C \parallel A \parallel_{L^2(supp \vari_0,g)}^2
	\parallel \vari_0 \parallel_{L^\infty(\sur)}^2
\end{displaymath}
and
\begin{displaymath}
	\int \limits_\sur \gpo{ik} \gpo{jl}
	e^{-2u} 2 \langle \partial_i \vari_0,
	\partial_j \vari_0 \rangle
	q^1_{kl}(\gp) \d \mu_\gp
	= \int \limits_{B_1(0)} 4 q^1 e^{-2v}
	\langle \partial_1 \vari_0,
	\partial_2 \vari_0 \rangle \d \Lt,
\end{displaymath}
where we identify $\ U(x_0) \cong B_1(0)\ $,
we get from (\ref{deg.basis.pm.aux})
\begin{equation} \label{deg.basis.pm.aux2}
	\Big| I(\vari_0)
	- \int \limits_{B_1(0)} 4 q^1 e^{-2v}
	\langle \partial_1 \vari_0,
	\partial_2 \vari_0 \rangle \d \Lt \Big|
	\leq C \parallel A \parallel_{L^2(supp \vari_0,g)}^2
	\parallel \vari_0 \parallel_{L^\infty(\sur)}^2,
\end{equation}
where $\ C\ $ does not depend on $\ \vari_0\ $.

Perturbing $\ x_0 \cong 0 \mbox{ in }
U(x_0) \cong B_1(0)\ $ slightly,
we may assume that $\ 0\ $
is a Lebesgue point for $\ \nabla f
\mbox{ and } v\ $.
We select a vector $\ \nor_0 \in \rel^n
\mbox{ normal in } f(0) \mbox{ at } f\ $
and define via normal projection
$\ \nor := \pro_f^\perp \nor_0 \in W^{1,2}(B_1(0),\rel^n)\ $,
see (\ref{direct.pro}).
Clearly, $\ \nor(0) = \nor_0\ $,
and $\ 0\ $ is a Lebesgue point of $\ \nor\ $.
For $\ \eta \in C^\infty_0(B_1(0))\ $,
we put $\ \eta_\varrho(y) := \eta(\varrho^{-1} y)
\mbox{ and } \vari_0 := \eta_\varrho \nor\ $
and calculate
\begin{displaymath}
	\int \limits_{B_1(0)}
	\langle \partial_1 \vari_0 , \partial_2 \vari_0 \rangle
	e^{-2v} \d \Lt
	= \int \limits_{B_1(0)}
	\partial_1 \eta_\varrho \partial_2 \eta_\varrho
	|\nor|^2 e^{-2v} \hspace{-.1cm} \d \Lt
	+ \int \limits_{B_1(0)}
	(\partial_1 \eta_\varrho) \eta_\varrho
	\langle \nor , \partial_2 \nor \rangle
	e^{-2v} \hspace{-.1cm} \d \Lt +
\end{displaymath}
\begin{displaymath}
	+ \int \limits_{B_1(0)}
	\eta_\varrho \partial_2 \eta_\varrho
	\langle \partial_1 \nor , \nor \rangle
	e^{-2v} \d \Lt
	+ \int \limits_{B_1(0)}
	\eta_\varrho^2
	\langle \partial_1 \nor , \partial_2 \nor \rangle
	e^{-2v} \d \Lt.
\end{displaymath}
As $\ \parallel \nabla \eta_\varrho \parallel_{L^2(B_\varrho(0))}
= \parallel \nabla \eta \parallel_{L^2(B_1(0))}\ $,
we see
\begin{displaymath}
	\lim \limits_{\varrho \rightarrow 0}
	\int \limits_{B_1(0)}
	\langle \partial_1 \vari_0 , \partial_2 \vari_0 \rangle
	e^{-2v} \d \Lt =
\end{displaymath}
\begin{displaymath}
	= \lim \limits_{\varrho \rightarrow 0}
	\int \limits_{B_1(0)}
	\partial_1 \eta(y) \partial_2 \eta(y)
	|\nor(\varrho y)|^2 e^{-2v(\varrho y)} \d y
	= e^{-2 v(0)}
	\int \limits_{B_1(0)}
	\partial_1 \eta \partial_2 \eta \d \Lt.
\end{displaymath}
Since $\ \parallel \vari_0 \parallel_{L^\infty(\sur)}
\leq \parallel \eta_\varrho \parallel_{L^\infty(B_\varrho(0))}
= \parallel \eta \parallel_{L^\infty(B_1(0))}
\mbox{ and } \parallel A \parallel_{L^2(supp \vari_0,g)}
\rightarrow 0 \mbox{ for } \varrho \rightarrow 0\ $,
we get
\begin{displaymath}
	\lim \limits_{\varrho \rightarrow 0}
	I(\eta_\varrho \nor)
	= 4 q^1 e^{-2 v(0)}
	\int \limits_{B_1(0)}
	\partial_1 \eta \partial_2 \eta \d \Lt.
\end{displaymath}
Introducing new coordinates
$\ \tilde y_1 := (y_1 + y_2)/\sqrt{2},
\tilde y_2 := (-y_1 + y_2)/\sqrt{2}\ $,
that is we rotate the coordiantes by 45 degrees,
and putting $\ \eta(\tilde y_1,\tilde y_2)
= \xi(\tilde y_1) \tau(\tilde y_2)
\mbox{ with } \xi, \tau \in C^\infty_0(]-1/2,1/2[)\ $,
we see
\begin{displaymath}
	\int \limits_{B_1(0)}
	\partial_{y_1} \eta \partial_{y_2} \eta \d \Lt
	= \frac{1}{2} \int |\xi'|^2 \int |\tau|^2
	- \frac{1}{2} \int |\xi|^2 \int |\tau'|^2.
\end{displaymath}
Choosing $\ \tau \in C^\infty_0(]-1/2,1/2[), \tau \not\equiv 0
\mbox{ and } \xi(t) := \tau(2t)\ $,
we get $\ \int |\xi'|^2 = 2 \int |\tau'|^2,
2 \int |\xi|^2 = \int |\tau|^2\ $ and
\begin{displaymath}
	\int \limits_{B_1(0)}
	\partial_{y_1} \eta \partial_{y_2} \eta \d \Lt
	= \frac{3}{4} \int |\tau|^2 \int |\tau'|^2 > 0.
\end{displaymath}
Exchanging $\ \xi \mbox{ and } \tau\ $,
we produce a negative sign.
Choosing $\ \varrho\ $ small enough
and approximating $\ \vari_0 = \eta_\varrho \nor\ $
smoothly,
we obtain the desired $\ \vari_\pm \in C^\infty_0(\sur
- \{ x_1, \ldots, x_N\},\rel^n)\ $.
\proof
Next we extend the convegence criterion
in Proposition \ref{full.conv}
to the second variation.

\begin{proposition} \label{deg.conv}

Let $\ f: \sur \rightarrow \rel^n\ $ be a weak local bilipschitz immersion
approximated by smooth immersions $\ f_m\ $
with pull-back metrics $\ g = f^* \geu = e^{2u} \gp,
g_m = f_m^* \geu = e^{2u_m} \gp_m\ $
for some smooth unit volume constant
curvature metrics $\ \gp, \gpm\ $
and satisfying
\begin{equation} \label{deg.conv.ass}
\begin{array}{c}
	f_m \rightarrow f
	\quad \mbox{weakly in } W^{2,2}(\sur),
	\mbox{weakly}^* \mbox{ in } W^{1,\infty}(\sur), \\
	\Lambda^{-1} \gp \leq g_m \leq \Lambda \gp, \\
	\parallel u_m \parallel_{L^\infty(\sur)} \leq \Lambda
\end{array}
\end{equation}
for some $\ \Lambda < \infty\ $.
Then for any chart $\ \psi: U(\tp(\gp)) \subseteq \teich
\rightarrow \rel^{\dim \teich}\ $,
$\ \tpc := \psi \circ \tp,
\delta \tpc = d \psi \circ \delta \tp,
\delta^2 \tpc\ $ defined in (\ref{deg.sec-def})
and any $\ W \in W^{1,2}(\sur,\rel^n)\ $
\begin{displaymath}
\begin{array}{c}
	\delta \tpc_{f_m} . W \rightarrow \delta \tpc_f . W,
	\quad
	\delta^2 \tpc_{f_m}(W) \rightarrow \delta^2 \tpc_f(W)
\end{array}
\end{displaymath}
as $\ m \rightarrow \infty\ $.
\end{proposition}
{\pr Proof:} \\
By Proposition \ref{full.conv},
we know already $\ \delta \tpc_{f_m} . W \rightarrow
\delta \tpc_f . W\ $
and get further (\ref{full.conv.met-conv}) and (\ref{full.conv.conv}).

Next we select a slice $\ {\cal S}(\gp)\ $ of unit volume
constant curvature metrics for $\ \tp(\gp) =: \tau_0 \in \teich\ $
around $\ \gp\ $
with $\ \tp: {\cal S}(\gp) \cong U(\tau_0)\ $,
and $\ q^r(\tgp) \in S_2^{TT}(\tgp)
\mbox{ for } \tp(\tgp) \in U(\tau_0)\ $,
see \bcite{fisch.trom.conf}, \bcite{trom.teich}.

As $\ \tp(\gpm) = \tp(g_m) \rightarrow \tp(\gp) \in U(\tau_0)
\cong {\cal S}(\gp)\ $,
there exist for $\ m\ $ large enough
smooth diffeomorphisms $\ \diffeo_m
\mbox{ of } \sur\ $ homotopic to the identity
with $\ \diffeo_m^* \gpm =: \tgpm \in {\cal S}(\gp)\ $.
As $\ \tp(\tgpm) = \tp(\gpm) = \tp(g_m)
\rightarrow \tp(\gp) \mbox{ and } \tgpm \in {\cal S}(\gp)\ $,
we get
\begin{equation} \label{deg.conv.smooth}
	\tgpm \rightarrow \gp
	\quad \mbox{ smoothly}.
\end{equation}
The Theorem of Ebin and Palais,
see \bcite{fisch.trom.conf}, \bcite{trom.teich},
and the remarks following
imply by (\ref{full.conv.met-conv}), (\ref{deg.conv.ass})
and (\ref{deg.conv.smooth})
that after appropriately modifying $\ \diffeo_m\ $
\begin{equation} \label{deg.conv.diffeo-conv}
	\diffeo_m, \diffeo_m^{-1} \rightarrow id_\sur
	\quad \mbox{weakly in } W^{2,2}(\sur),
	\mbox{weakly}^* \mbox{ in } W^{1,\infty}(\sur),
\end{equation}
in particular,
\begin{equation} \label{deg.conv.diffeo-bound}
	\parallel D \diffeo_m
	\parallel_{W^{1,2}(\sur) \cap L^\infty(\sur)},
	\parallel D (\diffeo_m^{-1})
	\parallel_{W^{1,2}(\sur) \cap L^\infty(\sur)}
	\leq C
\end{equation}
with $\ C\ $ independent of $\ m \mbox{ and } \vari\ $.

For the second variations,
we see by (\ref{full.metric-vari})
and (\ref{deg.conv.ass})
\begin{displaymath}
	\partial_t ((f_m + t W)^*
	\geu)_{ij}\ _{|t = 0}
	= \langle \partial_i f_m,\partial_j W \rangle
	+ \langle \partial_j f_m,\partial_i W \rangle
	\rightarrow
\end{displaymath}
\begin{displaymath}
	\rightarrow \langle \partial_i f,\partial_j W \rangle
	+ \langle \partial_j f,\partial_i W \rangle
	= \partial_t ((f + t W)^*
	\geu)_{ij}\ _{|t = 0}
	\quad \mbox{weakly in } W^{1,2}(\sur).
\end{displaymath}
Following (\ref{deg.decomp}), we decompose
\begin{displaymath}
	\diffeo_m^* \partial_t ((\tilde f_m + t W)^*
	\geu)_{|t = 0}
	= \fu_m \tgpm + {\cal L}_{X_m} \tgpm
	+ q_m
\end{displaymath}
with $\ q_m \in S_2^{TT}(\tgpm)\ $
and moreover recalling (\ref{deg.conv.diffeo-bound}),
we can achieve
\begin{displaymath}
	\parallel \fu_m \parallel_{W^{1,2}(\sur)},
	\parallel X_m \parallel_{W^{2,2}(\sur)},
	\parallel q_m \parallel_{C^2(\sur)}
	\leq C.
\end{displaymath}
For a subsequence, we get $\ \fu_m \rightarrow \fu
\mbox{ weakly in } W^{1,2}(\sur),
X_m \rightarrow X \mbox{ weakly in } W^{2,2}(\sur), \\
q_m \rightarrow q \mbox{ strongly in } C^1(\sur)\ $
with $\ q \in S_2^{TT}(\gp)\ $ and by (\ref{deg.conv.diffeo-conv})
\begin{equation} \label{deg.conv.decomp}
	\partial_t ((f+ t W)^* \geu)_{|t = 0}
	= \fu \gp + {\cal L}_X \gp + q.
\end{equation}
Observing from (\ref{full.metric-vari}) that
\begin{displaymath}
	(\partial_{tt} g_{t,ij})_{|t = 0}
	= 2 \langle \partial_i W , \partial_j W \rangle,
\end{displaymath}
we get from (\ref{deg.sec}), (\ref{deg.sec.coeff})
and the equivariance of $\ \tpc\ $
\begin{displaymath}
	\delta^2 \tpc_{f_m}(W) =
\end{displaymath}
\begin{displaymath}
	= \sum \limits_{r=1}^{\dim\ \teich}
	\alpha_{m,r}\ d \tpc_\tgpm . q^r(\tgpm)
	+ \sum \limits_{r,s=1}^{\dim\ \teich}
	\beta_{m,r} \beta_{m,s}
	\ d^2 \tpc_\tgpm(q^r(\tgpm),q^s(\tgpm)),
\end{displaymath}
where
\begin{displaymath}
\begin{array}{c}
	\alpha_{m,r}
	:= \int \limits_\sur g_m^{ik} g_m^{jl}
	2 \langle \partial_i W , \partial_j W \rangle
	q^r_{kl}(\gpm) \d \mu_{g_m} + \\
	- \int \limits_\sur \tgpmo{ik} \tgpmo{jl}
	\Big( {\cal L}_{X_m} {\cal L}_{X_m} \tgpm
	+ 2 \fu_m {\cal L}_{X_m} \tgpm \Big)_{ij}
	q^r_{kl}(\tgpm) \d \mu_\tgpm + \\
	- \int \limits_\sur \tgpmo{ik} \tgpmo{jl}
	\Big( 2 \fu_m q_m
	+ 2 {\cal L}_{X_m} q_m \Big)_{ij}
	q^r_{kl}(\tgpm) \d \mu_\tgpm,
\end{array}
\end{displaymath}
\begin{displaymath}
	\beta_{m,r}
	:= \int \limits_\sur \tgpmo{ik} \tgpmo{jl}
	q_{m,ij} q^r_{kl}(\tgpm) \d \mu_\tgpm.
\end{displaymath}
By the above convergences
in particular by (\ref{full.conv.met-conv}),
(\ref{full.conv.conv}) and (\ref{deg.conv.smooth}),
we obtain
\begin{displaymath}
\begin{array}{c}
	\alpha_{m,r} \rightarrow \alpha_r
	:= \int \limits_\sur g^{ik} g^{jl}
	2 \langle \partial_i W , \partial_j W \rangle
	q^r_{kl}(\gp) \d \mu_{g} + \\
	- \int \limits_\sur \gpo{ik} \gpo{jl}
	\Big( {\cal L}_{X} {\cal L}_{X} \gp
	+ 2 \fu {\cal L}_{X} \gp \Big)_{ij}
	q^r_{kl}(\gp) \d \mu_\gp + \\
	- \int \limits_\sur \gpo{ik} \gpo{jl}
	\Big( 2 \fu q
	+ 2 {\cal L}_{X} q \Big)_{ij}
	q^r_{kl}(\gp) \d \mu_\gp,
\end{array}
\end{displaymath}
\begin{displaymath}
	\beta_{m,r} \rightarrow \beta_r
	:= \int \limits_\sur \gpo{ik} \gpo{jl}
	q_{ij} q^r_{kl}(\gp) \d \mu_\gp.
\end{displaymath}
We see from (\ref{deg.conv.smooth})
\begin{displaymath}
	d^2 \tpc_\tgpm(q^r(\tgpm),q^s(\tgpm))
	\rightarrow d^2 \tpc_\gp(q^r(\gp),q^s(\gp)).
\end{displaymath}
Observing (\ref{deg.conv.decomp}),
we get from (\ref{deg.sec}), (\ref{deg.sec.coeff})
\begin{equation} \label{deg.conv.sec-conv}
	\delta^2 \tpc_{\tilde f_m}(W)
	\rightarrow \delta^2 \tpc_f(W)
	\quad \mbox{for any } W \in W^{1,2}(\sur,\rel^n),
\end{equation}
and the proposition is proved.
\proof
{\large \bf Remark:} \\
The bound on the conformal factor $\ u_m\ $
in the above proposition
is implied by Proposition \ref{conf.strong},
when we replace the weak convergence
of $\ f_m \rightarrow f \mbox{ in } W^{2,2}(\sur)\ $
by strong convergence.
\defin \\
Now we can extend the correction lemma \ref{full.corr}
to the degenerate case.

\begin{lemma} \label{deg.corr}

Let $\ f: \sur \rightarrow \rel^n\ $ be a weak local bilipschitz immersion
approximated by smooth immersions $\ f_m\ $
satisfying (\ref{direct.prop.teich}) - (\ref{direct.prop.imm}).

If $\ f\ $ is not of full rank in Teichm\"uller space,
then for arbitrary $\ x_0 \in \sur,
\mbox{ neighbourhood } U_\mark(x_0)
\subseteq \sur \mbox{ of } x_0,
\mbox{ and } \Lambda < \infty\ $,
there exists a neighbourhood $\ U(x_0) \subseteq U_\mark(x_0)
\mbox{ of } x_0\ $,
variations $\ \vari_1, \ldots, \vari_{\dim \teich - 1},
\vari_\pm \in C^\infty_0(\sur - \overline{U(x_0)},\rel^n)\ $,
satisfying (\ref{deg.basis.span}) and (\ref{deg.basis.pm}),
and $\ \delta > 0, C < \infty, m_0 \in \nat\ $
such that for any $\ \vari \in C^\infty_0(U(x_0),\rel^n)\ $
with $\ f_m + \vari\ $ a smooth immersion
for some $\ m \geq m_0\ $, and $\ \vari = 0\ $ or
\begin{equation} \label{deg.corr.ass}
\begin{array}{c}
	\Lambda^{-1} \gp \leq (f_m + \vari)^* \geu \leq \Lambda \gp, \\
	\parallel \vari \parallel_{W^{2,2}(\sur)}
	\leq \Lambda, \\
	\int \limits_{U_\mark(x_0)} |A_{f_m + \vari}|^2 \d \mu_{f_m + \vari}
	\leq \varepsilon_0(n),
\end{array}
\end{equation}
where $\ \varepsilon_0(n)\ $ is as in Lemma \ref{conf.lemma},
and any $\ \tau \in \teich\ $ with
\begin{equation} \label{deg.corr.teich-point}
	d_\teich(\tau,\tau_0) \leq \delta,
\end{equation}
there exists $\ \lambda \in \rel^{\dim \teich-1},
\mu_\pm \in \rel\ $, satisfying $\ \mu_+ \mu_- = 0\ $,
\begin{displaymath}
	\tp( (f_m + \vari
	+ 
	\lambda_r \vari_r + \mu_\pm \vari_\pm)^* \geu)
	= \tau,
\end{displaymath}
\begin{displaymath}
	|\lambda|, |\mu_\pm|
	\leq C d_\teich \Big( \tp( (f_m + \vari)^* \geu),
	\tau \Big)^{1/2}.
\end{displaymath}
Further for any $\ \lambda_0 > 0\ $,
one can choose $\ m_0, \delta\ $ in such a way
that for $\ m \geq m_0\ $,
\begin{displaymath}
	\parallel \vari \parallel_{W^{2,2}(\sur)} \leq \delta,
\end{displaymath}
there exists $\ \tilde \lambda \in \rel^{\dim \teich-1},
\tilde \mu_\pm \in \rel\ $, satisfying $\ \tilde \mu_+ \tilde \mu_- = 0\ $,
\begin{displaymath}
\begin{array}{c}
	\tp( (f_m + \vari
	+ 
	\tilde \lambda_r \vari_r + \tilde \mu_\pm \vari_\pm)^* \geu)
	= \tau, \\
	|\tilde \lambda|, |\tilde \mu_\pm|
	\leq \lambda_0,
\end{array}
\end{displaymath}
and
\begin{displaymath}
\begin{array}{c}
	\mu_+ \tilde \mu_+ \leq 0, \tilde \mu_- = 0,
	\quad \mbox{if } \mu_+ \neq 0, \\
	\mu_- \tilde \mu_- \leq 0, \tilde \mu_+ = 0,
	\quad \mbox{if } \mu_- \neq 0.
\end{array}
\end{displaymath}
\end{lemma}
{\pr Proof:} \\
By (\ref{direct.prop.w22}), (\ref{direct.prop.met}) and $\ \Lambda\ $
large enough, we may assume
\begin{equation} \label{deg.corr.d2}
\begin{array}{c}
	\parallel u_m, D f_m \parallel_{W^{1,2}(\sur)
	\cap L^\infty(\sur)} \leq \Lambda, \\
	\Lambda^{-1} \gp \leq f_m^* \geu \leq \Lambda \gp,
\end{array}
\end{equation}
in particular
\begin{equation} \label{deg.corr.fund}
	\int \limits_\sur |A_{f_m}|^2 \d \mu_{f_m}
	\leq C(\sur,\gp,\Lambda).
\end{equation}
Putting $\ \funda_m := |\nabla^2_\gp f_m|_\gp^2 \mu_{\gp}\ $,
we see $\ \funda_m(\sur) \leq C(\Lambda,\gp) $
and for a subsequence $\ \funda_m \rightarrow \funda\ $ 
weakly$^* \mbox{ in } C_0^0(\sur)^*\ $.
Clearly $\ \funda(\sur) < \infty\ $,
and there are at most finitely many $\ y_1, \ldots, y_N \in \sur
\mbox{ with } \funda(\{y_i\}) \geq \varepsilon_1\ $,
where we choose $\ \varepsilon_1 > 0\ $ below.

By Proposition \ref{deg.basis},
we can select $\ \vari_1, \ldots, \vari_{\dim \teich - 1},
\vari_+ = \vari_{\dim \teich},
\vari_- = \vari_{\dim \teich + 1}
\in C^\infty_0(\sur - \{x_0, y_1, \ldots, y_N\},\rel^n)\ $
satisfying (\ref{deg.basis.span}) and (\ref{deg.basis.pm}).

We choose a neighbourhood $\ U_0(x_0)
\subseteq U_\mark(x_0) \mbox{ of } x_0\ $
with a chart $\ \varphi_0: U_0(x_0)
\stackrel{\approx}{\longrightarrow} B_2(0),
\varphi_0(x_0) = 0\ $,
\begin{displaymath}
	supp\ \vari_r \cap \overline{U_0(x_0)} = \emptyset
	\quad \mbox{for } r = 1, \ldots, \dim \teich + 1,
\end{displaymath}
put $\ x_0 \in U_\varrho(x_0) = \varphi_0^{-1}(B_\varrho(0))
\mbox{ for } 0 < \varrho \leq 2\ $
and choose $\ x_0 \in U(x_0) \subseteq U_{1/2}(x_0)\ $
small enough, as we will see below.

Next for any $\ x \in \cup_{r=1}^{\dim \teich + 1} supp\ \vari_r\ $,
there exists a neighbourhood $\ U_0(x) \mbox{ of } x\ $
with a chart $\ \varphi_x: U_0(x)
\stackrel{\approx}{\longrightarrow} B_2(0),
\varphi_x(x) = 0,
\overline{U_0(x)} \cap \overline{U_0(x_0)} = \emptyset,
\funda(U_0(x)) < \varepsilon_1\ $ and
in the coordinates of the chart $\ \varphi_x\ $
\begin{equation} \label{deg.corr.christ}
	\int \limits_{U_0(x)} \gpo{ik} \gpo{jl} \gpu{rs}
	\Gamma_{\gp,ij}^r \Gamma_{\gp,kl}^s \d \mu_\gp
	\leq \varepsilon_1.
\end{equation}
Putting $\ x \in U_\varrho(x) = \varphi_x^{-1}(B_\varrho(0))
\subset \subset U_0(x) \mbox{ for } 0 < \varrho \leq 2\ $,
we see that there are finitely many $\ x_1, \ldots, x_\M
\in \cup_{r=1}^{\dim \teich + 1} supp\ \vari_r\ $
such that
\begin{displaymath}
	\cup_{r=1}^{\dim \teich + 1} supp\ \vari_r
	\subseteq \cup_{k=1}^\M U_{1/2}(x_k).
\end{displaymath}
Then there exists $\ m_0 \in \nat\ $
such that for $\ m \geq m_0\ $
\begin{displaymath}
	\int \limits_{U_1(x_k)} |\nabla^2_\gp f_m|^2_\gp
	\d \mu_{\gp} < \varepsilon_1
	\quad \mbox{for } k = 1, \ldots, \M.	
\end{displaymath}
For $\ \vari \mbox{ and } m \geq m_0\ $ as above,
we put $\ \tilde f_{m,\lambda,\mu} := f_m + \vari + \lambda_r \vari_r
+ \mu_\pm \vari_\pm\ $. Clearly
\begin{displaymath}
	supp(f_m - \tilde f_{m,\lambda,\mu})
	\subseteq \cup_{k=0}^\M U_{1/2}(x_k).
\end{displaymath}
By (\ref{deg.corr.ass}), (\ref{deg.corr.d2}), (\ref{deg.corr.christ}),
and $\ |\lambda|, |\mu| < \lambda_0 \leq 1/4\ $ small enough
independent of $\ m \mbox{ and } \vari\ $,
$\ \tilde f_{m,\lambda,\mu}\ $ is a smooth immersion with
\begin{equation} \label{deg.corr.met}
	(2 \Lambda)^{-1} \gp \leq \tilde g_{m,\lambda,\mu}
	:= \tilde f_{m,\lambda,\mu}^* \geu
	\leq 2 \Lambda \gp,
\end{equation}
and if $\ \vari \neq 0\ $ by (\ref{deg.corr.ass})
and the choice of $\ U_0(x)\ $ that
\begin{displaymath}
	\int \limits_{U_1(x_k)} |A_{\tilde f_{m,\lambda,\mu}}|^2
	\d \mu_{\tilde g_{m,\lambda,\mu}} \leq \varepsilon_0(n)
	\quad \mbox{for } k = 0, \ldots, \M
\end{displaymath}
for $\ C(\Lambda,\gp) (\varepsilon_1 + \lambda_0) \leq \varepsilon_0(n)\ $.
If $\ \vari = 0\ $, we see
$\ supp(f_m - \tilde f_{m,\lambda,\mu})
\subseteq \cup_{k=1}^\M U_{1/2}(x_k)\ $.
Further by (\ref{deg.corr.fund})
\begin{displaymath}
	\int \limits_\sur |K_{\tilde g_{m,\lambda,\mu}}|
	\d \mu_{\tilde g_{m,\lambda,\mu}}
	\leq \frac{1}{2} \int \limits_\sur |A_{\tilde f_{m,\lambda,\mu}}|^2
	\d \mu_{\tilde g_{m,\lambda,\mu}} \leq
\end{displaymath}
\begin{displaymath}
	\leq \frac{1}{2} \int \limits_\sur |A_{f_m}|^2 \d \mu_{f_m}
	+  \frac{1}{2} \sum \limits_{k=0}^\M \int \limits_{U_{1/2}(x_k)}
	|A_{\tilde f_{m,\lambda,\mu}}|^2 \d \mu_{\tilde g_{m,\lambda,\mu}}
	\leq C(\sur,\gp,\Lambda) + (\M + 1) \varepsilon_0(n).
\end{displaymath}
This verifies (\ref{conf.lemma.metric2})
and (\ref{conf.lemma.fund})
for $\ f = f_m, \tilde f = \tilde f_{m,\lambda,\mu}, g_0 = \gp\ $
and different, but appropriate $\ \Lambda\ $.
(\ref{conf.lemma.ub}) follows from (\ref{direct.prop.met})
and (\ref{deg.corr.fund}).
Then for the unit volume constant curvature metric
$\ \tgpmlm = e^{-2 \tilde u_{m,\lambda,\mu}} \tilde g_{m,\lambda,\mu}\ $
conformal to $\ \tilde g_{m,\lambda,\mu}\ $
by Poincar\'e's Theorem,
see \bcite{fisch.trom.conf}, \bcite{trom.teich},
we get from Lemma \ref{conf.lemma} that
\begin{equation} \label{deg.corr.conf}
	\parallel \tilde u_{m,\lambda,\mu} \parallel_{L^\infty(\sur)},
	\parallel \nabla \tilde u_{m,\lambda,\mu}
	\parallel_{L^2(\sur,\gp)} \leq C
\end{equation}
with $\ C < \infty\ $ independent of $\ m \mbox{ and } \vari\ $.

From (\ref{deg.corr.ass}),
we have a $\ W^{2,2} \cap W^{1,\infty}-$bound on $\ \tilde f_{m,0,0}\ $,
hence for $\ \tilde f_{m,\lambda,\mu}\ $.
On $\ \sur - U(x_0)\ $, we get
$\ \tilde f_{m,\lambda,\mu} = f_m + \lambda_r \vari_r
+ \mu_\pm \vari_\pm \rightarrow f
\mbox{ weakly in } W^{2,2}(\sur - U(x_0))
\mbox{ and weakly}^* \mbox{ in } W^{1,\infty}(\sur - U(x_0))
\mbox{ for } m_0 \rightarrow \infty, \lambda_0 \rightarrow 0\ $
by (\ref{direct.prop.conv}).
If $\ \vari = 0\ $,
then $\ \tilde f_{m,\lambda,\mu} = f_m + \lambda_r \vari_r
+ \mu_\pm \vari_\pm \rightarrow f
\mbox{ weakly in } W^{2,2}(\sur)
\mbox{ and weakly}^* \mbox{ in } W^{1,\infty}(\sur)
\mbox{ for } m_0 \rightarrow \infty, \lambda_0 \rightarrow 0\ $
by (\ref{direct.prop.conv}).
Hence letting $\ m_0 \rightarrow \infty,
\lambda_0 \rightarrow 0, U(x_0) \rightarrow \{x_0\}\ $,
we conclude
\begin{equation} \label{deg.corr.f-conv}
	\tilde f_{m,\lambda,\mu} \rightarrow f
	\quad \mbox{weakly in } W^{2,2}(\sur),
	\mbox{weakly}^* \mbox{ in } W^{1,\infty}(\sur),
\end{equation}
Then by (\ref{deg.corr.met}),
(\ref{deg.corr.conf}), (\ref{deg.corr.f-conv})
and Proposition \ref{deg.conv} for any $\ W \in W^{1,2}(\sur,\rel^n)\ $
\begin{equation} \label{deg.corr.vari-conv}
\begin{array}{c}
	\delta \tpc_{\tilde f_{m,\lambda,\mu}} . W \rightarrow \delta \tpc_f . W,
	\quad
	\delta^2 \tpc_{\tilde f_{m,\lambda,\mu}}(W) \rightarrow \delta^2 \tpc_f(W)
\end{array}
\end{equation}
for $\ m_0 \rightarrow \infty,
\lambda_0 \rightarrow 0, U(x_0) \rightarrow \{x_0\}\ $.

Further from (\ref{deg.corr.f-conv})
\begin{displaymath}
	\tilde g_{m,\lambda,\mu} \rightarrow f^* \geu = e^{2u} \gp =: g
	\left\{
	\begin{array}{c}
		\mbox{weakly in } W^{1,2}(\sur), \\
		\mbox{weakly}^*
		\mbox{ in } L^\infty(\sur),
	\end{array}
	\right.
\end{displaymath}
and by (\ref{deg.corr.conf}),
we see $\ \tilde u_{m,\lambda,\mu} \rightarrow \tilde u
\mbox{ weakly in } W^{1,2}(\sur)
\mbox{ and weakly}^* \mbox{ in } L^\infty(\sur)\ $,
hence
\begin{displaymath}
	\tgpmlm = e^{-2 \tilde u_{m,\lambda,\mu}} \tilde g_{m,\lambda,\mu}
	\rightarrow e^{2(u-\tilde u)} \gp
	\left\{
	\begin{array}{c}
		\mbox{weakly in } W^{1,2}(\sur), \\
		\mbox{weakly}^*
		\mbox{ in } L^\infty(\sur),
	\end{array}
	\right.
\end{displaymath}
which together with (\ref{deg.corr.met}) implies
\begin{equation} \label{deg.corr.teich-conv}
	\tp(\tilde g_{m,\lambda,\mu}) \rightarrow \tau_0.
\end{equation}
We select a chart $\ \psi: U(\tp(\gp)) \subseteq \teich
\rightarrow \rel^{\dim \teich}\ $
and put $\ \tpc := \psi \circ \tp,
\delta \tpc = d \psi \circ d \tp,
\hat {\cal V}_f := d \psi_{\tp \gp} . {\cal V}_f,
\delta^2 \tpc\ $ defined in (\ref{deg.sec-def}).
By (\ref{deg.corr.teich-conv}) for $\ m_0\ $ large enough,
$\ \lambda_0 \mbox{ and } U(x_0)\ $ small enough independent of $\ \vari\ $,
we get $\ \tp(\tilde g_{m,\lambda,\mu}) \in U(\tau_0)\ $
and define
\begin{displaymath}
        \Phi_m(\lambda,\mu) := \tpc( \tilde g_{m,\lambda,\mu}).
\end{displaymath}
We see by (\ref{deg.corr.vari-conv})
\begin{displaymath}
	D \Phi_m(\lambda,\mu)
	= (\delta \tpc_{\tilde f_{m,\lambda,\mu}} . \vari_r)
	_{r = 1, \ldots, \dim \teich + 1}
	\rightarrow (\delta \tpc_f . \vari_r)_{r = 1, \ldots, \dim \teich + 1}.
\end{displaymath}
After a change of coordinates, we may assume
\begin{displaymath}
	\hat {\cal V}_f = \rel^{\dim \teich - 1} \times \{0\}
\end{displaymath}
and $\ e := e_{\dim \teich} \perp \hat {\cal V}_f\ $.
Writing $\ \Phi_m(\lambda,\mu) = (\Phi_{m,0}(\lambda,\mu),\varphi_m(\lambda,\mu))
\in \rel^{\dim \teich - 1} \times \rel \times \rel\ $,
we get
\begin{displaymath}
	\partial_{\lambda} \Phi_{m,0}(\lambda,\mu)
	\rightarrow (\pro_{\hat {\cal V}_f}
	\delta \tpc_f . \vari_r)_{r = 1, \ldots, \dim \teich-1}
	=: A \in \rel^{(\dim \teich-1) \times (\dim \teich-1)}.
\end{displaymath}
From (\ref{deg.basis.span}),
we see that $\ A\ $ is invertible,
hence after a further change of coordiantes
we may assume that $\ A = I_{(\dim \teich - 1)}\ $
and
\begin{equation} \label{deg.corr.ana.inv}
	\parallel \partial_{\lambda} \Phi_{m,0}(\lambda,\mu)
	- I_{(\dim \teich - 1)} \parallel \leq 1/2
\end{equation}
for $\ m_0\ $ large enough,
$\ \lambda_0 \mbox{ and } U(x_0)\ $ small enough independent of $\ \vari\ $.
Next by (\ref{deg.basis.pm}), we obtain
\begin{equation} \label{deg.corr.vari1-conv}
\begin{array}{c}
	\partial_{\mu_\pm} \Phi_m(\lambda,\mu)
	\rightarrow \delta \tpc_f . \vari_\pm = 0, \\
	\nabla \varphi_m(\lambda,\mu) \rightarrow
	\langle (\delta \tpc_f . \vari_r)_{r = 1,
	\ldots, \dim \teich+1} , e \rangle = 0,
\end{array}
\end{equation}
as $\ \delta \tpc_f . \vari_r \in \hat {\cal V}_f \perp e\ $,
hence
\begin{equation} \label{deg.corr.ana-eps}
\begin{array}{c}
	|\partial_{\mu_\pm} \Phi_m(\lambda,\mu)|,
	|\nabla \varphi_m(\lambda,\mu)|
	\leq \varepsilon_2
\end{array}
\end{equation}
for any $\ \varepsilon_2 > 0\ $ chosen below,
if $\ m_0\ $ large enough,
$\ \lambda_0 \mbox{ and } U(x_0)\ $ small enough independent of $\ \vari\ $.

The second derivatives 
\begin{equation} \label{deg.conv.sec}
\begin{array}{c}
	\partial_{ss} \Phi_m(\lambda,\mu)
	= \delta^2 \tpc_{\tilde f_{m,\lambda,\mu}}(V_s), \\
	4 \partial_{sr} \Phi_m(\lambda,\mu)
	= \delta^2 \tpc_{\tilde f_{m,\lambda,\mu}}(V_s + V_r)
	- \delta^2 \tpc_{\tilde f_{m,\lambda,\mu}}(V_s - V_r)
\end{array}
\end{equation}
are given by the second variation in Teichm\"uller space.
From (\ref{deg.corr.vari-conv})
for $\ W = V_s, V_s \pm V_r \in C^\infty(\sur,\rel^n)\ $,
we conclude
for $\ m_0\ $ large enough,
$\ \lambda_0 \mbox{ and } U(x_0)\ $ small enough independent of $\ \vari\ $
that
\begin{equation} \label{deg.corr.ana.sec-bound}
	|D^2 \Phi_m(\lambda,\mu)| \leq \Lambda_1
\end{equation}
for some $\ 1 \leq \Lambda_1 < \infty\ $ and
\begin{equation} \label{deg.corr.sec-pos}
\begin{array}{c}
	\partial_{\mu_+ \mu_+} \Phi_m(\lambda,\mu)
	\rightarrow \delta^2 \tpc_f(\vari_+), \\
	\partial_{\mu_- \mu_-} \Phi_m(\lambda,\mu)
	\rightarrow \delta^2 \tpc_f(\vari_-),
\end{array}
\end{equation}
hence by (\ref{deg.basis.pm})
\begin{displaymath}
	\pm \lim \limits_{m,\lambda,\mu} \partial_{\mu_\pm \mu_\pm}
	\varphi_m(\lambda,\mu)
	= \pm \langle \delta^2 \tpc_f(\vari_\pm) , e \rangle > 0
\end{displaymath}
and
\begin{equation} \label{deg.corr.ana-pos}
\begin{array}{c}
	\pm \partial_{\mu_\pm \mu_\pm} \varphi_m(\lambda,\mu)
	\geq \gamma
\end{array}
\end{equation}
for some $\ 0 < \gamma \leq 1/4\ $
and $\ m_0\ $ large enough,
$\ \lambda_0 \mbox{ and } U(x_0)\ $ small enough independent of $\ \vari\ $.
Now, we choose $\ \varepsilon_2 > 0\ $ to satisfy
$\ C \Lambda_1 \varepsilon_2 \leq \gamma/4\ $.
Choosing $\ m_0\ $ even larger and $\ U(x_0)\ $ even smaller
we get by (\ref{deg.corr.teich-conv})
\begin{displaymath}
	d_\teich(\tp(\tilde g_0),\tau_0) \leq \delta,
\end{displaymath}
where we choose $\ \delta\ $ now. As
\begin{displaymath}
	d_\teich \Big( \tp( (f_m + \vari)^* \geu), \tau \Big)
	\leq d_\teich(\tp(\tilde g_0),\tau_0)
	+ d_\teich ( \tau_0,\tau )
	< 2 \delta
\end{displaymath}
by (\ref{deg.corr.teich-point}),
we choose $\ C_\psi \delta \leq \Lambda_1 \lambda_0^2, \lambda_0/8,
\gamma \lambda_0^2/32\ $
and conclude from Proposition \ref{ana.prop}
that there exists $\ \lambda \in \rel^{\dim \teich - 1},
\mu_\pm \in \rel \mbox{ with } \mu_+ \mu_- = 0\ $
and satisfying
\begin{displaymath}
\begin{array}{c}
	\tp( (f_m + \vari + \lambda_r \vari_r
	+ \mu_\pm \vari_\pm)^* \geu)
	= \tp(\tilde g_{m,\lambda,\mu})
	= \tau, \\
	|\lambda|, |\mu_\pm|
	\leq C d_\teich \Big( \tp( (f_m + \vari)^* \geu),
	\tau \Big)^{1/2}.
\end{array}
\end{displaymath}
To obtain the second conclusion,
we consider $\ \lambda_0 > 0\ $ such small
that $\ C \Lambda_1 \varepsilon_2 + C \Lambda_1 \lambda_0 \leq \gamma/2\ $
and fix this $\ \lambda_0\ $.
We assume $\ \parallel V \parallel_{W^{2,2}(\sur)} \leq \delta\ $
and see as in (\ref{deg.corr.f-conv})
that $\ \tilde f_{m,0,0} = f_m + \vari \rightarrow f
\mbox{ weakly in } W^{2,2}(\sur)
\mbox{ and weakly}^* \mbox{ in } W^{1,\infty}(\sur)
\mbox{ for } m_0 \rightarrow \infty, \delta \rightarrow 0\ $
by (\ref{direct.prop.conv}).
Again we get (\ref{deg.corr.vari-conv})
and (\ref{deg.corr.vari1-conv}) for $\ \lambda,\mu = 0\ $,
hence for $\ m_0\ $ large enough, $\ \delta \ $ small enough,
but fixed $\ U(x_0)\ $,
\begin{displaymath}
	|\nabla \varphi_m(0)| \leq \sigma
\end{displaymath}
with  $\ C \Lambda_1 \varepsilon_2 + C \sigma \lambda_0^{-1}
+ C \Lambda_1 \lambda_0 \leq \gamma\ $.
Then by Proposition \ref{ana.prop},
there exist further $\ \tilde \lambda \in \rel^{\dim \teich - 1},
\tilde \mu_\pm \in \rel \mbox{ with } \tilde \mu_+ \tilde \mu_- = 0\ $
and satisfying
\begin{displaymath}
\begin{array}{c}
	\tp( (f_m + \vari + \tilde \lambda_r \vari_r
	+ \tilde \mu_\pm \vari_\pm)^* \geu)
	= \tp(\tilde g_{\tilde \lambda,\tilde \mu})
	= \tau, \\
	|\tilde \lambda|, |\tilde \mu_\pm| \leq \lambda_0, \\
	\mu_+ \tilde \mu_+ \leq 0, \tilde \mu_- = 0,
	\quad \mbox{if } \mu_- = 0, \\
	\mu_- \tilde \mu_- \leq 0, \tilde \mu_+ = 0,
	\quad \mbox{if } \mu_+ = 0,
\end{array}
\end{displaymath}
and the lemma is proved.
\proof


\setcounter{equation}{0}

\section{Elementary properties of $\ \mini\ $}
\label{ele}

As a first application of our correction
Lemmas \ref{full.corr} and \ref{deg.corr}
we establish upper semicontinuity
of the minimal Willmore energy
under fixed Teichm\"uller class $\ \mini\ $.

\begin{proposition} \label{ele.upp}

$\ \mini: \teich \rightarrow [\beta_p^n,\infty]\ $
is upper semicontinuous.
Moreover $\ \mini\ $ is continuous at $\ \tau \in \teich, n = 3,4,\ $ with
\begin{displaymath}
	\mini(\tau) \leq \W_{n,p}.
\end{displaymath}
\end{proposition}
{\pr Proof:} \\
For the upper semicontinuity,
we have to prove
\begin{equation} \label{ele.upp.est}
	\limsup \limits_{\tau \rightarrow \tau_0}
	\mini(\tau) \leq \mini(\tau_0)
	\quad \mbox{for } \tau_0 \in \teich.
\end{equation}
It suffices to consider $\ \mini(\tau_0) < \infty\ $.
In this case, there exists for any $\ \varepsilon > 0\ $
a smooth immersion $\ f: \sur \rightarrow \rel^n
\mbox{ with } \tp(f^* \geu) = \tau_0
\mbox{ and } \W(f) < \mini(\tau_0) + \varepsilon\ $.

By Lemmas \ref{full.corr} and \ref{deg.corr}
applied to the constant sequence $\ f_m = f \mbox{ and } \vari = 0\ $,
there exist for any $\ \tau\ $ close enough to $\ \tau_0\ $
in Teichm\"uller space $\ \lambda_\tau \in \rel^{\dim \teich + 1}\ $
with
\begin{displaymath}
\begin{array}{c}
	\tp((f + \lambda_{\tau,r} \vari_r)^* \geu) = \tau, \\
	\lambda_\tau \rightarrow 0
	\mbox{ for } \tau \rightarrow \tau_0.
\end{array}
\end{displaymath}
Clearly $\ f + \lambda_{\tau,r} \vari_r \rightarrow f\ $ smoothly,
hence
\begin{displaymath}
	\limsup \limits_{\tau \rightarrow \tau_0} \mini(\tau)
	\leq \lim \limits_{\tau \rightarrow \tau_0}
	\W(f + \lambda_{\tau,r} \vari_r)
	= \W(f) \leq \mini(\tau_0) + \varepsilon
\end{displaymath}
and (\ref{ele.upp.est}) follows.

If $\ \mini\ $ were not continuous
at $\ \tau_0 \in \teich \mbox{ with } \mini(\tau_0) \leq \W_{n,p}\ $,
by upper semicontinuity of $\ \mini\ $ proved above,
there exists $\ \delta > 0\ $
and a sequence $\ \tau_m \rightarrow \tau_0 \mbox{ in } \teich\ $ with
\begin{displaymath}
	\mini(\tau_m) \leq \mini(\tau_0) - 2 \delta.
\end{displaymath}
We select smooth immersions $\ f_m: \sur \rightarrow \rel^n
\mbox{ with } \tp(f_m^* \geu) = \tau_m \rightarrow \tau_0\ $ and
\begin{equation} \label{mini.cont.sel}
	\W(f_m) \leq \mini(\tau_m) + 1/m,
\end{equation}
hence for $\ m\ $ large enough
$\ \W(f_m) \leq \mini(\tau_0) - \delta \leq \W_{n,p} - \delta\ $.
Replacing $\ f_m \mbox{ by } \moeb_m \circ f_m \circ \diffeo_m\ $
as in Proposition \ref{direct.prop} does neither change
the Willmore energy nor its projections in Teichm\"uller space,
and we may assume that $\ f_m,f\ $
satisfy (\ref{direct.prop.teich}) - (\ref{direct.prop.imm}).
By Lemmas \ref{full.corr} and \ref{deg.corr}
applied to $\ f_m, \vari = 0 \mbox{ and } \tau = \tau_0\ $,
there exist $\ \lambda_m \in \rel^{\dim \teich + 1}
\mbox{ for } m\ $ large enough with
\begin{displaymath}
\begin{array}{c}
	\tp((f_m + \lambda_{m,r} \vari_r)^* \geu) = \tau_0, \\
	\lambda_m \rightarrow 0
	\mbox{ for } m \rightarrow \infty.
\end{array}
\end{displaymath}
This yields
\begin{displaymath}
	\mini(\tau_0)
	\leq \liminf \limits_{m \rightarrow \infty}
	\W(f_m + \lambda_{m,r} \vari_r)
	\leq \liminf \limits_{m \rightarrow \infty}
	( \W(f_m) + C |\lambda_m| )
	\leq \mini(\tau_0) - \delta,
\end{displaymath}
which is a contradiction,
and the proposition is proved.
\proof
Secondly, we prove that the infimum taken
in Definition \ref{direct.def} over smooth
immersions is not improved by weak local
bilipschitz immersion in $\ W^{2,2}\ $.

\begin{proposition} \label{ele.inf}

Let $\ f: \sur \rightarrow \rel^n\ $
be a weak local bilipschitz $\ W^{2,2}-$immersion
with pull-back metric $\ f^* \geu = e^{2u} \gp\ $
conformal to a smooth unit volume constant curvature metric $\ \gp\ $.
Then
\begin{displaymath}
	\mini(\tp(\gp))
	\leq \W(f)
	= \frac{1}{4} \int \limits_\sur
	|A_f|^2 \d \mu_f + 2 \pi (1 - p).
\end{displaymath}
\end{proposition}
{\pr Proof:} \\
We approximate $\ f\ $ by smmoth immersions $\ f_m\ $
as in (\ref{direct.imm.strong})
and (\ref{direct.imm.comp}).
Putting $\ g := f^* \geu, g_m := f_m^* \geu\ $
and writing $\ g_m = e^{2u_m} \gpm\ $
for some unit volume constant curvature metric $\ \gpm\ $
by Poincar\'e's Theorem,
see \bcite{fisch.trom.conf}, \bcite{trom.teich},
we get by Proposition \ref{conf.strong}
\begin{equation} \label{ele.inf.conf}
	\parallel u_m \parallel_{L^\infty(\sur)},
	\parallel \nabla u_m
	\parallel_{L^2(\sur)} \leq C
\end{equation}
for some $\ C < \infty\ $ independent of $\ m\ $.
In local charts, we see
\begin{displaymath}
\begin{array}{c}
	g_m \rightarrow g
	\quad \mbox{strongly in } W^{1,2},
	\mbox{weakly}^* \mbox{ in } L^\infty, \\
	\Gamma_{g_m,ij}^k  \rightarrow \Gamma_{g,ij}^k
	\quad \mbox{strongly in } L^2,
\end{array}
\end{displaymath}
and
\begin{displaymath}
	A_{f_m,ij} = \nabla_i^{g_m} \nabla_j^{g_m} f_m
	\rightarrow \partial_{ij} f - \Gamma_{g,ij}^k \partial_k f =
\end{displaymath}
\begin{displaymath}
	= \nabla_i^{g} \nabla_j^{g} f = A_{f,ij}
	\quad \mbox{strongly in } L^2.
\end{displaymath}
Therefore by (\ref{intro.gauss})
\begin{equation} \label{ele.inf.will-conv}
	\W(f_m) \rightarrow \W(f)
	= \frac{1}{4} \int \limits_\sur
	|A_f|^2_g \d \mu_g + 2 \pi (1 - p)
\end{equation}
and
\begin{equation} \label{ele.inf.teich-conv}
	\tp(g_m) \rightarrow \tp(\gp).
\end{equation}
Next for a slice $\ {\cal S}(\gp)\ $ of unit volume
constant curvature metrics for $\ \tp(\gp) \in \teich\ $
around $\ \gp\ $,
see \bcite{fisch.trom.conf}, \bcite{trom.teich},
there exist unique $\ \tgpm \in {\cal S}(\gp)
\mbox{ with } \tp(\tgpm) = \tp(g_m) \rightarrow
\tp(\gp)\ $ for $\ m\ $ large enough,
hence
\begin{displaymath}
	\diffeo_m^* \gpm = \tgpm
	\rightarrow \gp
	\quad \mbox{smoothly}
\end{displaymath}
for suitable diffeomorphisms $\ \diffeo_m
\mbox{ of } \sur\ $ homotopic to the identity.

Next by (\ref{ele.inf.conf})
\begin{equation} \label{ele.inf.grad}
	\parallel D \diffeo_m \parallel_{L^\infty(\sur)},
	\parallel D \diffeo_m^{-1} \parallel_{L^\infty(\sur)}
	\leq C
\end{equation}
and
\begin{displaymath}
	\parallel \gpm \parallel_{W^{1,2}(\sur)}
	\leq \parallel e^{-2 u_m} \parallel_{W^{1,2}(\sur)
	\cap L^\infty(\sur)}
	\parallel g_m \parallel_{W^{1,2}(\sur) \cap L^\infty(\sur)}
	\leq C
\end{displaymath}
and for a subsequence $\ \gpm \rightarrow \tilde g,
u_m \rightarrow \tilde u \mbox{ weakly in } W^{1,2}(\sur)\ $,
in particular $\ \tilde g \leftarrow \gpm = e^{-2u_m} g_m
\rightarrow e^{-2 \tilde u} g\ $, and
\begin{equation} \label{ele.inf.limit}
	\tilde g = e^{-2 \tilde u} g = e^{2(u - \tilde u)} \gp
\end{equation}
is conformal to the smooth metric $\ \gp\ $.

The Theorem of Ebin and Palais,
see \bcite{fisch.trom.conf}, \bcite{trom.teich},
and the remarks following
imply after appropriately modifying $\ \diffeo_m\ $
\begin{equation} \label{ele.inf.id}
	\diffeo_m \rightarrow id_\sur
	\quad \mbox{weakly in } W^{2,2}(\sur),
	\mbox{weakly}^* \mbox{ in } W^{1,\infty}(\sur)
\end{equation}
and $\ \gp = id_\sur^* \tilde g
= e^{2(u - \tilde u)} \gp\ $
by (\ref{ele.inf.limit}), hence
\begin{equation} \label{ele.inf.factor}
	u_m \rightarrow \tilde u = u
	\quad \mbox{ weakly in } W^{1,2}(\sur),
	\mbox{weakly}^* \mbox{ in } L^\infty(\sur).
\end{equation}
Putting $\ \tilde f_m := f_m \circ \diffeo_m,
\tilde u_m := u_m \circ \diffeo_m\ $,
we see $\ \tilde f_m^* \geu = \diffeo_m^*(e^{2 u_m} \gpm)
= e^{2 \tilde u_m} \tgpm\ $
and by (\ref{ele.inf.grad}), (\ref{ele.inf.id})
and (\ref{ele.inf.factor})
\begin{displaymath}
\begin{array}{c}
	\tilde f_m \rightarrow f \mbox{ weakly in } W^{2,2}(\sur),
	\mbox{weakly}^* \mbox{ in } W^{1,\infty}(\sur), \\
	\tilde u_m \rightarrow u \mbox{ weakly in } W^{1,2}(\sur),
	\mbox{weakly}^* \mbox{ in } L^\infty(\sur).
\end{array}
\end{displaymath}
Therefore $\ \tilde f_m,f\ $ satisfy
(\ref{direct.prop.teich}) - (\ref{direct.prop.imm}).
By Lemmas \ref{full.corr} and \ref{deg.corr}
applied to $\ \vari = 0\ $,
there exist $\ \lambda_m \in \rel^{\dim \teich + 1}\ $
for $\ m\ $ large enough with
\begin{displaymath}
\begin{array}{c}
	\tp((\tilde f_m + \lambda_{m,r} \vari_r)^* \geu) = \tp(\gp), \\
	\lambda_m \rightarrow 0
	\mbox{ for } m \rightarrow \infty,
\end{array}
\end{displaymath}
as $\ \tp(\tilde f_m^* \geu) = \tp(g_m) \rightarrow \tp(\gp)\ $
by (\ref{ele.inf.teich-conv}).
Then
\begin{displaymath}
\begin{array}{c}
	\tp((f_m + \lambda_{m,r} (\vari_r \circ \diffeo_m^{-1}))^* \geu)
	= \tp(\gp), \\
	f_m + \lambda_{m,r} (\vari_r \circ \diffeo_m^{-1}) \rightarrow f
	\quad \mbox{strongly in } W^{2,2}(\sur),
	\mbox{weakly}^* \mbox{ in } W^{1,\infty}(\sur), \\
	c_0 \leq \Big( f_m
	+ \lambda_{m,r} (\vari_r \circ \diffeo_m^{-1}) \Big)^* \geu
	\leq C
\end{array}
\end{displaymath}
for some $\ 0 < c_0 \leq C < \infty
\mbox{ and } m\ $ large,
when observing that $\ \diffeo_m^{-1}\ $ is bounded
in $\ W^{2,2}(\sur) \cap W^{1,\infty}(\sur)\ $
by (\ref{ele.inf.grad}), (\ref{ele.inf.id})
and $\ |D^2(\diffeo_m^{-1})|
\leq C |D^2(\diffeo_m) \circ \diffeo_m^{-1}|\ $.
As in (\ref{ele.inf.will-conv}) we get
\begin{displaymath}
	\mini(\tp(\gp))
	\leq \lim \limits_{m \rightarrow \infty}
	\W(f_m + \lambda_{m,r} (\vari_r \circ \diffeo_m^{-1}))
	= \W(f),
\end{displaymath}
and recalling (\ref{ele.inf.will-conv}) again,
the proposition follows.
\proof
This proposition implies strong convergence
of minimizing sequences.

\begin{proposition} \label{ele.strong}

Let $\ f: \sur \rightarrow \rel^n\ $ be a weak local bilipschitz immersion
approximated by smooth immersions $\ f_m\ $
satisfying (\ref{direct.prop.teich}) - (\ref{direct.prop.imm}) and
\begin{equation} \label{ele.strong.mini}
	\W(f_m) \leq \mini(\tp(f_m^* \geu)) + \varepsilon_m
\end{equation}
with $\ \varepsilon_m \rightarrow 0\ $.
Then
\begin{equation} \label{ele.strong.strong}
	f_m \rightarrow f
	\quad \mbox{strongly in } W^{2,2}(\sur)
\end{equation}
and
\begin{equation} \label{ele.strong.minimizer}
	\W(f) = \mini(\tau_0).
\end{equation}
\end{proposition}
{\pr Proof:} \\
By (\ref{direct.prop.w22}), (\ref{direct.prop.met}),
(\ref{direct.prop.conv}) and $\ \Lambda_0\ $
large enough,
we may assume after relabeling the sequence $\ f_m\ $
\begin{equation} \label{ele.strong.d2}
\begin{array}{c}
	\parallel 
	D f_m \parallel_{W^{1,2}(\sur)
	\cap L^\infty(\sur)} \leq \Lambda_0, \\
	\frac{1}{2} \gp \leq \gpm \leq 2 \gp, \\
	\Lambda_0^{-1} \gp \leq g_m = f_m^* \geu \leq \Lambda_0 \gp, \\
	\int \limits_\sur |A_{f_m}|^2 \d \mu_{f_m}
	\leq \Lambda_0.
\end{array}
\end{equation}
In local charts, we see
\begin{displaymath}
\begin{array}{c}
	g_m \rightarrow g
	\quad \mbox{weakly in } W^{1,2},
	\mbox{weakly}^* \mbox{ in } L^\infty, \\
	\Gamma_{g_m,ij}^k  \rightarrow \Gamma_{g,ij}^k
	\quad \mbox{weakly in } L^2,
\end{array}
\end{displaymath}
and
\begin{displaymath}
	A_{f_m,ij} = \nabla_i^{g_m} \nabla_j^{g_m} f_m
	\rightarrow \partial_{ij} f - \Gamma_{g,ij}^k \partial_k f =
\end{displaymath}
\begin{displaymath}
	= \nabla_i^{g} \nabla_j^{g} f = A_{f,ij}
	\quad \mbox{weakly in } L^2.
\end{displaymath}
We conclude by Propositions \ref{ele.upp},
\ref{ele.inf} and (\ref{ele.strong.mini}),
as $\ \tp(f_m^* \geu) \rightarrow \tau_0\ $ by (\ref{direct.prop.teich}),
\begin{displaymath}
	\W(f)
	\leq \liminf \limits_{m \rightarrow \infty} \W(f_m)
	\leq \liminf \limits_{m \rightarrow \infty} \mini(f_m^* \geu)
	\leq \mini(\tau_0) \leq \W(f),
\end{displaymath}
hence $\ \mean_{f_m} \rightarrow \mean_f \mbox{ strongly in } L^2\ $.
This yields using (\ref{direct.prop.met})
\begin{displaymath}
	\Delta_\gpm f_m = e^{2 u_m} \mean_{f_m}
	\rightarrow e^{2u} \mean_f = \Delta_\gp f
	\quad \mbox{strongly in } L^2
\end{displaymath}
and
\begin{displaymath}
	\partial_i \Big( \gpo{ij} \sqrt{\gp} \partial_j(f_m - f) \Big)
	= \partial_i \Big( (\gpo{ij} \sqrt{\gp}
	- \gpmo{ij} \sqrt{\gpm}) \partial_j f_m \Big) +
\end{displaymath}
\begin{displaymath}
	+ \sqrt{\gpm} \Delta_\gpm f_m
	- \sqrt{\gp} \Delta_\gp f
	\rightarrow 0
	\quad \mbox{strongly in } L^2
\end{displaymath}
recalling that $\ \partial_j f_m\ $ is bounded in $\ L^\infty\ $
and $\ \gpm \rightarrow \gp\ $ smoothly by (\ref{direct.prop.met}),
hence
\begin{displaymath}
	f_m \rightarrow f
	\quad \mbox{strongly in } W^{2,2}(\sur),
\end{displaymath}
and the proposition is proved.
\proof


\setcounter{equation}{0}

\section{Decay of the second derivative} \label{dec}

In this section, we add to our assumptions
on $\ f, f_m\ $, as considered in \S \ref{direct} - \S \ref{deg},
that $\ f_m\ $ is approximately minimizing in its Teichm\"uller class,
see (\ref{dec.prop.mini}).
The aim is to prove in the following proposition
a decay for the second derivatives
which implies that the limits in $\ C^{1,\alpha}\ $.

\begin{proposition} \label{dec.prop}

Let $\ f: \sur \rightarrow \rel^n\ $ be a weak local bilipschitz immersion
approximated by smooth immersions $\ f_m\ $
satisfying (\ref{direct.prop.teich}) - (\ref{direct.prop.imm}) and
\begin{equation} \label{dec.prop.mini}
	\W(f_m) \leq \mini(\tp(f_m^* \geu)) + \varepsilon_m
\end{equation}
with $\ \varepsilon_m \rightarrow 0\ $.
Then there exists $\ \alpha > 0, C < \infty\ $ such that
\begin{equation} \label{dec.prop.dec}
	\int \limits_{B_\varrho^\gp(x)} |\nabla^2_\gp f|^2_\gp \d \mu_\gp
	\leq C \varrho^{2 \alpha}
	\quad \mbox{for any }  x \in \sur, \varrho > 0,
\end{equation}
in particular $\ f \in C^{1,\alpha}(\sur)\ $.
\end{proposition}
{\pr Proof:} \\
By (\ref{direct.prop.w22}), (\ref{direct.prop.met}) and $\ \Lambda_0\ $
large enough, we may assume after relabeling the sequence $\ f_m\ $
\begin{equation} \label{dec.prop.d2}
\begin{array}{c}
	\parallel u_m, D f_m \parallel_{W^{1,2}(\sur)
	\cap L^\infty(\sur)} \leq \Lambda_0, \\
	\frac{1}{2} \gp \leq \gpm \leq 2 \gp, \\
	\Lambda_0^{-1} \gp \leq g_m = f_m^* \geu \leq \Lambda_0 \gp, \\
	\int \limits_\sur |A_{f_m}|^2 \d \mu_{f_m}
	\leq \Lambda_0.
\end{array}
\end{equation}
Putting $\ \funda_m := |\nabla^2_\gp f_m|_\gp^2 \mu_{\gp}\ $,
we see $\ \funda_m(\sur) \leq C(\Lambda_0,\gp) $
and for a subsequence $\ \funda_m \rightarrow \funda\ $ 
weakly$^* \mbox{ in } C_0^0(\sur)^*
\mbox{ with } \funda(\sur) < \infty\ $.

We consider $\ x_0 \in \sur
\mbox{ with a neighbourhood } U_0(x_0)\ $
satisfying
\begin{equation} \label{dec.prop.good-punc}
	\funda(\overline{U_0(x_0)} - \{  x_0 \}) < \varepsilon_1,
\end{equation}
where we choose $\ \varepsilon_1 = \varepsilon_1(\Lambda_0,n) > 0\ $ below,
together with a chart $\ \varphi: U_0(x_0)
\stackrel{\approx}{\longrightarrow} B_{2 \varrho_0}(0),
\varphi(x_0) = 0, U_\varrho(x_0) := \varphi^{-1}(B_\varrho(0))
\mbox{ for } 0 < \varrho \leq 2 \varrho_0 \leq 2,\ $
such that
\begin{equation} \label{dec.prop.metric}
	\frac{1}{2} \geu \leq (\varphi^{-1})^* \gp \leq 2 \geu
\end{equation}
and in the coordinates of the chart $\ \varphi\ $
\begin{equation} \label{dec.prop.christ}
	\int \limits_{B_{\varrho_0}(0)} \gpo{ik} \gpo{jl} \gpu{rs}
	\Gamma_{\gp,ij}^r \Gamma_{\gp,kl}^s \d \mu_\gp
	< \varepsilon_1.
\end{equation}
Moreover we select
\begin{equation} \label{dec.prop.corr}
	U_{\varrho_1}(x_0)
	\subseteq U(x_0)
	\subseteq U_{\varrho_0}(x_0),
\end{equation}
variations $\ \vari_1, \ldots, \vari_{\dim \teich - 1},
\vari_\pm \in C^\infty_0(\sur - \overline{U(x_0)},\rel^n)\ $
and $\ \delta > 0, C = C_{x_0,\varphi} < \infty, m_0 \in \nat\ $,
as in Lemmas \ref{full.corr} and \ref{deg.corr}
for $\ x_0, U_{\varrho_0}(x_0)
\mbox{ and } \Lambda := C(\Lambda_0)\ $
defined below.
As $\ \tp(f_m^*) \rightarrow \tp(\gp)\ $ by (\ref{direct.prop.teich}),
we get for $\ m_0\ $ large enough
\begin{equation} \label{dec.prop.teich-close}
	d_\teich( \tp(f_m^* \geu),\tp(\gp) ) < \varepsilon
	\quad \mbox{for } m \geq m_0.
\end{equation}
Clearly for each $\ x_0 \in \sur\ $,
there exist $\ U_0(x_0), \varrho_0\ $ as above,
since $\ \nu(B_\varrho^\gp(x_0) - \{x_0\})
\rightarrow \nu(\emptyset) = 0
\mbox{ for } \varrho \rightarrow 0\ $.

For $\ x_0, U_0(x_0)\ $ as above and $\ 0 < \varrho \leq \varrho_0\ $,
there exists $\ m_1 \geq m_0 \mbox{ such that }
\funda_m(U_{\varrho_0}(x_0) - U_{\varrho/2}(x_0)) < \varepsilon_1
\mbox{ for } m \geq m_1\ $,
hence in the coordinates of the chart $\ \varphi\ $
\begin{displaymath}
	\int \limits_{B_{\varrho_0}(0) - B_{\varrho/2}(0)} |D^2 f_m|^2 \d \Lt
	\leq C \int \limits_{B_{\varrho_0}(0) - B_{\varrho/2}(0)}
	\gpo{ik} \gpo{jl} \langle \partial_{ij} f_m ,
	\partial_{kl} f_m \rangle \d \mu_\gp \leq
\end{displaymath}
\begin{displaymath}
	\leq 2 \hspace{-.3cm} \int \limits_{U_{\varrho_0}(0) - U_{\varrho/2}(x_0)}
	\hspace{-1.1cm} |\nabla_\gp f_m|^2_\gp \d \mu_\gp
	+ 2 \hspace{-.3cm} \int \limits_{B_{\varrho_0}(x)} \hspace{-.2cm} \gpo{ik} \gpo{jl}
	\Gamma_{\gp,ij}^r \Gamma_{\gp,kl}^s
	\langle \partial_r f_m , \partial_s f_m \rangle \d \mu_\gp \leq
\end{displaymath}
\begin{equation} \label{dec.prop.sec}
	\leq 2 \varepsilon_1
	+ C(\Lambda_0) \int \limits_{B_{\varrho_0}(0)} \gpo{ik} \gpo{jl} \gpu{rs}
	\Gamma_{\gp,ij}^r \Gamma_{\gp,kl}^s \d \mu_\gp
	\leq C(\Lambda_0) \varepsilon_1.
\end{equation}
There exists $\ \sigma \in ]3/\varrho/4,7 \varrho/8[\ $
satisfying by Co-Area formula, see \bcite{sim} \S 12,
\begin{displaymath}
	\int \limits_{\partial B_\sigma(0)} |D^2 f_m|^2 \d \Ho
	\leq 8 \varrho^{-1} \int \limits_{3 \varrho/4}^{7 \varrho/8}
	\int \limits_{\partial B_r(0)} |D^2 f_m|^2 \d \Ho \d r =
\end{displaymath}
\begin{equation} \label{dec.prop.sigma}
	= 8 \varrho^{-1} \int \limits_{B_{7 \varrho/8}(0)
	- B_{3 \varrho/4}(0)}
	|D^2 f_m|^2 \d \Lt
	\leq C(\Lambda_0) \varepsilon_1 \varrho^{-1}.
\end{equation}
First we conclude that
\begin{displaymath}
	osc_{\partial B_\sigma(0)} |D f_m|
	\leq C \int \limits_{\partial B_\sigma(0)} |D^2 f_m| \d \Ho
	\leq C(\Lambda_0) \varepsilon_1^{1/2},
\end{displaymath}
hence for any $\ x \in \partial B_\sigma(0)\ $
and the affine function $\ l(y)
:= f_m(x) + D f_m(x) (y-x)\ $
\begin{displaymath}
	\sigma^{-1} \parallel f_m - l
	\parallel_{L^\infty(\partial B_\sigma(0))}
	+ \parallel D (f_m - l)\parallel_{L^\infty(\partial B_\sigma(0))})
	\leq C(\Lambda_0) \varepsilon_1^{1/2}.
\end{displaymath}
Moreover by (\ref{dec.prop.d2}) and (\ref{dec.prop.metric})
\begin{displaymath}
	c_0(\Lambda_0) (\delta_{ij})_{ij}
	\leq (\langle \partial_i l , \partial_j l \rangle )_{ij}
	= g_m(x) \leq C(\Lambda_0) (\delta_{ij})_{ij}
\end{displaymath}
hence $\ Dl \in \rel^{2 \times 2}\ $ is invertible and
\begin{displaymath}
	\parallel D l \parallel,
	\parallel (Dl)^{-1} \parallel \leq C(\Lambda_0).
\end{displaymath}
Next by standard trace extension lemma,
there exists $\ \tilde f_m \in C^2(\overline{B_\sigma(0)})\ $
such that
\begin{equation} \label{dec.prop.trace}
\begin{array}{c}
        \tilde f_m = f_m, D \tilde f_m = D f_m
        \quad \mbox{on } \partial B_\sigma(0), \\
        \sigma^{-1} |\tilde f_m - l|
	+ |D (\tilde f_m - l)| \leq \\
        \leq C (\sigma^{-1} \parallel f_m - l
	\parallel_{L^\infty(\partial B_\sigma(0))}
	+ \parallel D (f_m - l) \parallel_{L^\infty(\partial B_\sigma(0))})
	\leq C(\Lambda_0) \varepsilon_1^{1/2}, \\
        \int \limits_{B_\sigma(0)}
        |D^2 \tilde f_m|^2 \d \Lt
        \leq C \sigma
	\int \limits_{\partial B_\sigma(0)}
        |D^2 f_m|^2 \d \Ho
	\leq C(\Lambda_0) \varepsilon_1.
\end{array}
\end{equation}
We see
\begin{displaymath}
	\parallel D \tilde f_m - D l \parallel
	\leq C(\Lambda_0) \varepsilon_1^{1/2}
	\leq \parallel (Dl)^{-1} \parallel^{-1} /2
\end{displaymath}
for $\ \varepsilon_1 = \varepsilon_1(\Lambda_0) \leq 1\ $ small enough,
and $\ D \tilde f_m\ $ is of full rank everywhere with
\begin{displaymath}
	\parallel D \tilde f_m \parallel,
	\parallel (D \tilde f_m)^{-1} \parallel
	\leq C(\Lambda_0).
\end{displaymath}
Extending $\ \tilde f_m = f_m \mbox{ on } \sur - U_\sigma(x_0)\ $,
we see that $\ \tilde f_m:\sur \rightarrow \rel^n\ $
is a $\ C^{1,1}-$immersion
with pullback metric satisfying by (\ref{dec.prop.metric})
\begin{equation} \label{dec.prop.vari-metric}
	c_0(\Lambda_0) \gp \leq c_0(\Lambda_0) \geu
	\leq \tilde g_m := \tilde f_m^* \geu
	\leq C(\Lambda_0) \geu \leq C(\Lambda_0) \gp.
\end{equation}
and by (\ref{dec.prop.christ}),
(\ref{dec.prop.sec}) and (\ref{dec.prop.trace})
\begin{equation} \label{dec.prop.vari-sec}
	\int \limits_{B_{\varrho_0}(0)}
	|D^2 \tilde f_m|^2 \d \Lt,
	\int \limits_{B_{\varrho_0}(0)}
	|\nabla_\gp \tilde f_m|^2_\gp \d \mu_\gp
	\leq C(\Lambda_0) \varepsilon_1,
\end{equation}
in particular, as $\ |A_{\tilde f_m,ij}|
\leq |\partial_{ij} \tilde f_m|\ $
in local coordinates,
\begin{equation} \label{dec.prop.vari-funda}
	\int \limits_{U_{\varrho_0}(x_0)}
	|A_{\tilde f_m}|^2 \d \mu_{\tilde g_m}
	\leq C(\Lambda_0) \varepsilon_1 \leq \varepsilon_0(n)
\end{equation}
for $\ \varepsilon_1 = \varepsilon_1(\Lambda_0,n)\ $ small enough.
Putting $\ V := \tilde f_m - f_m\ $,
we see $\ supp\ V \subset \subset U_\varrho(x_0)
\subseteq U(x_0) \mbox{ for } 0 < \varrho \leq \varrho_1\ $
from (\ref{dec.prop.corr}) and
\begin{displaymath}
	\parallel V \parallel_{W^{2,2}(\sur)}
	\leq C \parallel V \parallel_{W^{2,2}(B_\sigma(0))}
	\leq C(\Lambda_0).
\end{displaymath}
Together with (\ref{dec.prop.vari-metric})
and (\ref{dec.prop.vari-funda}),
this verifies (\ref{full.corr.ass}) and (\ref{deg.corr.ass})
for $\ \Lambda = C(\Lambda_0)\ $
in Lemmas \ref{full.corr} and \ref{deg.corr},
respectively.
After slightly smoothing $\ \vari\ $,
there exists $\ \lambda_m \in \rel^{\dim \teich + 1}\ $
by Lemmas \ref{full.corr} and \ref{deg.corr}
and (\ref{dec.prop.teich-close})
with
\begin{equation} \label{dec.prop.teich}
\begin{array}{c}
	\tp((\tilde f_m + \lambda_{m,r} \vari_r)^* \geu)
	= \tp(f_m^* \geu), \\
	|\lambda_m| \leq C_{x_0,\varphi} d_\teich \Big( \tp( \tilde f_m^* \geu),
	\tp(f_m^* \geu) \Big)^{1/2}.
\end{array}
\end{equation}
By the minimizing property (\ref{dec.prop.mini}),
and the Gau\ss-Bonnet Theorem in (\ref{intro.gauss}),
we get
\begin{displaymath}
	\frac{1}{4} \int \limits_\sur |A_{f_m}|^2 \d \mu_{f_m}
	- \varepsilon_m
	\leq \mini(\tp(f_m^* \geu)) + 2 \pi (p-1)
	\leq \frac{1}{4} \int \limits_\sur
	|A_{\tilde f_m + \lambda_{m,r} \vari_r}|^2
	\d \mu_{\tilde f_m + \lambda_{m,r} \vari_r},
\end{displaymath}
hence, as $\ \tilde f_m = f_m \mbox{ in } \sur - U_\sigma(x_0)
\mbox{ and } supp\ \vari_r \cap \overline{U_\sigma(x_0)} = \emptyset\ $,
\begin{equation} \label{dec.prop.iter-start}
	\int \limits_{U_\sigma(x_0)} |A_{f_m}|^2 \d \mu_{f_m}
	\leq \int \limits_{U_\sigma(x_0)} |A_{\tilde f_m}|^2
	\d \mu_{\tilde f_m} + C_{x_0,\varphi}(\Lambda_0,\gp) |\lambda_m|
	+ 4 \varepsilon_m.
\end{equation}
We continue, using $\ |A_{\tilde f_m,ij}|
\leq |\partial_{ij} \tilde f_m|\ $
in local coordinates,
(\ref{dec.prop.sigma}) and (\ref{dec.prop.trace}),
\begin{displaymath}
	\int \limits_{U_\sigma(x_0)} |A_{\tilde f_m}|^2
	\d \mu_{\tilde f_m} \leq
\end{displaymath}
\begin{displaymath}
	\leq C(\Lambda_0) \int \limits_{B_\sigma(0)} |D^2 \tilde f_m|^2 \d \Lt
	\leq C(\Lambda_0) \int \limits_{B_{7 \varrho/8}(0)
	- B_{3 \varrho/4}(0)} |D^2 f_m|^2 \d \Lt \leq
\end{displaymath}
\begin{displaymath}
	\leq C(\Lambda_0) \int \limits_{U_{7 \varrho/8}(x_0)
	- U_{3 \varrho/4}(x_0)} |\nabla_\gp^2 f_m|^2 \d \mu_\gp
	+  C(\Lambda_0)
	\parallel \Gamma_\gp \parallel_{L^\infty(B_{\varrho_0}(0))}^2
	\varrho^2 \leq
\end{displaymath}
\begin{equation} \label{dec.prop.iter-upp}
	\leq C(\Lambda_0) \nu_m \Big( U_{7 \varrho/8}(x_0)
	- U_{3 \varrho/4}(x_0) \Big) + C_{x_0,\varphi}(\Lambda_0,\gp) \varrho^2.
\end{equation}
We calculate in local coordinates in $\ B_\sigma(0)\ $ that
\begin{displaymath}
	\Delta_\gpm f_m = e^{2 u_m} \Delta_{g_m} f_m
	= e^{2 u_m} \mean_{f_m}
\end{displaymath}
and
\begin{displaymath}
	|\Delta_\gpm (f_m - \tilde f_m)|
	\leq C(\Lambda_0) ( |\mean_{f_m}|
	+ |D^2 \tilde f_m| + |\Gamma_\gpm|\ |D \tilde f_m| ).
\end{displaymath}
By standard elliptic theory, see \bcite{gil.tru} Theorem 8.8,
from (\ref{direct.prop.met}) for $\ m \geq m_1\ $ large enough,
as $\ f_m = \tilde f_m \mbox{ on } \partial B_\sigma(0)\ $,
we get
\begin{displaymath}
	\int \limits_{B_\sigma(0)}
	|D^2 f_m|^2 \d \Lt \leq
\end{displaymath}
\begin{displaymath}
	\leq C_{x_0,\varphi}(\Lambda_0,\gp)
	\Big( \int \limits_{U_\sigma(x_0)} |\mean_{f_m}|^2
	\d \mu_{f_m}
	+ \int \limits_{B_\sigma(0)}
	|D^2 \tilde f_m|^2 \d \Lt
	+ \int \limits_{B_\sigma(0)}
	|D \tilde f_m|^2 \d \Lt \Big) \leq
\end{displaymath}
\begin{equation} \label{dec.prop.iter-low}
	\leq C_{x_0,\varphi}(\Lambda_0,\gp)
	\Big( \int \limits_{U_\sigma(x_0)} |\mean_{f_m}|^2
	\d \mu_{f_m}
	+ \nu_m \Big( U_{7 \varrho/8}(x_0)
	- U_{3 \varrho/4}(x_0) \Big) + \varrho^2 \Big),
\end{equation}
where we have used (\ref{dec.prop.iter-upp}).
Putting (\ref{dec.prop.iter-start}),
(\ref{dec.prop.iter-upp}) and (\ref{dec.prop.iter-low})
together yields
\begin{displaymath}
	\nu_m(U_{\varrho/2}(x_0))
	= \int \limits_{U_{\varrho/2}(x_0)} |\nabla_\gp^2 f_m|^2 \d \mu_\gp \leq
\end{displaymath}
\begin{equation} \label{dec.prop.iter-conc}
	\leq C_{x_0,\varphi}(\Lambda_0,\gp)
	\bigg( \nu_m \Big( U_{7 \varrho/8}(x_0)
	- U_{3 \varrho/4}(x_0) \Big)
	+ \varrho^2 + |\lambda_m| + \varepsilon_m \bigg).
\end{equation}
To estimate $\ \lambda_m\ $,
we continue
observing that $\ \tilde f_m^* \geu = \tilde g_m
\mbox{ and } f_m^* \geu = g_m\ $ coincide on $\ \sur - U_\varrho(x_0)\ $,
\begin{displaymath}
	d_\teich(\tp(\tilde f_m^* \geu),\tp(f_m^* \geu))^2 \leq
\end{displaymath}
\begin{displaymath}
	\leq 2 d_\teich(\tp(\tilde f_m^* \geu),\tp(\gp))^2
	+ 2 d_\teich(\tp(f_m^* \geu),\tp(\gp))^2 \leq
\end{displaymath}
\begin{displaymath}
	\leq C_{\tau_0} \int \limits_{\sur}
	\Big( \frac{1}{2} \gpu{ij} \tilde g_m^{ij}
	\sqrt{\tilde g_m} - \sqrt{\gp} \Big) \d x
	+ C_{\tau_0} \int \limits_{\sur}
	\Big( \frac{1}{2} \gpu{ij} g_m^{ij}
	\sqrt{g_m} - \sqrt{\gp} \Big) \d x =
\end{displaymath}
\begin{displaymath}
	= 2 C_{\tau_0} \int \limits_{\sur - U_\varrho(x_0)}
	\Big( \frac{1}{2} \gpu{ij} g_m^{ij}
	\sqrt{g_m} - \sqrt{\gp} \Big) \d x +
\end{displaymath}
\begin{displaymath}
	+ C_{\tau_0} \int \limits_{U_\varrho(x_0)}
	\Big( \frac{1}{2} \gpu{ij} \tilde g_m^{ij}
	\sqrt{\tilde g_m} - \sqrt{\gp} \Big) \d x
	+ C_{\tau_0} \int \limits_{U_\varrho(x_0)}
	\Big( \frac{1}{2} \gpu{ij} g_m^{ij}
	\sqrt{g_m} - \sqrt{\gp} \Big) \d x \leq
\end{displaymath}
\begin{displaymath}
	\leq 2 C_{\tau_0} \int \limits_{\sur - U_\varrho(x_0)}
	\Big| \frac{1}{2} \gpu{ij} g_m^{ij}
	\sqrt{g_m} - \sqrt{\gp} \Big| \d x
	+ C_{\tau_0}(\Lambda_0) \varrho^2.
\end{displaymath}
where we have used (\ref{dec.prop.d2})
and (\ref{dec.prop.vari-metric}).
As $\ g_m \rightarrow g = e^{2u} \gp
\mbox{ pointwise and bounded on } \sur\ $,
we get from Lebesgue's convergence theorem
\begin{displaymath}
	\limsup \limits_{m \rightarrow \infty}
	d_\teich(\tp(\tilde f_m^* \geu),\tp(f_m^* \geu))
	\leq C_{\tau_0}(\Lambda_0) \varrho
\end{displaymath}
and from (\ref{dec.prop.teich})
\begin{displaymath}
	\limsup \limits_{m \rightarrow \infty} |\lambda_m|
	\leq C_{x_0,\varrho_0\gp}(\Lambda_0) \varrho^{1/4}.
\end{displaymath}
Plugging into (\ref{dec.prop.iter-conc})
and passing to the limit $\ m \rightarrow \infty\ $,
we obtain
\begin{displaymath}
	\nu(U_{\varrho/2}(x_0))
	\leq C_{x_0,\varphi}(\Lambda_0,\gp)
	\nu \Big( U_{\varrho}(x_0)
	- U_{\varrho/2}(x_0) \Big)
	+ C_{x_0,\varphi}(\Lambda_0,\gp) \varrho^{1/4}
\end{displaymath}
and by hole-fllling
\begin{displaymath}
	\nu(U_{\varrho/2}(x_0))
	\leq \gamma
	\nu(U_{\varrho}(x_0))
	+ C_{x_0,\varphi}(\Lambda_0,\gp) \varrho^{1/4}
\end{displaymath}
with $\ \gamma = C / (C + 1) < 1\ $.
Iterating with \bcite{gil.tru} Lemma 8.23,
we arrive at
\begin{equation} \label{dec.prop.dec-x0}
	\nu(B_\varrho^\gp(x_0))
	\leq C_{x_0,\varphi}(\Lambda_0,\gp)
	\varrho^{2 \alpha} \varrho_1^{-2 \alpha}
	\quad \mbox{for all } \varrho > 0
\end{equation}
and some $\ 0 < \alpha
= \alpha_{x_0,\varphi}(\Lambda_0,\gp) < 1\ $.
Since $\ \nu(B_\varrho^\gp(x_0)) \rightarrow
\nu(\{x_0\}) \mbox{ for } \varrho \rightarrow 0\ $,
we first conclude
\begin{equation} \label{dec.prop.good}
	\nu(\{x_0\}) = 0.
\end{equation}
Then we can improve the choice
of $\ U_0(x_0)\ $ in (\ref{dec.prop.good-punc}) to
\begin{displaymath}
	\funda(\overline{U_0(x_0)}) < \varepsilon_1,
\end{displaymath}
and we can repeat the above iteration
for any $\ x \in U_{\varrho_1/2}(x_0), 0 < \varrho \leq \varrho_1/2\ $
to obtain
\begin{displaymath}
	\nu(B_\varrho^\gp(x))
	\leq C_{x_0,\varphi}(\Lambda_0,\gp)
	\varrho^{2 \alpha} \varrho_1^{-2 \alpha}
	\quad \mbox{for all } x \in U_{\varrho_1/2}(x_0) \varrho > 0.
\end{displaymath}
By a finite covering, this yields (\ref{dec.prop.dec}).
Since in the coordinates of the chart $\ \varphi\ $
\begin{displaymath}
	\int \limits_{B_{\varrho}(x)} |D^2 f_m|^2 \d \Lt
	\leq C \int \limits_{B_{\varrho}(x)}
	\gpo{ik} \gpo{jl} \langle \partial_{ij} f_m ,
	\partial_{kl} f_m \rangle \d \mu_\gp \leq
\end{displaymath}
\begin{displaymath}
	\leq 2 \int \limits_{B_{\varrho}(x)}
	|\nabla_\gp f_m|^2_\gp \d \mu_\gp
	+ 2 \int \limits_{B_{\varrho}(x)} \gpo{ik} \gpo{jl}
	\Gamma_{\gp,ij}^r \Gamma_{\gp,kl}^s
	\langle \partial_r f_m , \partial_s f_m \rangle \d \mu_\gp \leq
\end{displaymath}
\begin{displaymath}
	\leq C_{x_0,\varphi}(\Lambda_0,\gp)
	\varrho^{2 \alpha} \varrho_1^{-2 \alpha}
	+ C(\Lambda_0,\gp) \varrho^2,
\end{displaymath}
we conclude by Morrey's lemma,
see \bcite{gil.tru} Theorem 7.19,
that $\ f \in C^{1,\alpha}(\sur)\ $,
and the proposition is proved.
\proof


\setcounter{equation}{0}

\section{The Euler-Lagrange equation} \label{euler}

The aim of this section is to prove the Euler-Lagrange equation
for the limit of immersions
approximately minimizing under fixed Teichm\"uller class.
From this we will conclude full regularity of the limit.

\begin{theorem} \label{euler.theo}

Let $\ f: \sur \rightarrow \rel^n\ $ be a weak local bilipschitz immersion
approximated by smooth immersions $\ f_m\ $
satisfying (\ref{direct.prop.teich}) - (\ref{direct.prop.imm}) and
\begin{equation} \label{euler.theo.mini}
	\W(f_m) \leq \mini(\tp(f_m^* \geu)) + \varepsilon_m
\end{equation}
with $\ \varepsilon_m \rightarrow 0\ $.
Then $\ f\ $ is a smooth minimizer
of the Willmore energy under fixed Teichm\"uller class
\begin{equation} \label{euler.theo.minimizer}
	\W(f) = \mini(\tau_0)
\end{equation}
and satisfies the Euler-Lagrange equation
\begin{equation} \label{euler.theo.euler}
        \Delta_g \mean + Q(A^0) \mean = g^{ik} g^{jl} A_{ij}^0 q_{kl}
	\quad \mbox{on } \sur,
\end{equation}
where $\ q\ $ is a smooth transverse traceless
symmetric 2-covariant tensor with respect
to $\ g = f^* \geu\ $.
\end{theorem}
{\pr Proof:} \\
By Propositions \ref{full.basis} and \ref{deg.basis},
we select variations
variations $\ \vari_1, \ldots, \vari_{\dim \teich}
\in C^\infty_0(\sur,\rel^n)\ $,
satisfying (\ref{full.basis.span})
or variations $\ \vari_1, \ldots, \vari_{\dim \teich - 1},
\vari_\pm \in C^\infty_0(\sur,\rel^n)\ $,
satisfying (\ref{deg.basis.span}) and (\ref{deg.basis.pm}),
depending on whether $\ f\ $ has full rank in Teichm\"uller space
or not.

For $\ \vari \in C^\infty(\sur,\rel^n)\ $
and putting $\ f_{m,t,\lambda,\mu} :=
f_m + t \vari + \lambda_r \vari_r + \mu_\pm \vari_\pm\ $,
we see for $\ |t| \leq t_0
\mbox{ for some } t_0 = t_0(\vari,\Lambda_0,n) > 0\ $ small enough
and $\ |\lambda|, |\mu| \leq \lambda_0
\mbox{ for some } \lambda_0
= \lambda_0(\vari_r,\vari_\pm,\Lambda_0) > 0\ $
small enough that
\begin{displaymath}
\begin{array}{c}
	\parallel D f_{m,t,\lambda,\mu} \parallel_{W^{1,2}(\sur)
	\cap L^\infty(\sur)} \leq 2 \Lambda_0, \\
	(2 \Lambda_0)^{-1} \gp
	\leq f_{m,t,\lambda,\mu}^* \geu
	\leq 2 \Lambda_0 \gp
\end{array}
\end{displaymath}
and all $\ m \in \nat\ $.
As $\ f_{m,t,\lambda,\mu} \rightarrow f
\mbox{ strongly in } W^{2,2}(\sur)
\mbox{ and weakly}^* \mbox{ in } W^{1,\infty}(\sur)\ $
for $\ m \rightarrow \infty, t, \lambda, \mu \rightarrow 0\ $
by Proposition \ref{ele.strong},
we get from Proposition \ref{deg.conv}
and the remark following
for any chart $\ \psi: U(\tp(\gp)) \subseteq \teich
\rightarrow \rel^{\dim \teich}\ $,
and put $\ \tpc := \psi \circ \tp,
\delta \tpc = d \psi \circ \delta \tp,
\hat {\cal V}_f := d \psi_{\tp \gp} . {\cal V}_f,
\delta^2 \tpc\ $ defined in (\ref{deg.sec-def}),
and any $\ W \in C^\infty(\sur,\rel^n)\ $
\begin{equation} \label{euler.theo.conv}
\begin{array}{c}
	\tp(f_{m,t,\lambda,\mu}^* \geu) \rightarrow \tau_0, \\
	\delta \tpc_{f_{m,\lambda,\mu}} . W \rightarrow \delta \tpc_f . W,
	\quad
	\delta^2 \tpc_{f_{m,t,\lambda,\mu}}(W) \rightarrow \delta^2 \tpc_f(W)
\end{array}
\end{equation}
as $\ m \rightarrow \infty, t, \lambda, \mu \rightarrow 0\ $.
If $\ f\ $ is of full rank in Teichm\"uller space,
then $\ \hat {\cal V}_f = \rel^{\dim \teich}\ $,
and we put $\ d = \dim \teich\ $.
If $\ f\ $ is not of full rank in Teichm\"uller space,
we may assume fter a change of coordinates,
\begin{displaymath}
	\hat {\cal V}_f = \rel^{\dim \teich - 1} \times \{0\}
\end{displaymath}
and put $\ d := \dim \teich - 1
\mbox{ and } e := e_{\dim \teich} \perp \hat {\cal V}_f\ $.
By (\ref{full.basis.span}) or (\ref{deg.basis.span}),
we see for the orthogonal projection
$\ \pro_{\hat {\cal V}_f}: \rel^{\dim \teich}
\rightarrow \hat {\cal V}_f\ $ that
\begin{equation} \label{euler.theo.invert}
        (\pro_{\hat {\cal V}_f}
        \delta \tpc_f . \vari_r)_{r = 1, \ldots, d}
        =: A \in \rel^{d \times d}
\end{equation}
is invertible, hence after a further change of variable,
we may assume that $\ A = I_d\ $.
In the degenerate case, we further know
\begin{equation} \label{euler.theo.null}
	\langle \delta \tpc_f . \vari_r , e \rangle = 0.
\end{equation}
By (\ref{deg.basis.pm})
\begin{equation} \label{euler.theo.pos}
	\pm \langle \delta^2 \tpc_f(\vari_\pm) , e \rangle \geq 2 \gamma,
	\quad \delta \tpc_f . \vari_\pm = 0,
\end{equation}
for some $\ \gamma > 0\ $.

Next we put for $\ m\ $ large enough
and $\ t_0, \lambda_0\ $ small enough
\begin{displaymath}
	\Phi_m(t,\lambda,\mu) := \tpc(f_{m,t,\lambda,\mu}^* \geu).
\end{displaymath}
Clearly, $\ \Phi_m\ $ is smooth.
We get from (\ref{euler.theo.conv}), (\ref{euler.theo.invert})
(\ref{euler.theo.null}) and (\ref{euler.theo.pos})
for some $\ \Lambda_1 < \infty
\mbox{ and any } 0 < \varepsilon \leq 1\ $ that
\begin{equation} \label{euler.theo.abl-full}
\begin{array}{c}
	\parallel D^2 \Phi_m(t,\lambda) \parallel \leq \Lambda_1, \\
	osc\ D^2 \Phi_m \leq \varepsilon, \\
	\parallel \partial_\lambda \Phi_m(t,\lambda) - I_d \parallel
	\leq \varepsilon \leq 1/2, \\
	\partial_\lambda \Phi_m(t,\lambda)
	\rightarrow (\delta \tpc_f . \vari_r)_{r = 1, \ldots, d},
\end{array}
\end{equation}
in the full rank case,
and writing $\ \Phi_m(t,\lambda,\mu) =
(\Phi_{m,0}(t,\lambda,\mu),\varphi_m(t,\lambda,\mu))\ $
in the degenerate case that
\begin{equation} \label{euler.theo.abl-deg}
\begin{array}{c}
	\parallel D^2 \Phi_m(t,\lambda,\mu) \parallel \leq \Lambda_1, \\
	osc\ D^2 \Phi_m \leq \varepsilon, \\
	\parallel \partial_\lambda \Phi_m(\lambda) - I_d \parallel
	\leq \varepsilon \leq 1/2, \\
	\pm \partial_{\mu_\pm \mu_\pm} \varphi_m(t,\lambda,\mu) \geq \gamma, \\
	|\partial_\mu \Phi_m(t,\lambda,\mu)|,
	|D \varphi_m(t,\lambda,\mu)| \leq \varepsilon, \\
	\partial_\lambda \Phi_m(t,\lambda,\mu)
	\rightarrow (\delta \tpc_f . \vari_r)_{r = 1, \ldots, d}, \\
	\partial_{\mu_\pm \mu_\pm} \varphi_m(t,\lambda,\mu)
	\rightarrow \langle \delta^2 \tpc_f(\vari_\pm) , e \rangle, \\
	\partial_\mu \Phi_m(t,\lambda,\mu),
	D \varphi_m(t,\lambda,\mu) \rightarrow 0,
\end{array}
\end{equation}
all for $\ m \geq m_0\ $ large enough
and $\ |t| \leq t_0, |\lambda|, |\mu| \leq \lambda_0\ $ small enough
or respectively $\ t,\lambda,\mu \rightarrow 0\ $.
We choose $\ \varepsilon, \lambda_0\ $ smaller
to satisfy $\ C \Lambda_1 \varepsilon
+ C \Lambda_1 \lambda_0 \leq \gamma/2\ $.
Moreover choosing $\ m_0\ $ large enough
and $\ t_0\ $ small enough, we can further achieve
\begin{displaymath}
	|D \varphi_m(t,0,0)| \leq \sigma
\end{displaymath}
with $\ C \Lambda_1 \varepsilon
+ C \sigma \lambda_0^{-1} + C \Lambda_1 \lambda_0
\leq \gamma\ $, and
\begin{displaymath}
	|\Phi_m(t,0,0) - \Phi_m(0,0,0)|
	\leq C t_0
	\leq \Lambda_1 \lambda_0^2, \lambda_0/8,
	\gamma \lambda_0^2/32.
\end{displaymath}
By Proposition \ref{ana.prop}
there exist $\ |\lambda_{m,r}(t)|, |\mu_{m,\pm}(t)|,
|\tilde \lambda_{m,r}(t)|, |\tilde \mu_{m,\pm}(t)| \leq \lambda_0\ $
with $\ \mu_{m,+}(t) \mu_{m,-}(t) = 0,
\tilde \mu_{m,+}(t) \tilde \mu_{m,-}(t) = 0\ $ and
\begin{equation} \label{euler.theo.solu}
	\Phi_m(t,\lambda_m(t),\mu_m(t)) = \Phi_m(0,0,0)
	= \Phi_m(t,\tilde \lambda_m(t),\tilde \mu_m(t)),
\end{equation}
which means
\begin{equation} \label{euler.theo.teich}
\begin{array}{c}
	\tp \Big(  (f_m + t \vari + \lambda_{m,r}(t) \vari_r
	+ \mu_{m,\pm}(t) \vari_\pm)^* \geu \Big)
	= \tp(f_m^* \geu) = \\
	= \tp \Big(  (f_m + t \vari + \tilde \lambda_{m,r}(t) \vari_r
	+ \tilde \mu_{m,\pm}(t) \vari_\pm)^* \geu \Big),
\end{array}
\end{equation}
and
\begin{equation} \label{euler.theo.teich-sign}
\begin{array}{c}
	\mu_{m,+}(t) \tilde \mu_{m,+}(t) \leq 0, \tilde \mu_{m,-}(t) = 0,
	\quad \mbox{if } \mu_{m,+}(t) \neq 0, \\
	\mu_{m,-}(t) \tilde \mu_{m,-}(t) \leq 0, \tilde \mu_{m,+}(t) = 0,
	\quad \mbox{if } \mu_{m,-}(t) \neq 0.
\end{array}
\end{equation}
By (\ref{euler.theo.abl-full}), (\ref{euler.theo.abl-deg})
and (\ref{euler.theo.solu}),
we get from a Taylor expansion of the smooth function
$\ \Phi_m \mbox{ at } 0\ $ that
\begin{displaymath}
	|t \partial_t \Phi_m(0) + \lambda_m(t) \partial_\lambda \Phi_m(0)
	+ \mu_m(t) \partial_\mu \Phi_m(0) +
\end{displaymath}
\begin{displaymath}
	+ \frac{1}{2} (t,\lambda_m(t),\mu_m(t))^T
	D^2 \Phi_m(0) (t,\lambda_m(t),\mu_m(t))| \leq
\end{displaymath}
\begin{displaymath}
	\leq C \varepsilon (|t|^2 + |\lambda_m(t)|^2 + |\mu_m(t)|^2 ),
\end{displaymath}
hence, as $\ \mu_{m,+}(t) \mu_{m,-}(t) = 0\ $,
\begin{displaymath}
	|t \partial_t \Phi_m(0) + \lambda_m(t) \partial_\lambda \Phi_m(0)
	+ \mu_{m,\pm}(t) \partial_{\mu_\pm} \Phi_m(0)
	+ \frac{1}{2} \mu_{m,\pm}(t)^2\ \partial_{\mu_\pm \mu_\pm} \Phi_m(0)| \leq
\end{displaymath}
\begin{displaymath}
	\leq C \Lambda_1 \varepsilon^{-1} (|t|^2 + |\lambda_m(t)|^2 )
	+ C \varepsilon |\mu_{m,\pm}(t)|^2 ).
\end{displaymath}
Passing to the limit $\ m \rightarrow \infty\ $,
we get for subsequences $\ \lambda_m(t) \rightarrow \lambda(t),
\mu_{m,\pm}(t) \rightarrow \mu(t)
\mbox{ with } |\lambda(t)|, |\mu(t)| \leq \lambda_0\ $
and by (\ref{euler.theo.abl-full}) and (\ref{euler.theo.abl-deg})
\begin{displaymath}
	|t \delta \tpc_f . \vari + \lambda_r(t) \delta \tpc_f . \vari_r
	+ \frac{1}{2} \mu_{\pm}(t)^2
	\ \delta^2 \tpc_f(\vari_\pm)| \leq
\end{displaymath}
\begin{equation} \label{euler.theo.taylor}
	\leq C \Lambda_1 \varepsilon^{-1} (|t|^2 + |\lambda(t)|^2 )
	+ C \varepsilon |\mu_\pm(t)|^2 ).
\end{equation}
In the degenerate case, we recall $\ \delta \tpc_f . V_{(r)}
\in \hat {\cal V}_f \perp e\ $,
hence by (\ref{euler.theo.pos}) and (\ref{euler.theo.taylor})
\begin{displaymath}
	\gamma \mu_\pm(t)^2
	\leq \frac{1}{2} \mu_\pm(t)^2 \langle
	\delta^2 \tpc_f(\vari_\pm) , e \rangle \leq
\end{displaymath}
\begin{displaymath}
	\leq |t \delta \tpc_f . \vari + \lambda_r(t) \delta \tpc_f . \vari_r
	+ \frac{1}{2} \mu_{m,\pm}(t)^2
	\ \delta^2 \tpc_f(\vari_\pm)| \leq	
\end{displaymath}
\begin{displaymath}
	\leq C \Lambda_1 \varepsilon^{-1} (|t|^2 + |\lambda(t)|^2 )
	+ C \varepsilon |\mu_\pm(t)|^2 )
\end{displaymath}
and for $\ C \varepsilon \leq \gamma/2\ $ small enough
\begin{equation} \label{euler.theo.mu-lin-aux}
	|\mu_\pm(t)| \leq C (|t| + |\lambda(t)|).
\end{equation}
Then again by (\ref{euler.theo.taylor})
\begin{equation} \label{euler.theo.square}
	|t \delta \tpc_f . \vari
	+ \lambda_r(t) \delta \tpc_f . \vari_r|
	\leq C (|t|^2 + |\lambda(t)|^2 ).
\end{equation}
In the full rank case, we have $\ \mu_\pm = 0\ $,
and (\ref{euler.theo.square}) directly follows
from (\ref{euler.theo.taylor}).
As $\ (\delta \tpc_f . \vari_r)_{r = 1, \ldots, d}\ $
are linearly independent by (\ref{euler.theo.invert}),
we continue
\begin{displaymath}
	|\lambda(t)| \leq C (|t| + |\lambda(t)|^2)
	\leq C |t| + C \lambda_0 |\lambda(t)|,
\end{displaymath}
hence for $\ \lambda_0\ $ small enough
$\ |\lambda(t)| \leq C |t|\ $.
We get from (\ref{euler.theo.mu-lin-aux})
\begin{equation} \label{euler.theo.mu-lin}
	|\mu_\pm(t)| \leq C |t|
\end{equation}
and from (\ref{euler.theo.square})
\begin{displaymath}
	|t \delta \tpc_f . \vari
	+ \lambda_r(t) \delta \tpc_f . \vari_r|
	\leq C |t|^2.
\end{displaymath}
Therefore $\ \lambda\ $ is differentiable at $\ t = 0\ $ with
\begin{equation} \label{euler.theo.lambda-diff}
	\lambda'_r(0) \delta \tpc_f . \vari_r
	= - \delta \tpc_f . \vari.
\end{equation}
Writing for the inverse $\ A^{-1} = (b_{rs})_{r,s = 1, \ldots,d}\ $
in (\ref{euler.theo.invert}), we continue with (\ref{full.vari-form})
\begin{displaymath}
	\lambda'_r(0) = - b_{rs} \langle \delta \tpc_f . \vari , e_s \rangle =
\end{displaymath}
\begin{displaymath}
	= \sum \limits_{\sigma=1}^{\dim\ \teich}
	2 \int \limits_\sur g^{ik} g^{jl}
	\langle A_{ij}^0 , \vari \rangle
	q^\sigma_{kl}(\gp) \d \mu_g
	\ \langle d \tpc_\gp . q^\sigma(\gp) , b_{rs} e_s \rangle,
\end{displaymath}
hence putting
\begin{displaymath}
	q^r_{kl} := \sum \limits_{\sigma=1}^{\dim\ \teich}
	2 q^\sigma _{kl}(\gp)
	\langle d \tpc_\gp . q^\sigma(\gp) , b_{rs} e_s \rangle
	\in S_2^{TT}(\gp),
\end{displaymath}
we get
\begin{equation} \label{euler.theo.lambda-abl}
	\lambda'_r(0)
	= \int \limits_\sur g^{ik} g^{jl}
	\langle A_{ij}^0 , \vari \rangle
	q^r_{kl} \d \mu_g.
\end{equation}
Moreover
\begin{displaymath}
	f_m + t \vari + \lambda_{m,r}(t) \vari_r
	+ \mu_{m,\pm}(t) \vari_\pm
	\rightarrow f + t \vari + \lambda_r(t) \vari_r
	+ \mu_\pm(t) \vari_\pm
\end{displaymath}
strongly in $\ W^{2,2}(\sur)\ $
and weakly$^* \mbox{ in } W^{1,\infty}(\sur)\ $,
hence recalling (\ref{euler.theo.teich})
\begin{displaymath}
	\W(f + t \vari + \lambda_r(t) \vari_r
	+ \mu_\pm(t) \vari_\pm)
	\leftarrow \W(f_m + t \vari + \lambda_{m,r}(t) \vari_r
	+ \mu_{m,\pm}(t) \vari_\pm) \geq
\end{displaymath}
\begin{displaymath}
	\geq \mini(\tp(f_m^* \geu))
	\geq  \W(f_m) - \varepsilon_m
	\rightarrow \W(f).
\end{displaymath}
Since $\ (t,\lambda,\mu) \mapsto \W(f + t \vari
+ \lambda_r \vari_r + \mu_\pm \vari_\pm)\ $ is smooth,
we get again by a Taylor expansion,
(\ref{euler.theo.mu-lin}) and  (\ref{euler.theo.lambda-diff}) 
\begin{displaymath}
	0 \leq \W(f + t \vari + \lambda_r(t) \vari_r
	+ \mu_\pm(t) \vari_\pm) - \W(f) =
\end{displaymath}
\begin{displaymath}
	= t \delta \W_f . \vari + \lambda_r(t) \delta \W_f . \vari_r
	+ \mu_\pm(t) \delta \W_f . \vari_\pm + O(|t|^2).
\end{displaymath}
As $\ \mu_+(t) \mu_-(t) = 0\ $
and by (\ref{euler.theo.teich-sign}),
we can adjust the sign of $\ \mu_\pm(t)\ $
according to the sign of $\ \delta \W_f . \vari_\pm\ $
and improve to
\begin{displaymath}
	0 \leq t \delta \W_f . \vari
	+ \lambda_r(t) \delta \W_f . \vari_r + O(|t|^2).
\end{displaymath}
Differentiating by $\ t \mbox{ at } t = 0\ $,
we conclude from (\ref{euler.theo.lambda-diff})
and (\ref{euler.theo.lambda-abl})
\begin{displaymath}
	\delta \W_f . \vari
	= - \lambda'_r(0) \delta \W_f . \vari_r
	= - \int \limits_\sur g^{ik} g^{jl}
	\langle A_{ij}^0 , \vari \rangle
	q^r_{kl}\ \delta \W_f . \vari_r \d \mu_g,
\end{displaymath}
hence putting $\ q_{kl} :=
q^r_{kl}\ \delta \W_f . \vari_r \in S_2^{TT}(\gp)\ $,
we get
\begin{equation} \label{euler.theo.equ-weak}
    \delta \W_f . \vari
    = \int \limits_\sur g^{ik} g^{jl}
    \langle A_{ij}^0 , \vari \rangle
    q_{kl} \d \mu_g
    \quad \mbox{for all } V \in C^\infty(\sur,\rel^n).
\end{equation}
As $\ f \in W^{2,2} \cap C^{1,\alpha}\ $
by Proposition \ref{dec.prop},
we can write $\ f\ $ as a graph,
more precisely for any $\ x_0 \in \sur\ $
there exists a neighbourhood
$\ U(x_0) \mbox{ of } x_0\ $ such that
after a translation, rotation
and a homothetie, which leaves $\ \W\ $
as conformal transformation invariant,
there is a $\ (W^{2,2} \cap C^{1,\alpha})-\mbox{inverse chart }
\diffeo: B_1(0) \stackrel{\approx}{\longrightarrow}U(x_0), \diffeo(0) = x_0\ $,
with $\ \tilde f(y) := (f \circ \diffeo)(y) = (y,u(y))\ $
for some $\ u \in (W^{2,2} \cap C^{1,\alpha})(B_1(0),\rel^{n-2})\ $.
Moreover, we may assume $\ |u|, |D u| \leq 1\ $
and from (\ref{dec.prop.dec}) that
\begin{equation} \label{euler.theo.dec}
	\int \limits_{B_\varrho} |D^2 u|^2 \d \Lt \leq C \varrho^{2 \alpha}
	\quad \mbox{for any Ball } B_\varrho.
\end{equation}
We calculate the square integral of the second fundamental form
for a graph as
\begin{displaymath}
	{\cal A}(u) := \int \limits_{B_1(0)} |A_{\tilde f}|^2 \d \mu_{\tilde f}
	= \int \limits_{B_1(0)}
	(\delta_{rs} - d_{rs}) g^{ij} g^{kl}
	\partial_{ik} u^r \partial_{jl} u^s \sqrt{g} \d \Lt,
\end{displaymath}
where $\ g_{ij} := \delta_{ij} + \partial_i u \partial_j u,
(g^{ij}) = (g_{ij})^{-1}, d_{rs} :=
g^{kl} \partial_k u^r \partial_l u^s\ $,
see \bcite{sim.will} p. 310.

By the Gau\ss-Bonnet Theorem in (\ref{intro.gauss}),
we see for any $\ v \in C_0^\infty(B_1(0),\rel^{n-2})\ $
that $\ 4 \W(graph(u+v)) = {\cal A}(u+v) + 8 \pi (1-p)\ $,
hence from (\ref{euler.theo.equ-weak})
\begin{displaymath}
	\delta {\cal A}_u . v
	= \int \limits_{B_1(0)} (g^{ik} g^{jl} - \frac{1}{2} g^{ij} g^{kl})
	\tilde q_{kl}
	\langle \partial_{ij} u - \Gamma_{ij}^m \partial_m u , v \rangle \d \mu_g,
\end{displaymath}
where $\ \tilde q = \diffeo^* q,
\Gamma_{ij}^m = g^{mk} \langle \partial_{ij} u ,
\partial_k u \rangle\ $.
From this we conclude that
\begin{equation} \label{boundary.equ}
	\partial_{jl} (2 a^{ijkl}_{rs} \partial_{ik} u^s)
	- \partial_j \Big( (\partial_{\partial_j u^r} a^{imkl}_{ts})
	\partial_{ik} u^s \partial_{ml} u^t \Big)
	= b_{rs}^{ijkl}(\partial u) \tilde q_{kl} \partial_{ij} u^s
\end{equation}
weakly for testfunctions $\ v \in C^\infty_0(B_1(0),\rel^{n-2})\ $,
where
\begin{displaymath}
\begin{array}{c}
	a^{ijkl}_{rs}(D u) := (\delta_{rs} - d_{rs})
	g^{ij} g^{kl} \sqrt{g}, \\
	b_{rs}^{ijkl}(D u)
	= (g^{ik} g^{jl} - \frac{1}{2} g^{ij} g^{kl})
	(\delta_{rs} - d_{rs}) \sqrt{g}.
\end{array}
\end{displaymath}
Then we conclude from \bcite{sim} Lemma 3.2
and (\ref{euler.theo.dec})
that $\ u \in (W^{3,2}_{loc} \cap C^{2,\alpha}_{loc})(B_1(0))\ $.

Full Regularity is now obtained by \bcite{adn1}, \bcite{adn2}.
First we conclude by finite differences
that $\ u \in W^{4,2}_{loc}(B_1(0))\ $ and
\begin{displaymath}
	2 a^{ijkl}_{rs}(Du) \partial_{ijkl} u^s
	+ \tilde b_r(D u,D^2u) * (1 + D^3 u)
	+ b_r(Du) * \tilde q(.) * D^2 u = 0
	\quad \mbox{strongly in } B_1(0),
\end{displaymath}
with $\ a^{ijkl}_{rs}, \tilde b_r, b_r\ $
are smooth in $\ Du \mbox{ and } D^2 u\ $,
whereas $\ \tilde q = \diffeo^* q \in (W^{1,2}
\cap C^{0,\alpha})(B_1(0))\ $.
As $\ W^{4,2} \hookrightarrow W^{3,p}
\mbox{ for all } 1 \leq p < \infty\ $,
we see $\ u \in W^{4,p}_{loc}(B_1(0))
\hookrightarrow C^{3,\alpha}_{loc}(B_1(0))\ $
and then $\ u \in C^{4,\alpha}_{loc}(B_1(0))\ $.

Now we proceed by induction
assuming $\ u, \tilde f \in C^{k,\alpha}_{loc}(B_1(0))
\mbox{ for some } k \geq 4\ $.
We see $\ \tilde g := \tilde f^* \geu \in C^{k-1,\alpha}_{loc}\ $,
hence we get locally conformal $\ C^{k,\alpha}-\mbox{charts }
\varphi: U \subseteq B_1(0) \stackrel{\approx}{\longrightarrow}
\Omega \subseteq \rel^2 \mbox{ with } \varphi^{-1,*} \tilde g = e^{2 v} \geu\ $.
On the other hand, as $\ \gp\ $ is smooth,
there exists a smooth conformal chart $\ \psi: U(x_0)
\stackrel{\approx}{\longrightarrow} \Omega_0 \subseteq \rel^2
 \mbox{ with } \psi^{-1,*} \gp = e^{2 u_0} \geu\ $,
when choosing $\ U(x_0)\ $ small enough.
As $\ g = f^* \geu = e^{2u} \gp\ $,
we see that $\ \psi \circ \diffeo \circ \varphi^{-1}:
\Omega \rightarrow \Omega_0\ $ is a regular conformal mapping
with respect to standard euclidian metric,
in particular holomorphic or anti-holomorphic,
hence smooth.
We conclude that $\ \diffeo \in C^{k,\alpha}_{loc}(B_1(0))
\mbox{ and } \tilde q = \diffeo^* q \in C^{k-1,\alpha}_{loc}(B_1(0))\ $,
as $\ q \in S_2^{TT}(\gp)\ $ is smooth.
Then we conlude $\ u \in C^{k+1,\alpha}_{loc}(B_1(0))\ $
and by induction $\ u, \tilde f, \diffeo
\mbox{ and } f = \tilde f \circ \diffeo^{-1}\ $ are smooth.

In \bcite{kuw.schae.will} \S2,
the first variation of the Willmore functional with a different factor
was calculated for variations $\ V\ $ to be 
\begin{displaymath}
	\delta \W(f) . V
	:= \frac{d}{dt} \W(f + t V)
	= \int \limits_\sur
	\frac{1}{2} \langle \Delta_g \mean 
	+ Q(A^0) \mean , V \rangle \d \mu_g,
\end{displaymath}
and we obtain from (\ref{euler.theo.equ-weak})
\begin{displaymath}
        \Delta_g \mean + Q(A^0) \mean = 2 g^{ik} g^{jl} A_{ij}^0 q_{kl}
	\quad \mbox{on } \sur,
\end{displaymath}
which is (\ref{euler.theo.euler})
up to a factor for $\ q\ $.
(\ref{euler.theo.minimizer}) was already
obtained in (\ref{ele.strong.minimizer}).
This concludes the proof of the theorem.
\proof
As a corollary we get minimizers
under fixed Teichm\"uller or conformal class,
when the infimum is smaller
than the bound $\ \W_{n,p}\ $ in (\ref{intro.moeb.energ}).

\begin{corollaryth} \label{euler.teich}

Let $\ \sur\ $ be a closed orientable surface
of genus $\ p \geq 1\ $
and $\ \tau_0 \in \teich\ $ satisfying
\begin{displaymath}
	\mini(\tau_0) <  \W_{n,p}
\end{displaymath}
where $\ \W_{n,p}\ $
is defined in (\ref{intro.moeb.energ})
and $\ n = 3,4\ $.

Then there exists a smooth immersion
$\ f: \sur \rightarrow \rel^n\ $
which minimizes the Willmore energy
in the fixed Teichm\"uller class $\ \tau_0 = \tp(f^* \geu)\ $
\begin{displaymath}
	\W(f) = \mini(\tau_0).
\end{displaymath}
Moreover $\ f\ $ satisfies the Euler-Lagrange equation
\begin{displaymath}
        \Delta_g \mean + Q(A^0) \mean = g^{ik} g^{jl} A_{ij}^0 q_{kl}
	\quad \mbox{on } \sur,
\end{displaymath}
where $\ q\ $ is a smooth transverse traceless
symmetric 2-covariant tensor with respect
to $\ g = f^* \geu\ $.
\end{corollaryth}
{\pr Proof:} \\
We select a minimizing sequence of smooth
immersions $\ f_m: \sur \rightarrow \rel^n
\mbox{ with } \tp(f_m^* \geu) = \tau_0\ $
\begin{equation} \label{euler.teich.mini}
	\W(f_m) \rightarrow \mini(\tau_0).
\end{equation}
We may assume that $\ \W(f_m) \leq \W_{n,p} - \delta
\mbox{ for some } \delta > 0\ $.
Replacing $\ f_m \mbox{ by }
\moeb_m \circ f_m \circ \diffeo_m\ $
for suitable M\"obius transformations $\ \moeb_m\ $
and diffeomorphisms $\ \diffeo_m \mbox{ of } \sur\ $
homotopic to the identity,
which does neither change the Willmore energy
nor the projection into the Teichm\"uller space,
we may further assume by Proposition \ref{direct.prop}
that $\ f_m \rightarrow f \mbox{ weakly in } W^{2,2}(\sur)\ $
and satisfies (\ref{direct.prop.teich}) - (\ref{direct.prop.imm}).
Since (\ref{euler.theo.mini})
is implied by (\ref{euler.teich.mini}),
Theorem \ref{euler.theo}) yields
that $\ f\ $ is a smooth immersion
which minimizes the Willmore energy
in the fixed Teichm\"uller class $\ \tau_0 = \tp(f^* \geu)\ $
and satisfies the above Euler-Lagrange equation.
\proof

\begin{corollaryth} \label{euler.conf}

Let $\ \sur\ $ be a closed Riemann surface
of genus $\ p \geq 1\ $
with
\begin{displaymath}
	\inf \{ \W(f)\ |
	\ f: \sur \rightarrow \rel^n
	\mbox{ conformal immersion}\ \}
	< \W_{n,p},
\end{displaymath}
where $\ \W_{n,p}\ $
is defined in (\ref{intro.moeb.energ})
and $\ n = 3,4\ $.

Then there exists a smooth conformal immersion
$\ f: \sur \rightarrow \rel^n\ $
which minimizes the Willmore energy
in the set of all conformal immersions.
Moreover $\ f\ $ satisfies
the Euler-Lagrange equation
\begin{equation} \label{euler.conf.euler}
	\Delta_g \mean + Q(A^0) \mean = g^{ik} g^{jl} A_{ij}^0 q_{kl}
	\quad \mbox{on } \sur,
\end{equation}
where $\ q\ $ is a smooth transverse traceless
symmetric 2-covariant tensor with respect
to the Riemann surface $\ \sur\ $,
that is with respect to $\ g = f^* \geu\ $.
\end{corollaryth}
{\pr Proof:} \\
We select select a smooth conformal metric $\ g_0\ $
for the Riemann surface $\ \sur\ $
and put $\ \tau_0 := \tp(g_0)\ $.
By invariance, we see
\begin{displaymath}
	\mini(\tau_0) = \inf \{ \W(f)\ |
	\ f: \sur \rightarrow \rel^n
	\mbox{ conformal immersion}\ \}
	< \W_{n,p}.
\end{displaymath}
Therefore by Corollary \ref{euler.teich},
there exists a smooth immersion $\ f: \sur \rightarrow \rel^n\ $
which minimizes the Willmore energy
in the fixed Teichm\"uller class $\ \tau_0 = \tp(f^* \geu)\ $
and satisfies the above Euler-Lagrange equation.
Moreover there exists a diffeomorphism $\ \diffeo\mbox{ of } \sur\ $
homotopic to the identity
such that $\ (f \circ \diffeo)^* \geu
\mbox{ is conformal to } g_0\ $.
Then $\ \tilde f := f \circ \diffeo\ $
is a smooth conformal immersion
of the Riemann surface $\ \sur\ $
which minimizes the Willmore energy
in the set of all conformal immersions
and moreover satisfies
the Euler-Lagrange equation.
\proof



\ \\
{\LARGE \bf Appendix}

\begin{appendix}
\renewcommand{\theequation}{\mbox{\Alph{section}.\arabic{equation}}}


\setcounter{equation}{0}
\section{Conformal factor} \label{conf}

\begin{lemma} \label{conf.lemma}

Let $\ \sur\ $ be a closed, orientable surface of genus $\ p \geq 1\ $,
$\ g_0\ $ a given smooth metric on $\ \sur\ $,
$\ x_1, \ldots, x_\M \in \sur\ $
with charts $\ \varphi_k: U(x_k) \stackrel{\approx}{\longrightarrow}
B_1(0), \varphi_k(x_k) = 0,
U_\varrho(x_k) := \varphi_k^{-1}(B_\varrho(0))
\mbox{ for } 0 < \varrho \leq 1,\ $
\begin{equation} \label{conf.lemma.metric}
        \Lambda^{-1} \geu \leq (\varphi_k^{-1})^* g_0 \leq \Lambda \geu
        \quad \mbox{for } k = 1, \ldots, \M.
\end{equation}
Let $\ f: \sur \rightarrow \rel^n\ $ be a smooth immersion
with $\ g := f^* \geu = e^{2u} \gp\ $
for some unit volume constant curvature metric $\ \gp\ $ and
\begin{equation} \label{conf.lemma.ub}
\begin{array}{c}
	\Lambda^{-1} g_0 \leq g \leq \Lambda g_0, \\
	\parallel u \parallel_{L^\infty(\sur)},
	\int \limits_\sur |K_g| \d \mu_g \leq \Lambda,
\end{array}
\end{equation}
and $\ \tilde f: \sur \rightarrow \rel^n\ $ a smooth immersion
with $\ \tilde g := \tilde f^* \geu = e^{2 \tilde u} \tgp\ $
for some unit volume constant curvature metric $\ \tgp\ $,
\begin{equation} \label{conf.lemma.metric2}
\begin{array}{c}
	\Lambda^{-1} g_0 \leq \tilde g \leq \Lambda g_0, \\
	\int \limits_\sur |K_{\tilde g}| \d \mu_{\tilde g} \leq \Lambda,
\end{array}
\end{equation}
$\ supp(f - \tilde f) \subseteq \cup_{k=1}^\M U_{1/2}(x_k)\ $ and
\begin{equation} \label{conf.lemma.fund}
	\int \limits_{U_1(x_k)} |\tilde A|^2 \d \mu_{\tilde g}
	\leq \varepsilon_0(n)
	\quad \mbox{for } k = 1, \ldots, \M
\end{equation}
for some universal $\ 0 < \varepsilon_0(n) < 1\ $.
Then
\begin{equation} \label{conf.lemma.tu}
	\parallel \tilde u \parallel_{L^\infty(\sur)},
	\parallel \nabla
	\tilde u \parallel_{L^2(\sur,\tilde g)}
	\leq C(\sur,g_0,\Lambda,p).
\end{equation}
\end{lemma}
{\pr Proof:} \\
We know
\begin{equation} \label{conf.lemma.equ}
	- \Delta_g u + K_p e^{-2u} = K_g,
	\quad
	- \Delta_{\tilde g} \tilde u + K_p e^{-2 \tilde u}
	= K_{\tilde g}
	\quad \mbox{on } \sur.
\end{equation}
Observing $\ \int_\sur |-K_p e^{-2 \tilde u}| \d \mu_{\tilde g}
= - K_p = 4 \pi (p-1)\ $
and multiplying (\ref{conf.lemma.equ}) by $\ \tilde u - \tilde \lambda
\mbox{ for any } \tilde \lambda \in \rel\ $,
we get recalling (\ref{conf.lemma.metric2})
\begin{displaymath}
	c_0(\Lambda) \int \limits_\sur |\nabla \tilde u|_{g_0}^2
	\d \mu_{g_0}
	\leq \int \limits_\sur |\nabla
	\tilde u|_{\tilde g}^2 \d \mu_{\tilde g} =
\end{displaymath}
\begin{displaymath}
	= \int \limits_\sur (-K_p e^{-2 \tilde u} + K_{\tilde g})
	(\tilde u - \tilde \lambda) \d \mu_{\tilde g}
	\leq C(\Lambda,p) \parallel \tilde u - \tilde \lambda
	\parallel_{L^\infty(\sur)},
\end{displaymath}
hence for $\ \tilde \lambda = \mint{\sur} \tilde u \d \mu_{g_0}\ $
by Poincar\'e inequality
\begin{equation} \label{conf.lemma.grad}
	\parallel \tilde u - \tilde \lambda \parallel_{L^2(\sur,g_0)}
	\leq C(\sur,g_0) \parallel \nabla \tilde u \parallel_{L^2(\sur,g_0)}
	\leq C(\sur,g_0,\Lambda,p) \sqrt{osc_\sur \tilde u}.
\end{equation}
Next by the uniformization theorem, see \bcite{far.kra} Theorem IV.4.1,
we can parametrize
$\ \tilde f \circ \varphi_k^{-1}: B_1(0) \rightarrow \rel^n\ $
conformally with respect to the euclidean metric on $\ B_1(0)\ $,
possibly after replacing $\ B_1(0)\ $ by a slightly smaller ball.
Then by \bcite{muell.sver} Theorem 4.2.1
for $\ \varepsilon_0(n)\ $ small enough,
there exist $\ v_k \in C^\infty(U_1(x_k))\ $ with
\begin{equation} \label{conf.lemma.v}
\begin{array}{c}
	- \Delta_{\tilde g} v_k = K_{\tilde g}
	\quad \mbox{on } U_1(x_k), \\
	\parallel v_k \parallel_{L^\infty(U_1(x_k))}
	\leq C_n \int \limits_{U_1(x_k)} |\tilde A|^2 \d \mu_{\tilde g}
	\leq C_n \varepsilon_0(n) \leq 1
\end{array}
\end{equation}
for $\ \varepsilon_0(n)\ $ small enough.
We get
\begin{equation} \label{conf.lemma.equ.uv}
	\left.
	\begin{array}{c}
		- \Delta_{\tilde g} (\tilde u - v_k)
		= - K_p e^{-2 \tilde u} \geq 0, \\
		- \Delta_{\tilde g} (\tilde u - v_k)
		+ K_p ( e^{-2 \tilde u} - e^{-2 v_k})
		= - K_p e^{-2 v_k} \geq 0,
	\end{array}
	\right\} \quad \mbox{on } U_1(x_k)
\end{equation}
for $\ k = 1, \ldots, \M\ $,
and, as $\ g = \tilde g \mbox{ on } \sur - \cup_{k=1}^\M
\overline{U_{1/2}(x_k)}\ $,
\begin{equation} \label{conf.lemma.equ.uu}
	\left.
	\begin{array}{c}
 		- \Delta_g (\tilde u - u)
 		+ K_p ( e^{-2 \tilde u} - e^{-2 u}) = 0, \\
		- \Delta_g (\tilde u - u)_+ \leq 0,
	\end{array}
	\right\}
	\quad \mbox{on } \sur - \cup_{k=1}^\M \overline{U_{1/2}(x_k)}.
\end{equation}
Therefore $\ \tilde u - u\ $ cannot have
positive interior maxima
nor negative interior minima
in $\ \sur - \cup_{k=1}^\M \overline{U_{1/2}(x_k)}\ $
as $\ K_p \leq 0\ $ by standard maximum principle,
hence putting $\ \Gamma := \cup_{k=1}^\M \partial U_{3/4}(x_k)\ $
we get
\begin{displaymath}
	\sup \limits_{\sur - \cup_{k=1}^\M \overline{U_{3/4}(x_k)}}
	(\tilde u - u)_\pm
	= \max \limits_\Gamma (\tilde u - u)_\pm.
\end{displaymath}
From above, we see from (\ref{conf.lemma.metric}),
(\ref{conf.lemma.metric2}), (\ref{conf.lemma.v}), (\ref{conf.lemma.equ.uv}),
\begin{equation} \label{conf.lemma.equ.uv-local}
	0 \leq - \partial_i \Big( \tilde g^{ij} \sqrt{\tilde g}
	\ \partial_j (\tilde u - v_k) \Big)
	+ K_p \sqrt{\tilde g} ( e^{-2 \tilde u} - e^{-2 v_k})
	= - K_p e^{-2 v_k} \sqrt{\tilde g}
	\leq C(\Lambda,p),
\end{equation}
hence using \bcite{gil.tru} Theorem 8.16
\begin{displaymath}
	\sup \limits_{U_{3/4}(x_k)} (\tilde u - v_k)_\pm
	\leq \max \limits_\Gamma (\tilde u - v_k)_\pm
	+ C(\Lambda,p).	
\end{displaymath}
Together, we get from (\ref{conf.lemma.ub})
and (\ref{conf.lemma.v})
\begin{equation} \label{conf.lemma.infty-omega}
	\max \limits_\sur \tilde u_\pm
	\leq \max \limits_\Gamma \tilde u_\pm + C(\Lambda,p).
\end{equation}
From (\ref{conf.lemma.metric}), (\ref{conf.lemma.ub}),
(\ref{conf.lemma.metric2}), and $\ \gp\ $ having  unit volume,
we get
\begin{equation} \label{conf.lemma.meas}
	c_0(\Lambda) \leq \mu_{g}(\sur), \mu_{g_0}(\sur),
	\mu_{\tilde g}(\sur) \leq C(\Lambda).
\end{equation}
Now if $\ \min_\sur \tilde u \leq - C(\Lambda,p)\ $,
there exists $\ x \in \Gamma \mbox{ with }
\tilde u(x) \leq \min_\sur \tilde u +  C(\Lambda,p) \leq 0\ $.
As $\ \tilde u - v_k \geq \min_\sur \tilde u - 1 =: \lambda\ $,
we get identifying $\ \varphi_k: U_1(x_k) \cong B_1(0)\ $
by the weak Harnack inequality, see \bcite{gil.tru} Theorem 8.18,
from (\ref{conf.lemma.equ.uv-local})
\begin{displaymath}
	\parallel \tilde u - v_k - \lambda \parallel_{L^2(B_{1/8}(x))}
	\leq C(\Lambda) \inf \limits_{B_{1/8}(x)}
	(\tilde u - v_k - \lambda) \leq
\end{displaymath}
\begin{displaymath}
	\leq C(\Lambda) (\tilde u(x) - \min \limits_\sur \tilde u + 2)
	\leq C(\Lambda,p),
\end{displaymath}
and
\begin{displaymath}
	\parallel \tilde u - \min_\sur \tilde u
	\parallel_{L^2(B_{1/8}(x))}
	\leq C(\Lambda,p).
\end{displaymath}
We see from (\ref{conf.lemma.metric})
and (\ref{conf.lemma.grad})
\begin{displaymath}
	c_0 |\tilde \lambda - \min_\sur \tilde u|
	\leq \parallel \tilde \lambda - \min_\sur \tilde u
	\parallel_{L^2(B_{1/8}(x))} \leq
\end{displaymath}
\begin{displaymath}
	\leq C(\Lambda) \parallel \tilde u - \tilde \lambda
	\parallel_{L^2(\sur,g_0)}
	+ \parallel \tilde u - \min_\sur \tilde u
	\parallel_{L^2(B_{1/8}(x))}
	\leq C(\sur,g_0,\Lambda,p) (1 + \sqrt{osc_\sur \tilde u}).
\end{displaymath}
Using (\ref{conf.lemma.meas}),
we have proved
\begin{displaymath}
	\min_\sur \tilde u \leq - C(\Lambda,p)
	\Longrightarrow
\end{displaymath}
\begin{equation} \label{conf.lemma.min}
	\parallel \tilde u - \min_\sur \tilde u \parallel_{L^2(\sur,g_0)}
	\leq C(\sur,g_0,\Lambda,p) (1 + \sqrt{osc_\sur \tilde u}).
\end{equation}
If further $\ \min_\sur \tilde u \ll - C(\sur,g_0,\Lambda,p)
(1 + \sqrt{osc_\sur \tilde u})\ $,
we put $\ A := [ \min_\sur \tilde u \leq \tilde u
\leq \min_\sur \tilde u]\ $
and see $\ \tilde u - \min_\sur \tilde u
\geq |\min_\sur \tilde u / 2| \mbox{ on } \sur - A\ $
and using (\ref{conf.lemma.metric2})
\begin{displaymath}
	\frac{1}{2} |\min_\sur \tilde u|\ \mu_{\tilde g}(\sur - A)
	\leq \int \limits_{\sur - A}
	|\tilde u - \min \limits_\sur \tilde u|
	\d \mu_{\tilde g} \leq
\end{displaymath}
\begin{displaymath}
	\leq C(\sur,g_0,\Lambda,p) \int \limits_\sur
	|\tilde u - \min \limits_\sur \tilde u| \d \mu_{g_0}
	\leq C(\sur,g_0,\Lambda,p) (1 + \sqrt{osc_\sur \tilde u}),
\end{displaymath}
hence
\begin{displaymath}
	\mu_{\tilde g}(\sur - A) \leq \frac{C(\sur,g_0,\Lambda,p)
	(1 + \sqrt{osc_\sur \tilde u})}{|\min \limits_\sur \tilde u|}
	\leq c_0(\Lambda)/2,
\end{displaymath}
if $\ |\min_\sur \tilde u| \gg C(\sur,g_0,\Lambda,p)
(1 + \sqrt{osc_\sur \tilde u})\ $ is large enough.
This yields using (\ref{conf.lemma.meas})
and $\ \mu_\tgp(\sur) = 1\ $
\begin{displaymath}
	c_0(\Lambda)/2 \leq \mu_{\tilde g}(A)
	= \int \limits_A e^{2 \tilde u} \d \mu_\tgp
	\leq \mu_\tgp(\sur) \exp(\min \limits_\sur \tilde u)
	= \exp(\min \limits_\sur \tilde u),
\end{displaymath}
and we conclude
\begin{equation} \label{conf.lemma.min2}
	\min_\sur \tilde u \geq - C(\sur,g_0,\Lambda,p)
	(1 + \sqrt{osc_\sur \tilde u}).
\end{equation}
In the same way as above,
if $\ \max_\sur \tilde u \geq - C(\Lambda,p)\ $,
we get from (\ref{conf.lemma.infty-omega}) that
there exists $\ x \in \Gamma \mbox{ with }
\tilde u(x) \geq \max_\sur \tilde u -  C(\Lambda,p) \geq 0\ $.
As $\ \tilde u - v_k \leq \max_\sur \tilde u + 1 =: \lambda\ $,
we get by the weak Harnack inequality, see \bcite{gil.tru} Theorem 8.18,
from (\ref{conf.lemma.equ.uv-local})
\begin{displaymath}
	\parallel \tilde u - v_k - \lambda \parallel_{L^2(B_{1/8}(x))}
	\leq C(\Lambda) \inf \limits_{B_{1/8}(x)}
	(\lambda - \tilde u + v_k) \leq
\end{displaymath}
\begin{displaymath}
	\leq C(\Lambda) (\max \limits_\sur \tilde u - \tilde u(x) + 2)
	\leq C(\Lambda,p),
\end{displaymath}
and
\begin{displaymath}
	\parallel \max_\sur \tilde u - \tilde u
	\parallel_{L^2(B_{1/8}(x))}
	\leq C(\Lambda,p).
\end{displaymath}
We see from (\ref{conf.lemma.metric})
and (\ref{conf.lemma.grad})
\begin{displaymath}
	c_0 |\max_\sur \tilde u - \tilde \lambda|
	\leq \parallel \max_\sur \tilde u - \tilde \lambda
	\parallel_{L^2(B_{1/8}(x))} \leq
\end{displaymath}
\begin{displaymath}
	\leq C(\Lambda) \parallel \tilde u - \tilde \lambda
	\parallel_{L^2(\sur,g_0)}
	+ \parallel \max_\sur \tilde u - \tilde u
	\parallel_{L^2(B_{1/8}(x))}
	\leq C(\sur,g_0,\Lambda,p) (1 + \sqrt{osc_\sur \tilde u}).
\end{displaymath}
Using (\ref{conf.lemma.meas}),
we have proved
\begin{displaymath}
	\max_\sur \tilde u \geq C(\Lambda,p)
	\Longrightarrow
\end{displaymath}
\begin{equation} \label{conf.lemma.max}
	\parallel \max_\sur \tilde u - \tilde u \parallel_{L^2(\sur,g_0)}
	\leq C(\sur,g_0,\Lambda,p) (1 + \sqrt{osc_\sur \tilde u}).
\end{equation}
Now if $\ \max_\sur \tilde u \gg C(\sur,g_0,\Lambda,p)
(1 + \sqrt{osc_\sur \tilde u})\ $,
we need a lower bound on $\ \tilde u\ $.
If $\ \min_\sur \tilde u \leq - C(\Lambda,p)\ $,
we see from (\ref{conf.lemma.min}) and (\ref{conf.lemma.max}) that
\begin{displaymath}
	\min \limits_\sur \tilde u
	\geq \max \limits_\sur \tilde u
	- C(\sur,g_0,\Lambda,p)
	(1 + \sqrt{osc_\sur \tilde u}) \geq 0,
\end{displaymath}
therefore $\ \min_\sur \tilde u \geq - C(\Lambda,p)\ $.
We put $\ A := [ \max_\sur \tilde u / 2 \leq \tilde u
\leq \max_\sur \tilde u]\ $
and see $\ \max_\sur \tilde u - \tilde u
\geq \max_\sur \tilde u / 2 \mbox{ on } \sur - A\ $,
hence using (\ref{conf.lemma.grad})
\begin{displaymath}
	\frac{1}{2} \max_\sur \tilde u\ \mu_\tgp(\sur - A)
	\leq \int \limits_{\sur - A}
	|\max \limits_\sur \tilde u - \tilde u|
	e^{-2 \tilde u} \d \mu_{\tilde g} \leq
\end{displaymath}
\begin{displaymath}
	\leq C(\sur,g_0,\Lambda,p) \int \limits_\sur
	|\max \limits_\sur \tilde u - \tilde u| \d \mu_{g_0}
	\leq C(\sur,g_0,\Lambda,p) (1 + \sqrt{osc_\sur \tilde u})
\end{displaymath}
and
\begin{displaymath}
	\mu_\tgp(\sur - A) \leq \frac{C(\sur,g_0,\Lambda,p)
	(1 + \sqrt{osc_\sur \tilde u})}{\max \limits_\sur \tilde u}
	\leq \frac{1}{2},
\end{displaymath}
if $\ \max_\sur \tilde u \gg C(\sur,g_0,\Lambda,p)
(1 + \sqrt{osc_\sur \tilde u})\ $ is large enough.
This yields $\ \mu_\tgp(A) \geq 1/2\ $
and by (\ref{conf.lemma.meas})
\begin{displaymath}
	C(\Lambda) \geq \mu_{\tilde g}(\sur)
	\geq \int \limits_A e^{2 \tilde u} \d \mu_\tgp
	\geq \mu_\tgp(A) \exp(\max \limits_\sur \tilde u)
	\geq  \exp(\max \limits_\sur \tilde u)/2,
\end{displaymath}
and we conclude
\begin{equation} \label{conf.lemma.max2}
	\max_\sur \tilde u \leq C(\sur,g_0,\Lambda,p)
	(1 + \sqrt{osc_\sur \tilde u}).
\end{equation}
(\ref{conf.lemma.min}) and (\ref{conf.lemma.max}) yield
\begin{displaymath}
	osc_\sur \tilde u
	= \max \limits_\sur \tilde u - \min \limits_\sur \tilde u \leq
\end{displaymath}
\begin{displaymath}
	\leq  C(\sur,g_0,\Lambda,p) (1 + \sqrt{osc_\sur \tilde u})
	\leq  C(\sur,g_0,\Lambda,p) + \frac{1}{2} osc_\sur \tilde u,
\end{displaymath}
hence $\ osc_\sur \tilde u \leq C(\sur,g_0,\Lambda,p)\ $.
Then (\ref{conf.lemma.tu}) follows from
(\ref{conf.lemma.min}), (\ref{conf.lemma.max})
and (\ref{conf.lemma.grad}),
and the lemma is proved.
\proof
Here, we use this lemma
to get a bound on the conformal factor
for sequences strongly converging in $\ W^{2,2}\ $.

\begin{proposition} \label{conf.strong}

Let $\ f: \sur \rightarrow \rel^n\ $ be a weak local bilipschitz immersion
approximated by smooth immersions $\ f_m\ $
with pull-back metrics $\ g = f^* \geu = e^{2u} \gp,
g_m = f_m^* \geu = e^{2u_m} \gp_m\ $
for some smooth unit volume constant
curvature metrics $\ \gp, \gpm\ $
and satisfying
\begin{displaymath}
\begin{array}{c}
	f_m \rightarrow f
	\quad \mbox{strongly in } W^{2,2}(\sur),
	\mbox{weakly}^* \mbox{ in } W^{1,\infty}(\sur), \\
	\Lambda^{-1} \gp \leq g_m \leq \Lambda \gp.
\end{array}
\end{displaymath}
Then
\begin{displaymath}
	\sup_{m \in \nat}
	\Big( \parallel u_m \parallel_{L^\infty(\sur)},
	\parallel \nabla u_m \parallel_{L^2(\sur,g_m)} \Big)
	< \infty.
\end{displaymath}
\end{proposition}
{\pr Proof:} \\
We want to apply Lemma \ref{conf.lemma}
to $\ f = f_1, \tilde f = f_m, g_0 = \gp\ $.
(\ref{conf.lemma.ub}) and (\ref{conf.lemma.metric2}) are immediate
by the above assumptions
for appropriate possibly larger $\ \Lambda < \infty\ $.

In local charts, we have
\begin{displaymath}
\begin{array}{c}
	g_m \rightarrow g
	\quad \mbox{strongly in } W^{1,2},
	\mbox{weakly}^* \mbox{ in } L^\infty(\sur), \\
	\Gamma_{g_m,ij}^k  \rightarrow \Gamma_{g,ij}^k
	\quad \mbox{strongly in } L^2,
\end{array}
\end{displaymath}
and
\begin{displaymath}
	A_{f_m,ij} = \nabla_i^{g_m} \nabla_j^{g_m} f_m
	\rightarrow \partial_{ij} f - \Gamma_{g,ij}^k \partial_k f =
\end{displaymath}
\begin{displaymath}
	= \nabla_i^{g} \nabla_j^{g} f = A_{f,ij}
	\quad \mbox{strongly in } L^2.
\end{displaymath}
Therefore $\ |A_{f_m}|^2_{g_m} \sqrt{g_m}
\rightarrow |A_f|^2_g \sqrt{g}
\mbox{ strongly in } L^1\ $,
hence for each $\ x \in \sur\ $
there exists a neighbourhood $\ U(x) \mbox{ of } x\ $ with
\begin{displaymath}
	\int \limits_{U(x)} |A_{f_m}|^2_{g_m} \d \mu_{g_m}
	\leq \varepsilon_0(n)
	\quad \mbox{for all } m \in \nat,
\end{displaymath}
where $\ \varepsilon_0(n)\ $ is as in Lemma \ref{conf.lemma}.
Choosing $\ U(x)\ $ even smaller, we may assume
that there are charts $\ \varphi_x: U(x) \stackrel{\approx}{\longrightarrow}
B_1(0) \mbox{ with } \varphi_x(x) = 0
\mbox{ and } c_{0,x} \geu \leq
(\varphi^{-1})^* \gp \leq C_x \geu\ $.
Selecting a finite cover $\ \sur
= \cup_{k=1}^M \varphi_{x_k}^{-1}(B_{1/2}(0))\ $,
we obtain (\ref{conf.lemma.metric}) and (\ref{conf.lemma.fund})
for appropriate $\ \Lambda < \infty\ $.
As clearly $\ supp(f_1 - f_m) \subseteq \sur\ $,
the assertion follows from Lemma \ref{conf.lemma}.
\proof


\def	\xxx	{{ \lambda }}
\def	\yyy	{{ T }}

\setcounter{equation}{0}
\section{Analysis} \label{ana}

\begin{proposition} \label{ana.prop}

Let $\ \Phi = (\Phi_0,\varphi): B_{\lambda_0}^{\M+2}(0) \rightarrow
\rel^{\M+1} = \rel^\M \times \rel, \M \in \nat,\ $
be twice differentiable satisfying
for $\ \xi = (\lambda,\mu,\nu) \in B_{\lambda_0}^{\M+2}(0)
\subseteq \rel^\M \times \rel \times \rel\ $
\begin{displaymath}
\begin{array}{c}
	\parallel \partial_\lambda \Phi_0 - I_\M \parallel \leq 1/2, \\
	| \partial_{(\mu,\nu)} \Phi_0 |,
	| \partial_\lambda \varphi | \leq \varepsilon, \\
	\parallel D^2 \Phi \parallel \leq \Lambda, \\
	\partial_{\mu \mu} \varphi, -\partial_{\nu \nu} \varphi
	\geq \pos
\end{array}
\end{displaymath}
with $\ 0 < \varepsilon, \pos, \lambda_0 \leq 1/4,
1 \leq \Lambda < \infty,\ $
\begin{equation} \label{ana.prop.eps}
	C \Lambda \varepsilon \leq \pos.
\end{equation}
Then for $\ \eta = (\eta_0,\bar \eta) \in \rel^\M \times \rel\ $ with
\begin{displaymath}
\begin{array}{c}
	|\Phi_0(0) - \eta_0| \leq \min(\Lambda \lambda_0^2,
	\lambda_0/8), \\
	|\varphi(0) - \bar \eta| \leq \pos \lambda_0^2/32,
\end{array}
\end{displaymath}
there exists $\ \xi \in B_{\lambda_0}^{\M+2}(0)
\mbox{ with } \mu \nu = 0\ $
and satisfying
\begin{displaymath}
	\Phi(\xi) = \eta
\end{displaymath}
with
\begin{equation} \label{ana.prop.err}
	|\xi| \leq C \gamma^{-1/2} |\Phi(0) - \eta|^{1/2}.
\end{equation}
If further
\begin{equation} \label{ana.prop.ass2}
\begin{array}{c}
	|D \varphi(0)| \leq \sigma, \\
	C \Lambda \varepsilon + C \sigma \lambda_0^{-1}
	+ C \Lambda \lambda_0 \leq \pos,
\end{array}
\end{equation}
there exists a solution $\ \tilde \xi \in B_{\lambda_0}^{\M+2}(0)
\mbox{ of } \Phi(\tilde \xi) = \eta,
\tilde \mu \tilde \nu = 0\ $ with
\begin{displaymath}
\begin{array}{c}
	\mu \tilde \mu \leq 0, \tilde \nu = 0,
	\quad \mbox{if } \nu = 0, \\
	\nu \tilde \nu \leq 0, \tilde \mu = 0,
	\quad \mbox{if } \mu = 0.
\end{array}
\end{displaymath}
\end{proposition}
{\pr Proof:} \\
After the choice in (\ref{ana.prop.ass}) below,
we will need only one variable $\ \mu \mbox{ or } \nu\ $.
Therefore to simplify the notation,
we put $\ \nu = 0\ $ and omit $\ \nu\ $.

First there exists a twice differentiable function
$\ \xxx: ]-\lambda_0/2,\lambda_0/2[ \rightarrow B_{\lambda_0/2}^\M(0)\ $
\begin{equation} \label{ana.prop.solve}
	\Phi_0(\xxx(\mu),\mu) = \eta_0.
\end{equation}
Indeed putting
\begin{displaymath}
	\yyy_\mu(\lambda) := \lambda - \Phi_0(\lambda,\mu) + \eta_0
	\quad \mbox{ for } |\lambda| \leq \lambda_0/2,
\end{displaymath}
we see
\begin{displaymath}
	\parallel \yyy'_\mu \parallel
	\leq \parallel I_\M - \partial_\lambda \Phi_0 \parallel
	\leq 1/2,
\end{displaymath}
hence
\begin{displaymath}
	|\yyy_\mu(\lambda_1) - \yyy_\mu(\lambda_2)|
	\leq \frac{1}{2} |\lambda_1 - \lambda_2|.
\end{displaymath}
As
\begin{displaymath}
	|\yyy_\mu(0)| = |\Phi_0(0,\mu) - \eta_0|
	\leq |\Phi_0(0,\mu) - \Phi_0(0,0)| + |\Phi_0(0) - \eta_0| \leq
\end{displaymath}
\begin{displaymath}
	\leq \parallel \partial_\mu \Phi_0 \parallel |\mu| + \lambda_0/8
	< (\varepsilon/2 + 1/8) \lambda_0 \leq \lambda_0/4,
\end{displaymath}
we get
\begin{displaymath}
	|\yyy_\mu(\lambda)| \leq |\lambda|/2 + |\yyy_\mu(0)|
	< \lambda_0/2
\end{displaymath}
and by Banach's fixed point theorem,
(\ref{ana.prop.solve}) has a unique solution
$\ \lambda = \xxx(\mu) \in B_{\lambda_0/2}^\M(0)\ $.
In particular
\begin{equation} \label{ana.prop.first}
	|\xxx(0)| \leq 2 |\yyy_0(0)|
	\leq 2 |\Phi_0(0) - \eta_0| \leq 2 \Lambda \lambda_0^2.
\end{equation}
By implict function theorem, $\ \xxx\ $ is twice differentiable and
\begin{displaymath}
\begin{array}{c}
	\partial_\lambda \Phi_0 \partial_\mu \xxx + \partial_\mu \Phi_0 = 0, \\
	\partial_\mu \xxx^T \partial_{\lambda \lambda} \Phi_0
	\partial_\mu \xxx
	+ 2 \partial_{\lambda \mu} \Phi_0 \partial_\mu \xxx
	+ \partial_\lambda \Phi_0 \partial_{\mu \mu} \xxx
	+ \partial_{\mu \mu} \Phi_0 = 0.
\end{array}
\end{displaymath}
This yields
\begin{equation} \label{ana.prop.varphi}
\begin{array}{c}
	|\partial_\mu \xxx| \leq 2 |\partial_\mu \Phi_0|
	\leq 2 \varepsilon \leq 1/2, \\
	|\partial_{\mu \mu} \xxx|
	\leq 2 |\partial_\mu \xxx^T \partial_{\lambda \lambda} \Phi_0
	\partial_\mu \xxx
	+ 2 \partial_{\lambda \mu} \Phi_0 \partial_\mu \xxx
	+ \partial_{\mu \mu} \Phi_0| \leq C \Lambda.
\end{array}
\end{equation}
Next we put
\begin{displaymath}
	\psi(\mu) := \varphi(\xxx(\mu),\mu) - \bar \eta
	\quad \mbox{for } |\mu| < \lambda_0/2
\end{displaymath}
and see by (\ref{ana.prop.first})
\begin{displaymath}
	|\psi(0)| \leq |\varphi(\xxx(0),0) - \varphi(0)| + |\varphi(0) - \bar \eta| \leq
\end{displaymath}
\begin{equation} \label{ana.prop.psi0}
	\leq \parallel \partial_\lambda \varphi \parallel |\xxx(0)|
	+ \pos \lambda_0^2/32
	\leq C \Lambda \varepsilon \lambda_0^2 + \pos \lambda_0^2/32
	\leq \pos \lambda_0^2 / 16.
\end{equation}
Clearly $\ \psi\ $ is twice differentiable and
\begin{equation} \label{ana.prop.abl-psi}
\begin{array}{c}
	\partial_\mu \psi
	= \partial_\lambda \varphi \partial_\mu \xxx + \partial_\mu \varphi, \\
	\partial_{\mu \mu} \psi
	= \partial_\mu \xxx^T \partial_{\lambda \lambda} \varphi
	\partial_\mu \xxx
	+ 2 \partial_{\lambda \mu} \varphi \partial_\mu \xxx
	+ \partial_\lambda \varphi \partial_{\mu \mu} \xxx
	+ \partial_{\mu \mu} \varphi.
\end{array}
\end{equation}
This yields using (\ref{ana.prop.varphi})
\begin{equation} \label{ana.prop.psi-abl}
\begin{array}{c}
	|\partial_{\mu \mu} \psi - \partial_{\mu \mu} \varphi|
	\leq |\partial_\mu \xxx^T \partial_{\lambda \lambda} \varphi
	\partial_\mu \xxx
	+ 2 \partial_{\lambda \mu} \varphi \partial_\mu \xxx
	+ \partial_\lambda \varphi \partial_{\mu \mu} \xxx| \leq \\
	\leq C \parallel D^2 \varphi \parallel |\partial_\mu \xxx|
	+ \parallel \partial_\lambda \varphi \parallel |\partial_{\mu \mu} \xxx|
	\leq C \Lambda \varepsilon.
\end{array}
\end{equation}
Since the assumptions and conclusions
are the equivalent for $\ -\Phi \mbox{ and } - \eta\ $,
we assume after possibly exchanging
$\ \mu \mbox{ by } \nu\ $ that
\begin{equation} \label{ana.prop.ass}
	\psi(0) \leq 0.
\end{equation}
Replacing $\ \mu \mbox{ by } - \mu\ $,
we may further assume $\ \partial_\mu \psi(0) \geq 0\ $.
On the other hand by (\ref{ana.prop.psi0}) and (\ref{ana.prop.psi-abl})
\begin{displaymath}
	\liminf \limits_{\mu \linksarrow \lambda_0/2} \psi(\mu)
	\geq \psi(0) + \partial_\mu \psi(0) \lambda_0/2
	+ (\inf \partial_{\mu \mu} \psi/2) (\lambda_0/2)^2 \geq
\end{displaymath}
\begin{displaymath}
	\geq - \pos \lambda_0^2/16
	+ (\pos - C \Lambda \varepsilon) \lambda_0^2/8 > 0,
\end{displaymath}
when using (\ref{ana.prop.eps}).
Therefore there exists $\ 0 \leq \mu < \lambda_0/2
\mbox{ with } \psi(\mu) = 0\ $,
hence putting $\ \xi := (\xxx(\mu),\mu,0) \in B_{\lambda_0}^{\M+2}(0)\ $
\begin{displaymath}
	\Phi(\xi)
	= \Big( \Phi_0(\xxx(\mu),\mu) , \varphi(\xxx(\mu),\mu) \Big)
	= \Big( \eta_0 , \psi(\mu) + \bar \eta \Big)
	= \eta.
\end{displaymath}
Choosing $\ \lambda_0\ $ small enough,
we further obtain (\ref{ana.prop.err}).

If further (\ref{ana.prop.ass2}) is satisfied for the original $\ \lambda_0\ $,
we see from (\ref{ana.prop.first}), (\ref{ana.prop.varphi})
and (\ref{ana.prop.abl-psi}) that
\begin{displaymath}
	|\partial_\mu \psi(0)|
	\leq C |D \varphi(\xxx(0),0)|
	\leq C |D \varphi(0)| + C \parallel D^2 \varphi \parallel |\xxx(0)|
	\leq C \sigma + C \Lambda \lambda_0^2.
\end{displaymath}
Proceeding as above with (\ref{ana.prop.psi0}) and (\ref{ana.prop.psi-abl}),
we calculate
\begin{displaymath}
	\liminf \limits_{\mu \rechtsarrow - \lambda_0/2} \psi(\mu)
	\geq \psi(0) + \partial_\mu \psi(0) \lambda_0/2
	+ (\inf \partial_{\mu \mu} \psi/2) (\lambda_0/2)^2 \geq
\end{displaymath}
\begin{displaymath}
	\geq - \pos \lambda_0^2/16
	-  (C \sigma + C \Lambda \lambda_0^2) \lambda_0
	+ (\pos - C \Lambda \varepsilon) \lambda_0^2/8 =
\end{displaymath}
\begin{displaymath}
	\geq (\pos/16 -  C \Lambda \varepsilon
	- C \sigma \lambda_0^{-1} - C \Lambda \lambda_0) \lambda_0^2 > 0,
\end{displaymath}
when using (\ref{ana.prop.ass2}).
Therefore there exists $\ - \lambda_0/2 < \tilde \mu \leq 0
\mbox{ with } \psi(\tilde \mu) = 0\ $,
hence putting $\ \tilde xi := (\xxx(\tilde \mu),\tilde \mu,0)
\in B_{\lambda_0}^{\M+2}(0)\ $,
we get $\ \Phi(\tilde \xi) = \eta \mbox{ and } \mu \tilde \mu \leq 0\ $.
\proof


\end{appendix}



\end{document}